\documentclass[preprint,12pt]{elsarticle}

\usepackage{amssymb}
\usepackage{amsmath}
\numberwithin{equation}{section}

\usepackage{color}
\usepackage{epstopdf}

\usepackage{amsmath}
\usepackage{dsfont}

\newtheorem{theorem}{Theorem}[section]

\newtheorem{corollary}[theorem]{Corollary}

\newtheorem{lemma}[theorem]{Lemma}

\journal{}

\begin{document}

\begin{frontmatter}




\title{The discontinuous Galerkin method for the Oseen eigenvalue problem} 

 \author{Lingling Sun\fnref{label1}}
 \ead{sunlingling@gmc.edu.cn}

 \author{Shixi Wang\fnref{label2}}
 \ead{wangshixi@gznu.edu.cn}

 \author{Hai Bi\fnref{label2}}
 \ead{bihaimath@gznu.edu.cn}

 \author{Yidu Yang\corref{cor1}\fnref{label2}}
 \ead{ydyang@gznu.edu.cn}
 \cortext[cor1]{Corresponding author}

\address[label1]{School of Biology and Engineering(School of Modern Industry for Health and Medicine), Guizhou Medical University, Guiyang, Guizhou 561113, China}
\address[label2]{School of Mathematical Sciences, Guizhou Normal University, Guiyang, Guizhou 550025, China}

\begin{abstract}
In this paper,  we focus on investigating symmetric and nonsymmetric discontinuous Galerkin (DG) methods
for solving the Oseen eigenvalue problem based on the
velocity-pressure formulation in $\mathbb{R}^{d}(d=2,3)$.
We derive the
a priori and a posteriori error estimates for the approximate eigenpairs
for each method.
We establish an adjoint-consistent symmetric DG method and derive optimal a priori error estimates,
 and prove the reliability and effectiveness of the error estimators for approximate eigenfunctions, as well as the reliability of the estimator for approximate eigenvalues.
Numerical experiments confirm our theoretical analysis and demonstrate that the symmetric DG method achieves the optimal order of convergence,
and that the nonsymmetric DG methods produce fewer spurious eigenvalues than the symmetric DG method for a fixed small penalty parameter $\gamma$.
\end{abstract}


\end{frontmatter}
\section{Introduction}

\indent The Oseen eigenvalue problem, as a linearization of the Navier-Stokes equation, has important applications in fluid mechanics.
Some references, such as \cite{Benzi2007}, discuss the computation of the spectrum of the Oseen operators, which plays a role in designing
appropriate solvers for the Navier-Stokes equations. Additionally, for preconditioners of the Navier-Stokes equations, the computation of
the spectrum of the Oseen problem is relevant, as shown in the work of \cite{Fan2016}. Currently, the study of numerical methods for the Oseen
eigenvalue problem has attracted significant attention in the academic community.
For the Oseen eigenvalue problem based on the velocity-pressure formulation, Lepe et al. \cite{Lepe2024}
derived the a priori and a posteriori error estimates using the Mini-element family and the Taylor-Hood
$(P_{2}-P_{1})$ element on two-dimensional and three-dimensional domains,
and Adak et al. \cite{Adak2025} developed the a
priori error estimates using a nonconforming virtual element discretization on two-dimensional domains.
For the Oseen eigenvalue problem based on the velocity-pseudostress formulation,
Lepe et al. \cite{Lepe2024new} employed the Raviart-Thomas (RT) and Brezzi-Douglas-Marini (BDM)
tensorized elements and proved the a priori error estimates on both two and three dimensional domains.
However, to the best of our knowledge, no reports exist on the discontinuous Galerkin (DG) method for the Oseen eigenvalue problem.\\
\indent The DG method,  first proposed by Reed and Hill
 \cite{ReedandHill1973} in 1973, has evolved into a widely used approach for solving differential
 equations (see, e.g., introductory textbooks \cite{Antonio2012,Riviere2008,Hesthaven2008,Cockburn1999}).
Among them, Chapter 6 in \cite{Antonio2012}, Chapter 6 in \cite{Riviere2008}, Chapter 7 in
\cite{Hesthaven2008}, and Chapter 5 in \cite{Cockburn1999} elaborate on the DG method for the Stokes problem.
The articles \cite{Hansbo2002,Stokes-flow2010,Navier-Stokes2010,Badia2014,Pietro2010,Houston2005,Schotzau,
Kanschat2008,Cockburn2002} have investigated the DG method for the Stokes problem.
In 2006, Antonietti et al. \cite{Antonietti2006} first applied the DG method to solve the Laplace eigenvalue problem. This paper is key for IPDG methods for eigenvalue problems. In \cite{Antonietti2006}, the authors employed the theory developed by Descloux, Nassif and Rappaz in \cite{Descloux1978a,Descloux1978b} for their analysis.
Since then, many scholars have employed the DG method for the Stokes eigenvalue problem (see, e.g., \cite{Cliff2010,Sun2023Cicp,DG+Gedicke2020,DG+Lepe2020,Lepe2023+DG+Stokes}). Among these works, Lepe et al. \cite{DG+Lepe2020,Lepe2023+DG+Stokes} adopted the analytical approach of \cite{Antonietti2006}.
In this paper, we focus on investigating the DG $P_{k}-P_{k-1}~(k\geq 1)$ element method,
including the classic symmetric interior penalty ({\bf SIP}) method \cite{Arnold1982+SIP}, the nonsymmetric interior penalty ({\bf NIP}) method
\cite{Baumann+NIP,Riviere+NIP} and the incomplete interior penalty ({\bf IIP}) method \cite{Dawson+IIP},
for solving the Oseen eigenvalue problem based on the
velocity-pressure formulation in $\mathbb{R}^{d}~(d=2,3)$.
Inspired by the works of Boffi on non-conforming methods (see Section 11 in \cite{Boffi2010}), we employ the Babu\v{s}ka-Osborn theory for the analysis in this paper.
\\
\indent The Oseen eigenvalue problem is non-self-adjoint, and its error analysis and numerical computation are closely related to the DG formulations of both
the primal and adjoint eigenvalue problems. The adjoint consistency is crucial to derive the optimal error estimates for the DG approximate solutions
(see Section 6 in [32] and Section 3 in [33]). We discuss adjoint consistency for each method and, in particular, prove the adjoint consistency of the {\bf SIP} method.\\
\indent Based on \cite{Nedelec1986,Acosta2006}, we prove the discrete inf-sup condition (\ref{s2.35}) and, via a subtle analysis, establish the discrete coercivity property (\ref{s2.19}).
We then prove that the discrete solution operator converges to the continuous one. Consequently, applying the Babu\v{s}ka-Osborn spectral approximation
 theory \cite{Babuska1991book}, we derive error estimates for each method. In particular, we obtain optimal a priori error estimates for the approximate
  eigenpairs associated with the {\bf SIP} method. However, due to the lack of adjoint consistency, we cannot derive optimal a priori error estimates for the {\bf NIP} and {\bf IIP} methods.\\
\indent In practical finite element computations, it is desirable to perform computations in an adaptive fashion (see, e.g., \cite{Babuska1978,Verfurth2013,oden2011,Adaptive2023,Adaptive2024}
and the references therein). For the Oseen eigenvalue problem,
Lepe et al. \cite{Lepe2024} investigated the a posteriori error estimation and adaptive computation for low-order Mini-element and low-order Taylor-Hood elements. Although high-order Taylor-Hood
elements have been established by Boffi et al. (see Section 8.8.2 in \cite{Boffi2013}) and applied to the Stokes eigenvalue problems \cite{Boffi2025},
but there are no reports on  a posteriori error estimation and adaptive method for the Oseen eigenvalue problem.
In this paper,
we propose an adaptive DG method using $P_{k}-P_{k-1}~(k\geq 1)$ elements.
By the inf-sup condition (\ref{s2.7}), we establish the a posteriori error formula (Theorem \ref{th3-1}),
and propose the a posteriori error estimators of residual type.
By employing the enriching operator (see \cite{Brenner2003,Ohannes2003}) and the bubble function technique (see \cite{Verfurth2013}),
 we derive the a posteriori error estimates for each method. For the {\bf SIP} method, we prove the reliability and efficiency of the estimator
for approximate eigenfunctions, as well as the reliability of the estimator for approximate eigenvalues.\\
\indent To validate our approach, we conduct numerical computations on both uniform and adaptively refined meshes. The numerical results demonstrate that
the {\bf SIP} method is capable of achieving the optimal order of convergence and capable of yielding high-accuracy approximate eigenvalues.
Nevertheless, the numerical results demonstrate that although the {\bf NIP} and {\bf IIP} methods cannot achieve the optimal convergence order,
these two methods yield fewer spurious eigenvalues than the {\bf SIP} method for a fixed small penalty parameter $\gamma$.\\
\indent Throughout this article, $C$ in different positions represents different positive
constants, which are independent of mesh size $h$. We use $a\lesssim b$ to represent $a\leq C b$, and $a\backsimeq b$
to represent $a\lesssim b$ and $b\lesssim a$.

\section{ The DG approximation and its a priori error estimate}

\indent Consider the following Oseen eigenvalue problem:
\begin{equation}\label{s2.1}
\begin{cases}
-\mu\Delta \mathbf{u}+(\beta\cdot\nabla)\mathbf{u}+\nabla p&=\lambda \mathbf{u},~~~~ in~~\Omega,\\
~~~~~~~~~~~~~~~~~~~~~~~~~~div \mathbf{u}&=0,~~~~~~in ~~\Omega,\\
~~~~~~~~~~~~~~~~~~~~~~~~~~~\int_{\Omega}p&=0,~~~~~~in ~~\Omega,\\
~~~~~~~~~~~~~~~~~~~~~~~~~~~~\mathbf{u}&=0,~~~~~~on~~\partial\Omega,
\end{cases}
\end{equation}
where $\Omega\subset \mathbb{R}^{d}(d=2,3)$ is a bounded polyhedral domain,
$\mathbf{u}=(u_{1},...,u_{d})^{t}$ is the velocity of the flow, $p$ is the pressure,
$\boldsymbol{\beta}$ is a given  convection field, representing a steady flow velocity such that
 $\boldsymbol{\beta}\in W^{1,\infty}({\Omega})^{d}$ is divergence free, $\mu>0$ is the kinematic viscosity parameter of the fluid.\\
\indent The standard assumptions on the coefficient are the following (see \cite{John2016}):\\
\indent $\cdot$ $\|\boldsymbol{\beta}\|_{\infty,\Omega}\sim 1$ if $\mu\leq \|\boldsymbol{\beta}\|_{\infty,\Omega}$,\\
\indent $\cdot$ $\mu\sim 1$ if $\|\boldsymbol{\beta}\|_{\infty,\Omega}< \mu$,\\
where the first situation is the interesting one which appears in applications.\\
\indent  Let $H^{\rho}(D)$ be the complex Sobolev space of order $\rho\geq 0$ on $D\subseteq \Omega$
with the norm $\|\cdot\|_{\rho,D}$.
We denote $\|\mathbf{v}\|_{\rho,D}=\sum\limits_{i=1}^{d}\|v_{i}\|_{\rho,D}$ for $\mathbf{v}=(v_{1},\cdots,v_{d})\in H^{\rho}(D)^{d}$.
If  $D=\Omega$, we omit the subscript $D$ from $\|\cdot\|_{\rho,D}$.
We denote by $(\cdot,\cdot)$ the inner product in $L^{2}(\Omega)^{d}$ which is given by $(\mathbf{u},\mathbf{v})=\int_{\Omega}\mathbf{u}\cdot\overline{\mathbf{v}} dx$
and  $(\nabla\mathbf{u},\nabla\mathbf{v})=\int_{\Omega}\nabla\mathbf{u}:\overline{\nabla\mathbf{v}} dx$.
We denote $H_{0}^{1}(\Omega)=\{v\in H^{1}(\Omega), v|_{\partial\Omega}=0\}$,
\begin{align*}
\mathbb{X}=H_{0}^{1}(\Omega)^{d},~~~
\mathbb{W} = L_{0}^{2}(\Omega)=\{\varrho\in L^{2}(\Omega):(\varrho,1)=0\}.
\end{align*}
\indent The weak formulation of (\ref{s2.1}) is given by:
 find $(\lambda,\mathbf{u},p)\in \mathbb{C}\times \mathbb{X} \times \mathbb{W}$, $\|\mathbf{u}\|_{0}=1$, such that
\begin{align}
A(\mathbf{u},\mathbf{v})+B(\mathbf{v},p)&=\lambda (\mathbf{u},\mathbf{v}),~~~\forall \mathbf{v}\in \mathbb{X},\label{s2.2}\\
\overline{B(\mathbf{u},q)}&=0, ~~~\forall q\in \mathbb{W},\label{s2.3}
\end{align}
where
\begin{align*}
A(\mathbf{u},\mathbf{v})=\mu(\nabla\mathbf{u}, {\nabla\mathbf{v}})+((\boldsymbol{\beta}\cdot\nabla)\mathbf{u},{\mathbf{v}}),~~~
B(\mathbf{v},p)=-(p, \nabla\cdot\mathbf{v}).
\end{align*}
\indent Thanks to Remark 5.6 in \cite{John2016},  we deduce that $Re((\boldsymbol{\beta}\cdot\nabla)\mathbf{v},\mathbf{v})=0$. Therefore,
we can obtain
\begin{align}\label{s2.4}
Re A(\mathbf{v},\mathbf{v})\gtrsim\|\mathbf{v}\|^{2}_{1},
~~~ \forall \mathbf{v}\in \mathbb{X}.
\end{align}
\indent Thanks to Remark 53.10 in \cite{Ern2021} we know that the sesquilinear form $\mathbb{B}(\cdot,\cdot)$ satisfies the following inf-sup condition: there is a constant $\alpha>0$ such that
\begin{align}\label{s2.5}
\inf\limits_{q\in \mathbb{W}}\sup\limits_{\mathbf{v}\in \mathbb{X}}\frac{|B(\mathbf{v},q)|}{\|q\|_{0}\|\mathbf{v}\|_{1}}\geq \alpha.
\end{align}

\indent Let $\mathbb{V}=\mathbb{X}\times \mathbb{W}$ with norm $\|[\mathbf{v},q]\|_{\mathbb{V}}=\|\mathbf{v}\|_{1}+\|q\|_{0}$. Referring to  (49.39) in \cite{Ern2021}, we introduce the following sesquilinear form
\begin{align}\label{s2.6}
\mathcal{A}([\mathbf{w},r], [\mathbf{v},q])=A(\mathbf{w},\mathbf{v})+B(\mathbf{v},r)
+\overline{B(\mathbf{w},q)},~~~ \forall [\mathbf{w},r], [\mathbf{v},q]\in \mathbb{V}.
\end{align}
And from Theorem 49.15 in \cite{Ern2021}, we have
\begin{align}\label{s2.7}
\inf\limits_{[\mathbf{w},r]\in \mathbb{V}}\sup\limits_{[\mathbf{v},q]\in \mathbb{V}}\frac{|\mathcal{A}([\mathbf{w},r], [\mathbf{v},q])|}{\|[\mathbf{w},r]\|_{\mathbb{V}}\|[\mathbf{v},q]\|_{\mathbb{V}}}\geq \alpha,~~~
\inf\limits_{[\mathbf{v},q]\in \mathbb{V}}\sup\limits_{[\mathbf{w},r]\in \mathbb{V}}\frac{|\mathcal{A}([\mathbf{w},r], [\mathbf{v},q])|}{\|[\mathbf{w},r]\|_{\mathbb{V}}\|[\mathbf{v},q]\|_{\mathbb{V}}}\geq \alpha.
\end{align}
The formulation (\ref{s2.7}) plays a crucial role in the proof of Theorem \ref{th3-1}.\\
 \indent Let $\mathcal{T}_{h}=\{\tau\}$ be a regular simplex partition of $\Omega$ with the mesh diameter $h=\max\limits_{\tau\in \mathcal{T}_{h}} h_{\tau}$ where $h_{\tau}$ is the diameter of element $\tau$.
We use $\mathcal{E}_{h}^{i}$ and $\mathcal{E}_{h}^{b}$ to denote the set of interior faces (edges) and the set of faces (edges) on $\partial\Omega$, respectively.
 $\mathcal{E}_{h}=\mathcal{E}_{h}^{i}\cup\mathcal{E}_{h}^{b}$. We use $h_{F}$ to denote the measure of $F\in \mathcal{E}_{h}$. We denote by $(\cdot, \cdot)_{\tau}$ and $(\cdot, \cdot)_{F}$ the inner product in $L^{2}(\tau)$ and $L^{2}(F)$, respectively.
We use $\omega(\tau)$ to represent the union of all elements, which share at least one edge (face) with $\tau$,
and use $\omega(F)$ to represent the union of the elements having in common with $F$.\\
\indent The broken Sobolev space is defined by:
$$H^{\xi}(\mathcal{T}_{h})=\{v\in L^{2}(\Omega): v|_{\tau}\in H^{\xi}(\tau),~\forall \tau\in\mathcal{T}_{h}\}~~(\xi>\frac{1}{2}).$$
\indent For any $F\in\mathcal{E}_{h}^{i}$, there are two simplices  $\tau^{+}$ and $\tau^{-}$ such that $F=\tau^{+}\cap\tau^{-}$ (e.g. see \cite{Badia2014}). Let
$\mathbf{n}^{+}$ be the unit normal of $F$ pointing from $\tau^{+}$ to $\tau^{-}$ and let $\mathbf{n}^{-}=-\mathbf{n}^{+}$. \\
\indent For any $\varphi\in H^{1}(\mathcal{T}_{h})$, we denote its jump and mean on $F\in\mathcal{E}_{h}^{i}$
by $[\![\varphi]\!]=\varphi^{+}\mathbf{n}^{+}+\varphi^{-}\mathbf{n}^{-}$ and $\mathbf{\{}\varphi\mathbf{\}}=\frac{1}{2}(\varphi^{+}+\varphi^{-})$, respectively, where $\varphi^{\pm}=\varphi|_{\tau^{\pm}}$.
For $\boldsymbol{\varphi}\in H^{1}(\mathcal{T}_{h})^{d}$, we denote by $[\![\boldsymbol{\varphi}]\!]=\boldsymbol{\varphi}^{+}\cdot\mathbf{n}^{+}+\boldsymbol{\varphi}^{-}\cdot\mathbf{n}^{-}$ the jump and $\{\boldsymbol{\varphi}\}=\frac{1}{2}(\boldsymbol{\varphi}^{+}+\boldsymbol{\varphi}^{-})$ the mean of $\boldsymbol{\varphi}$ on $F\in\mathcal{E}_{h}^{i}$.
We also denote by $[\![\underline{\boldsymbol{\varphi}}]\!]=\boldsymbol{\varphi}^{+}\otimes \mathbf{n}^{+}+\boldsymbol{\varphi}^{-}\otimes \mathbf{n}^{-}$ the full jump of $\boldsymbol{\varphi}$ on $F\in\mathcal{E}_{h}^{i}$, where $\boldsymbol{\varphi}\otimes \mathbf{n}=[\varphi_{i}n_{j}]_{1\leq i,j\leq d}$, $\boldsymbol{\varphi}=(\varphi_{i}), \mathbf{n}=(n_{j})$.
For tensors $\boldsymbol{\chi}\in H^{1}(\mathcal{T}_{h})^{d\times d}$, we denote by $[\![\boldsymbol{\chi}]\!]=\boldsymbol{\chi}^{+}\mathbf{n}^{+}+\boldsymbol{\chi}^{-}\mathbf{n}^{-}$ the jump and $\{\boldsymbol{\chi}\}=\frac{1}{2}(\boldsymbol{\chi}^{+}+\boldsymbol{\chi}^{-})$ the mean on $F\in\mathcal{E}_{h}^{i}$. When $F\in\mathcal{E}_{h}^{b}$, by modifying the above definitions appropriately, we obtain the jump and the mean on $\partial\Omega$. That is to say, we modify $\varphi^{-}=0$ (similarly, $\boldsymbol{\varphi}^{-}=0$ and $\boldsymbol{\chi}^{-}=0$) to obtain the definition of jump on $\partial\Omega$, and modify $\varphi^{-}=\varphi^{+}$(similarly, $\boldsymbol{\varphi}^{-}=\boldsymbol{\varphi}^{+}$ and $\boldsymbol{\chi}^{-}=\boldsymbol{\chi}^{+}$) to obtain the definition of mean on $\partial\Omega$.\\
\indent The discrete velocity and pressure spaces are defined as follows:
\begin{align*}
&\mathbb{X}_{h}=\{\mathbf{v}_{h}\in L^{2}(\Omega)^{d}: \mathbf{v}_{h}|_{\tau}\in P_{k}(\tau)^{d},~\forall \tau \in\mathcal{T}_{h}\},\\
&\mathbb{W}_{h}=\{q_{h}\in \mathbb{W}:q_{h}|_{\tau}\in P_{k-1}(\tau),~\forall \tau \in\mathcal{T}_{h}\},
\end{align*}
where $P_{k}(\tau)$ is the space of polynomials of degree less than or equal to $k\geq 1$ on $\tau$.
And let $\mathbb{X}(h)=\mathbb{X}+\mathbb{X}_{h}$. $\mathbb{V}_{h}=\mathbb{X}_{h}\times \mathbb{W}_{h}$.
\\
\indent The DG formula for the problem (\ref{s2.1}) is to find $(\lambda_{h}, \mathbf{u}_{h}, p_{h})\in \mathbb{C}\times\mathbb{X}_{h}\times \mathbb{W}_{h}$,
$\|\mathbf{u}_{h}\|_{0}=1$, such that
\begin{align}
A_{h}(\mathbf{u}_{h},\mathbf{v}_{h})+B_{h}(\mathbf{v}_{h},p_{h})&=\lambda_{h}(\mathbf{u}_{h},\mathbf{v}_{h}),  ~~\forall \mathbf{v}_{h}\in\mathbb{X}_{h},\label{s2.8}\\
\overline{B_{h}(\mathbf{u}_{h},q_{h})}&=0,~~~\forall q_{h}\in \mathbb{W}_{h},\label{s2.9}
\end{align}
where
\begin{align}
&A_{h}(\mathbf{u}_{h},\mathbf{v}_{h})=\sum\limits_{\tau\in\mathcal{T}_{h}}\int_{\tau}\mu\nabla\mathbf{u}_{h}:\overline{\nabla\mathbf{v}}_{h}dx
+ \sum\limits_{\tau\in\mathcal{T}_{h}}\int_{\tau}(\boldsymbol{\beta}\cdot\nabla)\mathbf{u}_{h}\cdot \overline{\mathbf{v}}_{h}dx \nonumber\\
&~~~
-\sum\limits_{F\in\mathcal{E}_{h}}\int_{F}\mu\{\nabla\mathbf{u}_{h}\}:\overline{[\![\underline{{{\mathbf{v}}_{h}}}]\!]}ds
+\theta\sum\limits_{F\in\mathcal{E}_{h}}\int_{F}\mu\{\overline{\nabla\mathbf{v}}_{h}\}:[\![\underline{\mathbf{u}_{h}}]\!]ds\nonumber\\
&~~~+\sum\limits_{F\in\mathcal{E}_{h}}\int_{F}\frac{\gamma}{h_{F}}[\![\underline{\mathbf{u}_{h}}]\!]:
\overline{[\![{\underline{\mathbf{v}}_{h}}]\!]}ds
-\sum\limits_{F\in\mathcal{E}^{i}_{h}}\int_{F} [\![\mathbf{u}_{h}\otimes\boldsymbol{\beta}]\!]\cdot\overline{\{ \mathbf{v}_{h}\}} ds\nonumber\\
&~~~-\frac{1}{2}\sum\limits_{F\in\mathcal{E}_{h}\bigcap\partial\Omega}\int_{F} (\mathbf{u}_{h}\otimes\boldsymbol{\beta})\mathbf{n}\cdot\overline{ \mathbf{v}_{h}} ds,\label{s2.10}\\
&B_{h}(\mathbf{v}_{h},q_{h})=
-\sum\limits_{\tau\in\mathcal{T}_{h}}\int_{\tau}q_{h}\overline{\nabla\cdot\mathbf{v}_{h}}dx
+\sum\limits_{F\in\mathcal{E}_{h}}\int_{F}\{q_{h}\}\overline{[\![\mathbf{v}_{h}]\!]}ds.\label{s2.11}
\end{align}
The term $\frac{\gamma}{h_{F}}$ is the interior penalty parameter.\\
\indent The {\bf SIP} method is obtained for $\theta=-1$, the {\bf NIP} method for $ \theta=1$, and the {\bf IIP} method for $\theta=0$.

\indent Define the DG-norm as follows:
\begin{align}
\|\mathbf{v}_{h}\|_{h}^{2}&=\sum\limits_{\tau\in\mathcal{T}_{h}}\| \nabla\mathbf{v}_{h}\|_{0,\tau}^{2}+
\sum\limits_{F\in\mathcal{E}_{h}}\int_{F}{\frac{\gamma}{h_{F}}}|[\![\underline{\mathbf{v}_{h}}]\!]|^{2}ds,
~on~(H^{1}(\mathcal{T}_{h}))^{d};\label{s2.12}\\
|||\mathbf{v_{h}}|||^{2}&=\|\mathbf{v}_{h}\|_{h}^{2}+\sum\limits_{F\in \mathcal{E}_{h}}\int_{F}h_{F}|\{\nabla \mathbf{v}_{h}\}|^{2}ds,~on~(H^{1+s}(\mathcal{T}_{h}))^{d},(s>\frac{1}{2}).\label{s2.13}
\end{align}
Note that $\|\cdot\|_{h}$ is equivalent to $|||\cdot|||$ on $\mathbb{X}_{h}$.\\
 \indent Using the trace theorem on the reference
element and the scaling argument, we deduce for any
 $\tau\in \mathcal{T}_{h}$ that the trace inequality
\begin{align} \label{s2.14}
\|w\|_{0,\partial\tau}\lesssim h_{\tau}^{-\frac{1}{2}}\|w\|_{0,\tau}+h_{\tau}^{s-\frac{1}{2}}|w|_{s,\tau},~~\forall w\in H^{s}(\tau),~s\in (\frac{1}{2},1].
\end{align}
 From the Cauchy-Schwarz inequality and (\ref{s2.14}), we obtain
 continuity of the sesquilinear forms $A_{h}$ and $B_{h}$:
\begin{align}\label{s2.15}
&~~~|A_{h}(\mathbf{u}_{h},\mathbf{v}_{h})|\lesssim |||\mathbf{u}_{h}|||~|||\mathbf{v}_{h}|||\nonumber\\
&\lesssim (\|\mathbf{u}_{h}\|_{h}+(\sum\limits_{\tau\in \mathcal{T}_{h}}h_{\tau}^{2s}\|\nabla\mathbf{u}_{h}\|_{s,\tau}^{2})^{\frac{1}{2}})
~ (\|\mathbf{v}_{h}\|_{h}+(\sum\limits_{\tau\in \mathcal{T}_{h}}h_{\tau}^{2s}\|\nabla\mathbf{v}_{h}\|_{s,\tau}^{2})^{\frac{1}{2}}),\nonumber\\
&~~~~~~
\forall [\mathbf{u}_{h},\mathbf{v}_{h}]\in (H^{1+s}(\mathcal{T}_{h}))^{d}\times  (H^{1+s}(\mathcal{T}_{h}))^{d},~~(s>\frac{1}{2}).
\end{align}
\begin{align}\label{s2.16}
&|B_{h}(\mathbf{v}_{h},q_{h})|=
|-\sum\limits_{\tau\in\mathcal{T}_{h}}\int_{\tau}q_{h}\overline{\nabla\cdot\mathbf{v}_{h}}dx
+\sum\limits_{F\in\mathcal{E}_{h}}\int_{F}\{q_{h}\}\overline{[\![\mathbf{v}_{h}]\!]}ds|\nonumber\\
&~~~\lesssim \|\mathbf{v}_{h}\|_{h} \|q_{h}\|_{0}+(\sum\limits_{F\in\mathcal{E}_{h}}\int_{F}\frac{\gamma}{h_{F}}
|[\![\underline{\mathbf{v}_{h}}]\!]|^{2}ds)^{\frac{1}{2}}(\sum\limits_{F\in \mathcal{E}_{h}}\int_{F}\frac{h_{F}}{\gamma}|\{ q_{h}\}|^{2}ds)^{\frac{1}{2}}\nonumber\\
&~~~\lesssim \|\mathbf{v}_{h}\|_{h}\|q_{h}\|_{0} + \|\mathbf{v}_{h}\|_{h}(\|q_{h}\|_{0}+(\sum\limits_{\tau\in\mathcal{T}_{h}}h_{\tau}^{2s}\|q_{h}\|_{s,\tau}^{2})^{^{\frac{1}{2}}})\\
&~~~\lesssim \|\mathbf{v}_{h}\|_{h}(\|q_{h}\|_{0} +(\sum\limits_{\tau\in\mathcal{T}_{h}}h_{\tau}^{2s}\|q_{h}\|_{s,\tau}^{2})^{^{\frac{1}{2}}})
~~~\forall [\mathbf{v}_{h},q_{h}]\in (H^{1}(\mathcal{T}_{h}))^{d}\times H^{s}(\mathcal{T}_{h}).\nonumber
\end{align}
\indent The following two estimates in $\mathbb{X}_{h}$ are a key in the proof of Lemma 2.1 below. \\
\indent{\bf Discrete Sobolev embedding} (see the Theorem 6.1 in \cite{Pietro2010}):
\begin{align}\label{s2.17}
\|\mathbf{v}_{h}\|_{0}^{2}\leq \sigma \|\mathbf{v}_{h}\|_{h}^{2}
~~~ \forall \mathbf{v}_{h} \in \mathbb{X}_{h}.
\end{align}
\indent{\bf The trace inequality }(see (2.8) and (2.9) in \cite{Riviere2008}):
\begin{align}\label{s2.18}
\sum\limits_{F\in\mathcal{E}_{h}}\|h^{\frac{1}{2}}\{\mathbf{v}_{h}\}\|_{0,F}^{2}\leq C_{k}\|\mathbf{v}_{h}\|_{0}
~~~ \forall \mathbf{v}_{h} \in \mathbb{X}_{h},
\end{align}
where $C_{k}$ is a function of the polynomial degree.\\
\begin{lemma}\label{lem2-1}
When the interior penalty parameter $\gamma$ is large enough, there exists a constant $\mu_{0}>0$ independent of $h$ such that the following discrete
coercivity property holds
\begin{align}\label{s2.19}
 Re A_{h}(\mathbf{u}_{h},\mathbf{u}_{h})\geq \mu_{0}  \|\mathbf{u}_{h}\|^{2}_{h},
~~~ \forall \mathbf{u}_{h} \in \mathbb{X}_{h}.
\end{align}
\end{lemma}

\noindent{\bf Proof.} By virtue of Green's formula, we deduce that
\begin{align*}
&Re\sum\limits_{\tau\in\mathcal{T}_{h}}\int_{\tau}(\boldsymbol{\beta}\cdot\nabla)\mathbf{u}_{h}\cdot \overline{\mathbf{u}_{h}}dx\nonumber\\
&~~~=\frac{1}{2}Re(\sum\limits_{F\in\mathcal{E}_{h}}\int_{F}  {\overline{\{\mathbf{u}_{h}\otimes\boldsymbol{\beta}\}}}: \underline{[\![\mathbf{u}_{h}]\!]} ds
+ \sum\limits_{F\in\mathcal{E}^{i}_{h}}\int_{F}  \{ \mathbf{u}_{h}  \} {\overline{[\![ \mathbf{u}_{h}\otimes\boldsymbol{\beta}]\!]} } ds),
\end{align*}
Hence, we obtain from Young's inequality with $\epsilon_{1},\epsilon_{2}>0$ and (\ref{s2.18}),
\begin{align}\label{s2.20}
&|Re(\sum\limits_{\tau\in\mathcal{T}_{h}}\int_{\tau}(\boldsymbol{\beta}\cdot\nabla)\mathbf{u}_{h}\cdot \overline{\mathbf{u}_{h}}dx-\frac{1}{2}\sum\limits_{\tau\in\mathcal{T}_{h}}\sum_{F\subset\partial\tau}\int_{F} [\![\mathbf{u}_{h}\otimes\boldsymbol{\beta}]\!]\cdot\overline{\{ \mathbf{u}_{h}\}} ds)|
 \nonumber\\
&=|\frac{1}{2}Re (\sum\limits_{F\in\mathcal{E}_{h}}\int_{F}  {\overline{\{\mathbf{u}_{h}\otimes\boldsymbol{\beta}\}}}:\underline{ [\![\mathbf{u}_{h}]\!]}ds
+ \sum\limits_{F\in\mathcal{E}_{h}^{i}}\int_{F}  \{ \mathbf{u}_{h}  \} \cdot{\overline{[\![\mathbf{u}_{h}\otimes\boldsymbol{\beta}]\!]} } ds\nonumber\\
&~~~-\sum_{F\in\mathcal{E}_{h}^{i}}\int_{F} [\![\mathbf{u}_{h}\otimes\boldsymbol{\beta}]\!]\cdot\overline{\{ \mathbf{u}_{h}\}} ds
-\frac{1}{2}\sum_{F\in\mathcal{E}_{h}^{b}}\int_{F} [\![\mathbf{u}_{h}\otimes\boldsymbol{\beta}]\!]\cdot\overline{\{ \mathbf{u}_{h}\}} ds) |  \nonumber\\
&\leq \sum\limits_{F\in\mathcal{E}_{h}}(\epsilon_{1}\|h_{F}^{\frac{1}{2}}\overline{\{\mathbf{u}_{h}\otimes\boldsymbol{\beta}\}}\|_{0,F}^{2}
+\frac{1}{\epsilon_{1}}
\|h_{F}^{-\frac{1}{2}}{\underline{[\![{\mathbf{u}_{h}}]\!]}}\|_{0,F}^{2}) + \sum\limits_{F\in\mathcal{E}_{h}}(\epsilon_{2}\|h_{F}^{\frac{1}{2}}{\{\mathbf{u}_{h}\}}\|_{0,F}^{2}\nonumber\\
&~~~+\frac{1}{\epsilon_{2}}
\|h_{F}^{-\frac{1}{2}}{\overline{[\![\mathbf{u}_{h}\otimes\boldsymbol{\beta}]\!]}}\|_{0,F}^{2})
\nonumber\\
&\leq \sum\limits_{\tau\in\mathcal{T}_{h}}C_{k}(\epsilon_{1}\|\boldsymbol{\beta}\|_{L^{\infty}}+\epsilon_{2}){\|\mathbf{u}_{h}\|}_{0,\tau}^{2}
+\sum\limits_{F\in\mathcal{E}_{h}}(\frac{1}{\epsilon_{1}}h^{-1}_{F}+ \frac{1}{\epsilon_{2}}h^{-1}_{F}\|\boldsymbol{\beta}\|_{L^{\infty}}  )\|{[\![\underline{\mathbf{u}_{h}}]\!]}\|_{0,F}^{2}.
\end{align}
Similarly, by Young's inequality with $\epsilon_{3}>0$ and (\ref{s2.18}), we have
\begin{align}\label{s2.21}
&~~~|Re(\sum\limits_{F\in\mathcal{E}_{h}}\int_{F}-\mu\{\nabla\mathbf{u}_{h}\}:\overline{\underline{[\![\mathbf{u}_{h}]\!]}}ds
+\theta\sum\limits_{F\in\mathcal{E}_{h}}\int_{F}\mu\{\overline{\nabla\mathbf{u}}_{h}\}:[\![\underline{\mathbf{u}_{h}}]\!]ds)|\nonumber\\
&\leq(1-\theta)\mu\sum\limits_{F\in\mathcal{E}_{h}}\|h_{F}^{\frac{1}{2}}
\{\nabla\mathbf{u}_{h}\}\|_{0,F}\|h_{F}^{-\frac{1}{2}}\overline{[\![\underline{\mathbf{u}_{h}}]\!]}\|_{0,F}\nonumber\\
&\leq \frac{1}{2}(1-\theta)\mu\sum\limits_{F\in\mathcal{E}_{h}}(\epsilon_{3}\|h_{F}^{\frac{1}{2}}
\{\nabla\mathbf{u}_{h}\}\|_{0,F}^{2}+\frac{1}{\epsilon_{3}}\|h_{F}^{-\frac{1}{2}}
\overline{[\![\underline{\mathbf{u}_{h}}]\!]}\|_{0,F}^{2})\nonumber\\
&\leq \frac{1}{2}(1-\theta)\mu(\sum\limits_{\tau\in\mathcal{T}_{h}}C_{k-1}\epsilon_{3}\|\nabla\mathbf{u}_{h}\|_{0,\tau}^{2}
+\sum\limits_{F\in\mathcal{E}_{h}}\frac{1}{\epsilon_{3}}\|h_{F}^{-\frac{1}{2}}\overline{[\![\underline{\mathbf{u}_{h}}]\!]}\|_{0,F}^{2}).
\end{align}
From (\ref{s2.20}), (\ref{s2.21}) and (\ref{s2.17}) we can derive
\begin{align*}
&Re(A_{h}(\mathbf{u}_{h},\mathbf{u}_{h}))\geq \mu\sum\limits_{\tau\in\mathcal{T}_{h}}\|\nabla \mathbf{u}_{h}\|_{0,\tau}^{2}-C_{k}(\epsilon_{1}\|\boldsymbol{\beta}\|_{L^{\infty}}+\epsilon_{2})\sum\limits_{\tau\in\mathcal{T}_{h}}{\|\mathbf{u}_{h}\|}_{0,\tau}^{2}\\
&~~~-\sum\limits_{F\in\mathcal{E}_{h}}(\frac{1}{\epsilon_{1}}h^{-1}_{F}+ \frac{1}{\epsilon_{2}}h^{-1}_{F} \|\boldsymbol{\beta}\|_{L^{\infty}} )\|{[\![\underline{\mathbf{u}_{h}}]\!]}\|_{0,F}^{2}\\
&~~~-\frac{1}{2}(1-\theta)\mu(C_{k-1}\epsilon_{3} \sum\limits_{\tau\in\mathcal{T}_{h}}\|\nabla \mathbf{u}_{h}\|_{0,\tau}^{2}
+\frac{1}{\epsilon_{3} } \sum\limits_{F\in\mathcal{E}_{h}}h_{F}^{-1}\|[\![\underline{\mathbf{u}_{h}}]\!]\|_{0,F}^{2})\\
&~~~+\sum\limits_{F\in\mathcal{E}_{h}}
\int_{F}\gamma h_{F}^{-1}
[\![\underline{\mathbf{u}_{h}}]\!]:\overline{[\![\underline{\mathbf{u}_{h}}]\!]}ds\\
& \geq
\left( \mu(1-\frac{1}{2}(1-\theta)C_{k-1}\epsilon_{3})-C_{k}(\epsilon_{1}\|\boldsymbol{\beta}\|_{L^{\infty}} +\epsilon_{2})\sigma \right)\sum\limits_{\tau\in\mathcal{T}_{h}}\|\nabla \mathbf{u}_{h}\|_{0,\tau}^{2}\\
&~~~+\left(\gamma -(\frac{1}{\epsilon_{1}}+\frac{1}{\epsilon_{2}}\|\boldsymbol{\beta}\|_{L^{\infty}} +\frac{1}{2}(1-\theta)\mu\frac{1}{\epsilon_{3}} )\right) \sum\limits_{F\in\mathcal{E}_{h}}h_{F}^{-1}\|{[\![\underline{{\mathbf{u}_{h}}}]\!]}\|_{0,F}^{2}.
\end{align*}
Choosing $\epsilon_{i}<1~(i=1,2,3)$ such that
\begin{align}\label{s2.22}
\mu(1-\frac{1}{2}(1-\theta)C_{k-1}\epsilon_{3})-C_{k}(\epsilon_{1}\|\boldsymbol{\beta}\|_{L^{\infty}} +\epsilon_{2})\sigma \geq \frac{\mu}{2},
\end{align}
 and choosing $\gamma$ such that
\begin{align}\label{s2.23}
\gamma -(\frac{1}{\epsilon_{1}}+\frac{1}{\epsilon_{2}}\|\boldsymbol{\beta}\|_{L^{\infty}} +\frac{1}{2}(1-\theta)\mu\frac{1}{\epsilon_{3}} )\geq \frac{1}{2}\gamma,
\end{align}
we obtain (\ref{s2.19}) with $\mu_{0}=\min\{\frac{1}{2},\frac{\mu}{2}\}$.
~~~~$\Box$\\

\indent The (\ref{s2.22}) and (\ref{s2.23}) provide the value of the penalty parameter  $\gamma$  and show that, to ensure the discrete coercivity of $A_{h}$, the ${\bf SIP}$ method requires a larger $\gamma$ compared with the ${\bf NIP}$ and ${\bf IIP}$ methods.\\
\indent The weak formulation of the boundary value problem (Oseen equation) associated with (\ref{s2.1}) reads:  find $[\mathbf{u}^{f},p^{f}]\in \mathbb{X}\times \mathbb{W}$ such that
\begin{align}
A(\mathbf{u}^{f},\mathbf{v})+B(\mathbf{v},p^{f})&=(\boldsymbol{f},\mathbf{v}),
~~~~\forall\mathbf{v}\in\mathbb{X},\label{s2.24}\\
\overline{B(\mathbf{u}^{f},q)}&=0,~~~~\forall q\in \mathbb{W}.\label{s2.25}
\end{align}
The DG formula of (\ref{s2.24})-(\ref{s2.25}) reads:
find  $(\mathbf{u}^{f}_{h},p^{f}_{h})\in \mathbb{X}_{h}\times \mathbb{W}_{h}$ such that
\begin{align}
A_{h}(\mathbf{u}^{f}_{h},\mathbf{v}_{h})+B_{h}(\mathbf{v}_{h},p^{f}_{h})&=(\boldsymbol{f},\mathbf{v}_{h}),
~~~~\forall\mathbf{v}_{h}\in\mathbb{X}_{h},\label{s2.26}\\
\overline{B_{h}(\mathbf{u}^{f}_{h},q_{h})}&=0,~~~~\forall q_{h}\in \mathbb{W}_{h}.\label{s2.27}
\end{align}

\indent We introduce the following sesquilinear form
\begin{align}\label{s2.28}
\mathcal{A}_{h}([\mathbf{w},r], [\mathbf{v},q])=A_{h}(\mathbf{w},\mathbf{v})+B_{h}(\mathbf{v},r)
+\overline{B_{h}(\mathbf{w},q)},~~~ \forall [\mathbf{w},r], [\mathbf{v},q]\in \mathbb{V}_{h}.
\end{align}
\indent Using the standard argument (see  Lemma 6.5 and Proposition 2.9 of \cite{Riviere2008}), we can derive the consistency lemma of (\ref{s2.26})-(\ref{s2.27}).

\begin{lemma}\label{lem2-2}
Let $[\mathbf{u}^{f},p^{f}]$ be the solution to the boundary value problem (\ref{s2.24})-(\ref{s2.25}). Then
for any $\theta\in \{-1,0,1\}$ there holds
\begin{align}
\mathcal{A}_{h}([\mathbf{u}^{f},p^{f}], [\mathbf{v}_{h},q_{h}])=(\boldsymbol{f},\mathbf{v}_{h}),  ~~\forall [\mathbf{v}_{h},q_{h}]\in \mathbb{X}(h)\times \mathbb{W},\label{s2.29}
\end{align}
\end{lemma}
\indent For the proof of the following discrete inf-sup condition and a priori error estimate in $\mathbb{R}^{d}$, we need
to introduce the projection operator $\Pi_{h}$ of N\'{e}d\'{e}lec, see \cite{Nedelec1986},
and the local $L^{2}$ projection operator $\vartheta_{h}$.\\
\indent The N\'{e}d\'{e}lec $H(div)$ space of index $k$  is given by
\begin{align*}
V_{h}^{N}=\{\mathbf{v}\in H(\Omega;div):  \mathbf{v}|_{\tau}\in[P_{k}(\tau)]^{d}, ~\forall \tau\in \mathcal{T}_{h}\}.
\end{align*}
\indent The projection $\Pi_{h}:[H^{1}(\Omega)]^{d}\rightarrow V_{h}^{N}$ is defined as follows for every $\tau\in \mathcal{T}_{h}$:
\begin{align}
\int_{\tau}(\mathbf{v}-\Pi_{h} \mathbf{v})\cdot \overline{\mathbf{w}}&=0,  \forall \mathbf{w}\in N_{k-1}(\tau),\label{s2.30}\\
\int_{F}(\mathbf{v}-\Pi_{h} \mathbf{v})\cdot \mathbf{n}\overline{q}&=0, \forall q\in P_{k}(F), \forall F\subset \partial\tau,\label{s2.31}
\end{align}
where $N_{k-1}(\tau)$ is the N\'{e}d\'{e}lec space of index $k-1$ given by
\begin{align*}
N_{k-1}(\tau)=[P_{k-2}(\tau)]^{d}+\{\mathbf{v}\in[\widetilde{P}_{k-1}(\tau)]^{d}:\mathbf{v}\cdot \mathbf{x}=0\},
\end{align*}
and $\widetilde{P}_{k-1}(\tau)$ is the space of homogeneous polynomials of degree $k-1$.\\
\indent  We obtain from  (\ref{s2.30}) and (\ref{s2.31})
\begin{align}\label{s2.32}
&\int_{\tau}\nabla\cdot (\Pi_{h}\mathbf{v}-\mathbf{v})\overline{q}dx=-\int_{\tau} (\Pi_{h}\mathbf{v}-\mathbf{v})\cdot \overline{\nabla q}dx\nonumber\\
&~~~~~~
+\int_{\partial\tau} (\Pi_{h}\mathbf{v}-\mathbf{v})\cdot\mathbf{n}\overline{q} ds=
0,~~ \forall \mathbf{v}\in H^{1}(\Omega)^{d}, ~~\forall q\in P_{k-1}(\tau).
\end{align}
From Proposition 1 in \cite{Nedelec1986}, if $\mathbf{u}\in (H^{1+s}(\Omega))^{d}$ for any $k\geq s\geq 0$ holds
\begin{align}\label{s2.33}
\|\mathbf{u}-\Pi_{h} \mathbf{u}\|_{h}\lesssim h^{s}\|\mathbf{u}\|_{1+s}.
\end{align}
\indent Let $\vartheta_{h}:~\mathbb{W}\rightarrow \mathbb{W}_{h}$ be the local $L^{2}$ projection operator satisfying
$\vartheta_{h}p|_{\tau}\in \mathbb{P}_{k-1}(\tau)$ and
\begin{eqnarray*}
 \int_{\tau}(p-\vartheta_{h}p)\overline{q}dx=0,~~\forall q\in P_{k-1}(\tau),~~~\forall \tau\in\mathcal{T}_{h};
\end{eqnarray*}
Then, from Section 6.1.5 in \cite{Riviere2008}, the following holds when $p\in H^{s}(\tau)$ $(0<s\leq k)$
\begin{eqnarray}\label{s2.34}
 \|p-\vartheta_{h}p\|_{0,\tau}\lesssim h^{s}\|p\|_{s,\tau},~~~\forall \tau\in\mathcal{T}_{h}.
\end{eqnarray}

\indent The following lemma is classical, with a standard proof argument (see Section VI.4 in \cite{Brezzi1991} and Proposition 10 in \cite{Hansbo2002}); for the convenience of readers, we present the proof for the complex space $\mathbb{X}_{h}\times \mathbb{W}_{h}$ and in $\mathbb{R}^{d}$.
\begin{lemma}\label{lem2-3}
There exists a constant $\alpha>0$ independent of $h$ such that the following discrete inf-sup condition holds
\begin{align}\label{s2.35}
\inf\limits_{0\neq q_{h}\in \mathbb{W}_{h}}\sup\limits_{\mathbf{0}\neq\mathbf{v}_{h}\in {\mathbb{X}}_{h}  }\frac{Re B_{h}(\mathbf{v}_{h},q_{h})}{\|q_{h}\|_{0}\|\mathbf{v}_{h}\|_{h}}\geq \alpha.
\end{align}
\end{lemma}
\noindent {\bf Proof.}
The following commuting property holds (\cite{Nedelec1986}):
\begin{align}\label{s2.36}
\nabla\cdot \Pi_{h}\mathbf{v}=\vartheta_{h}\nabla\cdot \mathbf{v}.
\end{align}
In addition, the following bound holds:
\begin{align}\label{s2.37}
\|\Pi_{h}\mathbf{v}\|_{h}\leq C \|\mathbf{v}\|_{1}, \forall \mathbf{v}\in[H_{0}^{1}(\Omega)]^{d}.
\end{align}
From \cite{Acosta2006}, $\forall p_{h}\in \mathbb{W}_{h}$, there exists a $\mathbf{v}\in [H_{0}^{1}(\Omega)]^{d}$, such that
\begin{align*}
\nabla\cdot \mathbf{v}=p_{h},~~~in~\Omega,~~~
\mathbf{v}=0,~~~on~\partial\Omega,
\end{align*}
\begin{align}\label{s2.38}
\|\mathbf{v}\|_{1}\leq C\|p_{h}\|_{0}.
\end{align}
By (\ref{s2.36}), we have that $\nabla\cdot \Pi_{h}\mathbf{v}=p_{h}$, and from (\ref{s2.30}) and (\ref{s2.31}), it is easy to know that $[\![\Pi_{h}\mathbf{v}]\!]=0$, $\forall F\in \mathcal{E}_{h}^{i}$, and together with $\Pi_{h}\mathbf{v}\cdot \mathbf{n}=0$ on $\partial\Omega$, then
\begin{align}\label{s2.39}
\|p_{h}\|_{0}^{2}=\int_{\Omega}p_{h}\overline{p_{h}}dx=\int_{\Omega}p_{h}\overline{\nabla\cdot \Pi_{h}\mathbf{v} }dx
=-B_{h}(\Pi_{h}\mathbf{v},p_{h}).
\end{align}
We get from (\ref{s2.37}) and (\ref{s2.38})
\begin{align*}
-\frac{B_{h}(\Pi_{h}\mathbf{v},p_{h})}{\|\Pi_{h}\mathbf{v}\|_{h}}\gtrsim -\frac{B_{h}(\Pi_{h}\mathbf{v},p_{h})}{\|\mathbf{v}\|_{1}} \gtrsim -\frac{B_{h}(\Pi_{h}\mathbf{v},p_{h})}{\|p_{h}\|_{0}}=\|p_{h}\|_{0}, ~~~\forall p_{h}\in \mathbb{W}_{h}.
\end{align*}
Hence, there exists a $\alpha>0$ independent of $h$, such that (\ref{s2.35}) holds.~~~$\Box$\\

\indent From (\ref{s2.4}) and (\ref{s2.5}), we know that (\ref{s2.24})-(\ref{s2.25}) are  well-posed.
Then we define the solution operators
$\mathbf{T}: L^{2}(\Omega)^{d}\rightarrow \mathbb{X},~\mathbf{T}\boldsymbol{f}=\mathbf{u}^{f}$, and
$\mathbf{S}: L^{2}(\Omega)^{d}\rightarrow \mathbb{W}, ~\mathbf{S}\boldsymbol{f}=p^{f}$ as follows:
\begin{align}
A(\mathbf{T}\boldsymbol{f},\mathbf{v})+B(\mathbf{v},\mathbf{S}\boldsymbol{f})&=(\boldsymbol{f},\mathbf{v}),
~~~~\forall\mathbf{v}\in\mathbb{X},\label{s2.40}\\
\overline{B(\mathbf{T}\boldsymbol{f},q)}&=0,~~~~\forall q\in \mathbb{W}.\label{s2.41}
\end{align}
And it is valid that
\begin{align}
&\|\mathbf{T}\boldsymbol{f}\|_{1}+\|\mathbf{S}\boldsymbol{f}\|_{0}\lesssim\|\boldsymbol{f}\|_{0}.~~~\label{s2.42}
\end{align}
Because of the compact embedding $\mathbb{X}\hookrightarrow L^{2}(\Omega)^{d}$, $\mathbf{T}: L^{2}(\Omega)^{d}\rightarrow L^{2}(\Omega)^{d}$
is compact.\\
\indent By (\ref{s2.19}) and (\ref{s2.35}), we know that (\ref{s2.26})-(\ref{s2.27}) are well-posed, and we define the discrete solution operators
$\mathbf{T}_{h}: L^{2}(\Omega)^{d}\rightarrow \mathbb{X}_{h}, ~\mathbf{T}_{h}\boldsymbol{f}=\mathbf{u}^{f}_{h}$, and
$\mathbf{S}_{h}: L^{2}(\Omega)^{d}\rightarrow \mathbb{W}_{h}, ~\mathbf{S}_{h}\boldsymbol{f}=p^{f}_{h}$
as follows:
\begin{align}
A_{h}(\mathbf{T}_{h}\boldsymbol{f},\mathbf{v})+B_{h}(\mathbf{v},\mathbf{S}_{h}\boldsymbol{f})&=(\boldsymbol{f},\mathbf{v}),
~~~~\forall\mathbf{v}\in\mathbb{X}_{h},\label{s2.43}\\
\overline{B_{h}(\mathbf{T}_{h}\boldsymbol{f},q)}&=0,~~~~\forall q\in \mathbb{W}_{h}.\label{s2.44}
\end{align}
And there holds
\begin{align}
&||\mathbf{T}_{h}\boldsymbol{f}||_{h}+\|\mathbf{S}_{h}\boldsymbol{f}\|_{0}\lesssim\|\boldsymbol{f}\|_{0}.~~~\label{s2.45}
\end{align}
 Thus, (\ref{s2.2})-(\ref{s2.3}) and (\ref{s2.8})-(\ref{s2.9}) have the following equivalent operator forms, respectively:
\begin{align}
& \lambda \mathbf{T} \mathbf{u}=\mathbf{u},~~~~~~~~~~\mathbf{S} (\lambda \mathbf{u})=p,\label{s2.46}\\
&  \lambda_{h} \mathbf{T}_{h} \mathbf{u}_{h}=\mathbf{u}_{h},~~~~\mathbf{S}_{h} (\lambda_{h} \mathbf{u}_{h})=p_{h}.\label{s2.47}
\end{align}
\indent Referring to \cite{Babuska1991book}, we give the adjoint eigenvalue problem of (\ref{s2.2})-(\ref{s2.3}) and its
 DG formula:
find $(\lambda^{*},\mathbf{u}^{*},p^{*})\in  \mathbb{C}\times \mathbb{X}\times \mathbb{W}$,
$\|\mathbf{u}^{*}\|_{0}=1$, such that
\begin{align}\label{s2.48}
\mathcal{A}([\mathbf{v},q], [\mathbf{u}^{*},p^{*}])=\overline{\lambda^*}(\mathbf{v}, \mathbf{u}^{*}),  ~~\forall [\mathbf{v},q]\in \mathbb{V},
\end{align}
 and find $({\lambda^{*}_{h}}, \mathbf{u}^{*}_{h}, p^{*}_{h})\in \mathbb{C}\times\mathbb{X}_{h}\times \mathbb{W}_{h}$,
$\|\mathbf{u}^{*}_{h}\|_{0}=1$, such that
\begin{align}\label{s2.49}
\mathcal{A}_{h}([\mathbf{v}_{h},q_{h}], [\mathbf{u}_{h}^{*},p_{h}^{*}])=\overline{\lambda_{h}^*}(\mathbf{v}_{h}, \mathbf{u}_{h}^{*}),  ~~\forall [\mathbf{v}_{h},q_{h}]\in \mathbb{V}_{h}.
\end{align}
\indent The source problem associated with (\ref{s2.48}) is to find $[\mathbf{u}^{f*},p^{f*}]\in \mathbb{V}$ such that
\begin{align}\label{s2.50}
\mathcal{A}([\mathbf{v},q], [\mathbf{u}^{f*},p^{f*}])=(\mathbf{v}, \boldsymbol{f}),  ~~\forall [\mathbf{v},q]\in \mathbb{V},
\end{align}
and its DG formula reads:
find  $[\mathbf{u}^{f*}_{h},p_{h}^{f*}]\in \mathbb{V}_{h}$ such that
\begin{align}\label{s2.51}
\mathcal{A}_{h}([\mathbf{v}_{h},q_{h}], [\mathbf{u}_{h}^{f*},p_{h}^{f*}])&=(\mathbf{v}_{h}, \boldsymbol{f}),  ~~\forall [\mathbf{v}_{h},q_{h}]\in \mathbb{V}_{h}.
\end{align}

\indent From (\ref{s2.4}) and (\ref{s2.5}), we can define the solution operators $\mathbf{T}^{*}: L^{2}(\Omega)^{d}\rightarrow \mathbb{X}\subset L^{2}(\Omega)^{d}$,
$\mathbf{S}^{*}: L^{2}(\Omega)^{d}\rightarrow \mathbb{W}$,
and
it is valid that
\begin{align}
&\|\mathbf{T}^{*}\boldsymbol{f}\|_{1}+\|\mathbf{S}^{*}\boldsymbol{f}\|_{0}\lesssim\|\boldsymbol{f}\|_{0}.~~~\label{s2.52}
\end{align}
\indent By (\ref{s2.19}) and (\ref{s2.35}), we can define the discrete solution operators $\mathbf{T}_{h}^{*}: L^{2}(\Omega)^{d}\rightarrow \mathbb{X}_{h}\subset L^{2}(\Omega)^{d}$, $\mathbf{S}_{h}^{*}: L^{2}(\Omega)^{d}\rightarrow \mathbb{W}_{h}$,
and it is valid that
\begin{align}
&||\mathbf{T}_{h}^{*}\boldsymbol{f}||_{h}+\|\mathbf{S}_{h}^{*}\boldsymbol{f}\|_{0}\lesssim\|\boldsymbol{f}\|_{0}.~~~\label{s2.53}
\end{align}
And we can deduce that
\begin{align}\label{s2.54}
&(\mathbf{T}\mathbf{v},\mathbf{z})= (\mathbf{v},\mathbf{T}^{*}\mathbf{z}),
~~~ \forall  \mathbf{v},\mathbf{z}\in (L^{2}(\Omega))^{d},\\\label{s2.55}
&(\mathbf{T}_{h}\mathbf{v},\mathbf{z})=(\mathbf{v},\mathbf{T}_{h}^{*}\mathbf{z}),~~~\forall  \mathbf{v},\mathbf{z}\in (L^{2}(\Omega))^{d}.
\end{align}
Then $\mathbf{T}^{*}$ and $\mathbf{T}_{h}^{*}$ are adjoint operators of $\mathbf{T}$ and $\mathbf{T}_{h}$, respectively.
Therefore the prime and adjoint eigenvalues are connected via $\lambda = \overline{\lambda^*}$ and $\lambda_h = \overline{\lambda_h^*}$.\\
\indent  Using the well-posed results (\ref{s2.42}) and (\ref{s2.52}), and referring to the Stokes  regularity result (see
 \cite{Fabes1988,Savare1998, Kellogg1976} and Theorem 1 in \cite{DG+Gedicke2020} for instance), we assume that the following regularity result $\mathbf{R(\Omega)}$ is valid.\\
  \indent $\mathbf{R(\Omega)}$.~For any $\boldsymbol{f}\in L^{2}(\Omega)^{d}$, there exists $r\in (\frac{1}{2},1]$, such that
  $[\mathbf{T}\boldsymbol{f}, \mathbf{S}\boldsymbol{f}]$, $[\mathbf{T}^{*}\boldsymbol{f}, \mathbf{S}^{*}\boldsymbol{f}]\in H^{1+r}(\Omega)^{d}\times H^{r}(\Omega)$ and
\begin{eqnarray}\label{s2.56}
\|\mathbf{T}\boldsymbol{f}\|_{1+r}+\|\mathbf{S}\boldsymbol{f}\|_{r}\lesssim C_{\Omega}\|\boldsymbol{f}\|_{0},~~~\|\mathbf{T}^*\boldsymbol{f}\|_{1+r}
+\|\mathbf{S}^*\boldsymbol{f}\|_{r}\lesssim C_{\Omega}\|\boldsymbol{f}\|_{0},
\end{eqnarray}
where $C_{\Omega}$ is a positive constant independent of $\boldsymbol{f}$.\\
\indent Referring to Proposition 2.9 of \cite{Riviere2008}, next we will discuss the adjoint-consistency of the DG method.
\begin{lemma}\label{lem2-4}
Let $[\mathbf{u}^{f*},p^{f*}]$ be the solution to the boundary value problem (\ref{s2.50}).
 Then, for any $\theta\in \{-1,0,1\}$, there holds
\begin{align}\label{s2.57}
&\mathcal{A}_{h}([\mathbf{v}_{h},q_{h}], [\mathbf{u}^{f*},p^{f*}])-
(1+\theta)\sum\limits_{F\in\mathcal{E}_{h}}\int_{F}\mu\overline{\{{\nabla\mathbf{u}}^{f*}\}}:[\![\underline{\mathbf{v}_{h}}]\!]ds
=(\mathbf{v}_{h}, \boldsymbol{f}),\nonumber\\
&~~~~~~ \forall [\mathbf{v}_{h},q_{h}]\in \mathbb{X}(h)\times
\mathbb{W}.
\end{align}
\end{lemma}
\noindent {\bf Proof.}
A simple calculation shows that on element $\tau$, there hold
\begin{align*}
(\boldsymbol{\beta}\cdot\nabla)\mathbf{v}_{h}=div(\mathbf{v}_{h}\otimes\boldsymbol{\beta}),~~~
(\mathbf{v}_{h}\otimes\boldsymbol{\beta}):\overline{\nabla \mathbf{u}^{f*}}=\mathbf{v}_{h}\cdot\overline{(\boldsymbol{\beta}\cdot\nabla)\mathbf{u}^{f*}}.
\end{align*}
From the above two equations and Green's formula, we deduce that
\begin{align*}
&\sum\limits_{\tau\in\mathcal{T}_{h}}\int_{\tau}(\boldsymbol{\beta}\cdot\nabla)\mathbf{v}_{h}\cdot\overline{{\mathbf{u}^{f*}}}dx
=\sum\limits_{\tau\in\mathcal{T}_{h}}\int_{\tau}div(\mathbf{v}_{h}\otimes\boldsymbol{\beta})\cdot\overline{{\mathbf{u}^{f*}}}dx\nonumber\\
&~~~=-\sum\limits_{\tau\in\mathcal{T}_{h}}\int_{\tau}\mathbf{v}_{h}\cdot\overline{(\boldsymbol{\beta}\cdot\nabla)\mathbf{u}^{f*}}dx
+\frac{1}{2}\sum\limits_{\tau\in\mathcal{T}_{h}}\sum\limits_{F\subset\partial\tau}\int_{F}
[\![(\mathbf{v}_{h}\otimes\boldsymbol{\beta})]\!]
\cdot\overline{\{\mathbf{u}^{f*}\}}ds.\nonumber\\
&\sum\limits_{\tau\in\mathcal{T}_{h}}\int_{\tau}\mu\nabla \mathbf{v}_{h}:\overline{\nabla\mathbf{u}^{f*}}dx
=\sum\limits_{\tau\in\mathcal{T}_{h}}(-\int_{\tau}\mathbf{v}_{h}:\overline{\mu\Delta\mathbf{u}^{f*}}dx+
\sum\limits_{F\subset\partial\tau}\int_{F}\mathbf{v}_{h}\otimes\mathbf{n}:\mu\overline{{\nabla\mathbf{u}}^{f*}}ds)\nonumber\\
&~~~=-\sum\limits_{\tau\in\mathcal{T}_{h}}\int_{\tau}\mathbf{v}_{h}:\overline{\mu\Delta\mathbf{u}^{f*}}dx+
\sum\limits_{F\in\mathcal{E}_{h}}\int_{F}[\![\underline{\mathbf{v}_{h}}]\!]:\mu\overline{\{{\nabla\mathbf{u}}^{f*}\}}ds.
\end{align*}
Using (\ref{s2.10}), (\ref{s2.11}) and the above two equations and $[\![\underline{\mathbf{u}^{f*}}]\!]=0$ we deduce
\begin{align*}
&\mathcal{A}_{h}([\mathbf{v}_{h},q_{h}], [\mathbf{u}^{f*},p^{f*}])=A_{h}(\mathbf{v}_{h}, \mathbf{u}^{f*})+\overline{B_{h}(\mathbf{v}_{h}, p^{f*})}\nonumber\\
&=-\sum\limits_{\tau\in\mathcal{T}_{h}}\int_{\tau}\mathbf{v}_{h}:\overline{\mu\Delta\mathbf{u}^{f*}}dx+
\sum\limits_{F\in\mathcal{E}_{h}}\int_{F}[\![\underline{\mathbf{v}_{h}}]\!]:\mu\overline{\{{\nabla\mathbf{u}}^{f*}\}}ds\nonumber\\
&~~~-\sum\limits_{\tau\in\mathcal{T}_{h}}\int_{\tau}\mathbf{v}_{h}\cdot\overline{(\boldsymbol{\beta}\cdot\nabla)\mathbf{u}^{f*}}dx
+\frac{1}{2}\sum\limits_{\tau\in\mathcal{T}_{h}}\sum\limits_{F\subset\partial\tau}\int_{F}
[\![(\mathbf{v}_{h}\otimes\boldsymbol{\beta})]\!]
\cdot\overline{\{\mathbf{u}^{f*}\}}ds\nonumber\\
&~~~-\sum\limits_{F\in\mathcal{E}_{h}}\int_{F}\mu\{\nabla\mathbf{v}_{h}\}:\overline{[\![\underline{{{\mathbf{u}}^{f*}}}]\!]}ds
+\theta\sum\limits_{F\in\mathcal{E}_{h}}\int_{F}\mu\overline{\{{\nabla\mathbf{u}}^{f*}\}}:[\![\underline{\mathbf{v}_{h}}]\!]ds\nonumber\\
&~~~+\sum\limits_{F\in\mathcal{E}_{h}}\int_{F}\frac{\gamma}{h_{F}}[\![\underline{\mathbf{v}_{h}}]\!]:\overline{[\![\underline{{{\mathbf{u}}^{f*}}}]\!]}ds
-\frac{1}{2}\sum\limits_{\tau\in\mathcal{T}_{h}}\sum_{F\subset\partial\tau}\int_{F} [\![\mathbf{v}_{h}\otimes\boldsymbol{\beta}]\!]\cdot\overline{\{ \mathbf{u}^{f*}\}} ds\nonumber\\
&~~~+\sum\limits_{\tau\in\mathcal{T}_{h}}(\int_{\tau}\overline{\nabla p^{f*}}\mathbf{v}_{h}dx-
\sum\limits_{F\subset\partial\tau}
\int_{F}\overline{ p^{f*}}\mathbf{v}_{h}\cdot\mathbf{n}ds)
+\sum\limits_{F\in\mathcal{E}_{h}}\int_{F}\{\overline{p^{f*}}\}{[\![\mathbf{v}_{h}]\!]}ds\nonumber\\
&=-\sum\limits_{\tau\in\mathcal{T}_{h}}\int_{\tau}\mathbf{v}_{h}\cdot\overline{\mu\Delta\mathbf{u}^{f*}}dx
-\sum\limits_{\tau\in\mathcal{T}_{h}}\int_{\tau}\mathbf{v}_{h}\cdot\overline{(\boldsymbol{\beta}\cdot\nabla)\mathbf{u}^{f*}}dx-0+0
\nonumber\\
&~~~+\sum\limits_{\tau\in\mathcal{T}_{h}}(\int_{\tau}\overline{\nabla p^{f*}}\cdot\mathbf{v}_{h}dx-
\sum\limits_{F\subset\partial\tau}
\int_{F}\overline{ p^{f*}}\mathbf{v}_{h}\cdot\mathbf{n}ds)
+\sum\limits_{F\in\mathcal{E}_{h}}\int_{F}\{\overline{p^{f*}}\}{[\![\mathbf{v}_{h}]\!]}ds\nonumber\\
&~~~+(1+\theta)\sum\limits_{F\in\mathcal{E}_{h}}\int_{F}\mu\overline{\{{\nabla\mathbf{u}}^{f*}\}}:[\![\underline{\mathbf{v}_{h}}]\!]ds\nonumber.
\end{align*}
Then, (\ref{s2.57}) holds.
~~~$\Box$\\
\indent  Lemma 2.4 shows that the {\bf SIP} method is adjoint consistent, but the {\bf NIP} and {\bf IIP} methods are not adjoint consistent.\\
\indent
From the inverse inequality and the projection estimate, we deduce
\begin{align}\label{s2.58}
&h_{\tau}^{s}|\nabla(\mathbf{u}-\mathbf{v}_{h})|_{s,\tau} \lesssim
(|\nabla(\mathbf{u}-\mathbf{v}_{h})|_{0,\tau}
+h_{\tau}^{s}
|\mathbf{u}-\Pi_{h}\mathbf{u}|_{1+s,\tau}),\\\label{s2.59}
&h_{\tau}^{s}|p-p_{h}|_{s,\tau} \lesssim
(|p-p_{h}|_{0,\tau}
+h_{\tau}^{s}
|p-\vartheta_{h}p|_{s,\tau}),~\forall \mathbf{v}_{h}\in \mathbb{X}_{h}, p_{h}\in\mathbb{W}_{h}.
\end{align}
 The above two estimates are useful in the proofs of Theorems 2.5 and 3.7.
\begin{theorem}\label{th2-5}
Let $[\mathbf{u}^{f}, p^{f}]$ and $[\mathbf{u}_{h}^{f}, p_{h}^{f}]$ be the solution of (\ref{s2.24})-(\ref{s2.25})
and (\ref{s2.26})-(\ref{s2.27}), respectively,
and let $[\mathbf{u}^{f*}, p^{f*}]$ and $[\mathbf{u}_{h}^{f*}, p_{h}^{f*}]$ be the solution of (\ref{s2.50})
and (\ref{s2.51}), respectively.
Assume that $[\mathbf{u}^{f}, p^{f}],~[\mathbf{u}^{f*}, p^{f*}]\in H^{1+s}(\Omega)^{d}\times H^{s}(\Omega)$, $\frac{1}{2}< s\leq k$,
then, we have
\begin{align}\label{s2.60}
&|||\mathbf{u}^{f}-\mathbf{u}^{f}_{h}|||+\|p^{f}-p^{f}_{h}\|_{0}\lesssim h^{s}(|\mathbf{u}^{f}|_{1+s}+|p^{f}|_{s})~~(\theta\in \{-1,0,1\}),\\
\label{s2.61}
&|||\mathbf{u}^{f*}-\mathbf{u}^{f*}_{h}|||+\|p^{f*}-p^{f*}_{h}\|_{0}\lesssim h^{s}(|\mathbf{u}^{f*}|_{1+s}+|p^{f*}|_{s})~~~(\theta=-1).
\end{align}
\end{theorem}
\noindent {\bf Proof.}
 By the consistency (\ref{s2.29}), (\ref{s2.15}), (\ref{s2.16}), (\ref{s2.58}), (\ref{s2.59}), (\ref{s2.33}) and (\ref{s2.34})
we deduce that
\begin{align*}
&\mathcal{A}_{h}( [\mathbf{u}_{h}^{f}-\Pi_{h}\mathbf{u}^{f}, p_{h}^{f}-\vartheta_{h}p^{f}], [\psi_{h},\varphi_{h}])
=\mathcal{A}_{h}( [\mathbf{u}^{f}-\Pi_{h}\mathbf{u}^{f}, p^{f}-\vartheta_{h}p^{f}], [\psi_{h},\varphi_{h}])\nonumber\\
&~~~=A_{h}( (\mathbf{u}^{f}-\Pi_{h}\mathbf{u}^{f}, \psi_{h})+B_{h}(\psi_{h}, p^{f}-\vartheta_{h}p^{f})
+\overline{B_{h}(\mathbf{u}^{f}-\Pi_{h}\mathbf{u}^{f}, \varphi_{h})}\nonumber\\
&~~~\lesssim h^{s}(\|\mathbf{u}^{f}\|_{1+s}+\|p^{f}\|_{s})(\|\psi_{h}\|_{h}+\|\varphi_{h}\|_{0}),~~~\forall (\psi_{h},\varphi_{h})\in \mathbb{V}_{h}.
\end{align*}
From Lemmas \ref{lem2-1} and \ref{lem2-3} we can deduce that
$\mathcal{A}_{h}$ satisfy the discrete inf-sup condition. Thus we obtain
\begin{align*}
&\alpha(\|\mathbf{u}_{h}^{f}-\Pi_{h}\mathbf{u}^{f}\|_{h}+\| p_{h}^{f}-\vartheta_{h}p^{f}\|_{0})\leq
\sup\limits_{(\psi_{h},\varphi_{h})\in \mathbb{V}_{h}}\frac{|\mathcal{A}_{h}( [\mathbf{u}_{h}^{f}-\Pi_{h}\mathbf{u}^{f}, p_{h}^{f}-\vartheta_{h}p^{f}], [\psi_{h},\varphi_{h}])|}{\|\psi_{h}\|_{h}+\|\varphi_{h}\|_{0}}\nonumber\\
&~~~\lesssim h^{s}(\|\mathbf{u}^{f}\|_{1+s}+\|p^{f}\|_{s}) .
\end{align*}
Applying triangle inequality and the above estimate yield (\ref{s2.60}).\\
\indent When $\theta=-1$, by following the steps of argument as for (\ref{s2.60}),
 we can deduce that the estimate (\ref{s2.61}) holds.~~~~$\Box$\\
\indent Due to the lack of conjugate consistency, we are unable to prove that error estimate (\ref {s2.61}) holds for either the {\bf NIP} or {\bf IIP} method.
\begin{theorem}\label{th2-6}
Assume that $\mathbf{R(\Omega)}$ holds.
 Then for $\theta=-1$ there hold
\begin{align}\label{s2.62}
&\|\mathbf{u}^{f*}-\mathbf{u}^{f*}_{h}\|_{0}\lesssim h^{r}(|||\mathbf{u}^{f*}-\mathbf{u}^{f*}_{h}|||+\|p^{f*}-p^{f*}_{h}\|_{0}),\\
\label{s2.63}
&\|\mathbf{u}^{f}-\mathbf{u}^{f}_{h}\|_{0}\lesssim h^{r}(|||\mathbf{u}^{f}-\mathbf{u}_{h}^{f}|||+\|p^{f}-p^{f}_{h}\|_{0}).
\end{align}
\end{theorem}
\noindent {\bf Proof.}~
Let $\boldsymbol{f}=\mathbf{u}^{f*}-\mathbf{u}^{f*}_{h}$, from (\ref{s2.29}), (\ref{s2.57}),
 we derive
\begin{align*}
&\|\mathbf{u}^{f*}-\mathbf{u}^{f*}_{h}\|_{0}^{2}=(\boldsymbol{f},\mathbf{u}^{f*}-\mathbf{u}^{f*}_{h})
=\mathcal{A}_{h}([\mathbf{T}\boldsymbol{f}, \mathbf{S}\boldsymbol{f}], [\mathbf{u}^{f*}-\mathbf{u}^{f*}_{h},p^{f*}-p^{f*}_{h}])\nonumber\\
&~~~=\mathcal{A}_{h}([\mathbf{T}\boldsymbol{f}-\mathbf{T}_{h}\boldsymbol{f}, \mathbf{S}\boldsymbol{f}-\mathbf{S}_{h}\boldsymbol{f}], [\mathbf{u}^{f*}-\mathbf{u}^{f*}_{h},p^{f*}-p^{f*}_{h}])\nonumber\\
&~~~\lesssim |||\mathbf{T}\boldsymbol{f}-\mathbf{T}_{h}\boldsymbol{f}|||~||| \mathbf{u}^{f*}
-\mathbf{u}^{f*}_{h}|||+\|\mathbf{T}_{h}\boldsymbol{f}-\mathbf{T}\boldsymbol{f}\|_{h} \| p^{f*}-p^{f*}_{h}\|_{0}\nonumber\\
&~~~+\|\mathbf{u}^{f*}-\mathbf{u}^{f*}_{h}\|_{h}(\|\mathbf{S}\boldsymbol{f}-\mathbf{S}_{h}\boldsymbol{f}\|_{0}+
(\sum\limits_{\tau\in\mathcal{T}_{h}}h_{\tau}^{2r}\|\mathbf{S}\boldsymbol{f}-\mathbf{S}_{h}\boldsymbol{f}\|_{r,\tau}^{2})^{\frac{1}{2}}),
\end{align*}
thus, from (\ref{s2.60}), (\ref{s2.61}), (\ref{s2.59}) and (\ref{s2.56}), we obtain (\ref{s2.62}).\\
\indent By the same argument as for (\ref{s2.62}), we obtain (\ref{s2.63}).~~~$\Box$\\
\indent Assume that $\mathbf{R(\Omega)}$ holds, then combining (\ref{s2.60}) and (\ref{s2.56}), and noting $\|\mathbf{T}_{h}^{*}-\mathbf{T}^{*}\|_{0}=\|\mathbf{T}_{h}-\mathbf{T}\|_{0}$,
we can deduce that for any $\theta\in \{-1,0,1\}$ there hold
\begin{align}
&\|\mathbf{T}_{h}-\mathbf{T}\|_{0}\lesssim h^{r},
~~~\|\mathbf{T}_{h}^{*}-\mathbf{T}^{*}\|_{0}\lesssim h^{r}.\label{s2.64}
\end{align}
\indent  Let $\{\lambda_{j}\}$ and $\{\lambda_{j,h}\}$ be enumerations of the eigenvalues of (\ref{s2.2})-(\ref{s2.3}) and (\ref{s2.8})-(\ref{s2.9}) respectively,
 by their modulus from small to large, each repeated as many times as its multiplicity. Suppose that $\lambda$ and $\lambda_{h}$ are the $k$th eigenvalue of (\ref{s2.2})-(\ref{s2.3}) and (\ref{s2.8})-(\ref{s2.9}), respectively.
Let $q$ be the algebraic multiplicity of $\lambda$, $\lambda=\lambda_{k}=\lambda_{k+1}=\cdots =\lambda_{k+q-1}$.
Assume that $t(\lambda)$ and $\mathbb{M}(\lambda)$ are generalized eigenspaces of $\mathbf{T}$ and (\ref{s2.2})-(\ref{s2.3}) associated with $\lambda$,
respectively, and $t^{*}(\lambda)$ and $\mathbb{M}^{*}(\lambda)$ are generalized adjoint eigenspace associated with $\lambda$,
and $t_{h}(\lambda)$ the direct sum of the generalized eigenspace of $\mathbf{T}_{h}$ associated with the eigenvalues $\lambda_{j,h}$
that converge to $\lambda$.\\
\indent For two closed subspaces $\mathfrak{A}$ and $\mathfrak{B}$ of $(L^{2}(\Omega))^{d}$, we define the gap $\widehat{\delta}(\mathfrak{A},\mathfrak{B})$ between $\mathfrak{A}$ and $\mathfrak{B}$ in the sense of norm $\|\cdot\|_{0}$ as
\begin{align*}
 \delta(\mathfrak{A},\mathfrak{B}) = \sup_{\phi\in \mathfrak{A},\|\phi\|_{0}=1}\inf_{\psi\in \mathfrak{B}}\|\phi-\psi\|_{0},\quad
\widehat{\delta}(\mathfrak{A},\mathfrak{B}) = \max\{\delta(\mathfrak{A},\mathfrak{B}),\delta(\mathfrak{B},\mathfrak{A})\}.
\end{align*}
\begin{lemma}\label{lem2-7}
Let $(\lambda,\mathbf{u}, p)$ and $(\lambda_{h},\mathbf{u}_{h}, p_{h})$ be the $j$th eigenpair of (\ref{s2.2})-(\ref{s2.3})
and (\ref{s2.8})-(\ref{s2.9}), and let $(\lambda^{*},\mathbf{u}^{*}, p^{*})$ and $(\lambda^{*}_{h},\mathbf{u}^{*}_{h}, p^{*}_{h})$ be the $j$th eigenpair of (\ref{s2.48})
and (\ref{s2.49}), respectively. Then there holds
\begin{align}\label{s2.65}
&|\lambda_{h}-\lambda|\lesssim |A_{h}(\mathbf{u}-\mathbf{u}_{h},\mathbf{u}^{*}-\mathbf{v}_{h})
+B_{h}(\mathbf{u}^{*}-\mathbf{v}_{h},p-p_{h})
+\overline{B_{h}(\mathbf{u}-\mathbf{u}_{h},p^{*}-q_{h})}
\nonumber\\
&~~~-\lambda(\mathbf{u}-\mathbf{u}_{h},\mathbf{u}^{*}-\mathbf{v}_{h})|+
|(1+\theta)\sum\limits_{F\in\mathcal{E}_{h}}\int_{F}\mu\overline{\{{\nabla\mathbf{u}}^{*}\}}:[\![\underline{\mathbf{u}_{h}}]\!]ds|,\nonumber\\
&~~~~~~~~~for~ [\mathbf{v}_{h},q_{h}]=
[\mathbf{u}_{h}^{*},p_{h}^{*}]~or~[\Pi_{h}\mathbf{u}^{*},\vartheta_{h}p^{*}].
\end{align}
\end{lemma}
\noindent{\bf Proof.}
From (\ref{s2.29}) ($\boldsymbol{f}=\lambda\mathbf{u}$ and $\mathbf{u}^{f}=\mathbf{u}$), (\ref{s2.57})($\boldsymbol{f}=\lambda^{*}\mathbf{u}^{*}$ and $\mathbf{u}^{f*}=\mathbf{u}^{*}$) and (\ref{s2.8})-(\ref{s2.9}), we derive
\begin{align*}
&~~~~A_{h}(\mathbf{u}-\mathbf{u}_{h},\mathbf{u}^{*}-\mathbf{v}_{h})
+B_{h}(\mathbf{u}^{*}-\mathbf{v}_{h},p-p_{h})
+\overline{B_{h}(\mathbf{u}-\mathbf{u}_{h},p^{*}-q_{h})}
\nonumber\\
&~~~-\lambda(\mathbf{u}-\mathbf{u}_{h},\mathbf{u}^{*}-\mathbf{v}_{h})\nonumber\\
&=A_{h}(\mathbf{u},\mathbf{u}^{*})-A_{h}(\mathbf{u},\mathbf{v}_{h})
-A_{h}(\mathbf{u}_{h},\mathbf{u}^{*})
+A_{h}(\mathbf{u}_{h},\mathbf{v}_{h})\nonumber\\
&~~~+B_{h}(\mathbf{u}^{*},p)
-B_{h}(\mathbf{v}_{h},p)-B_{h}(\mathbf{u}^{*},p_{h})+B_{h}(\mathbf{v}_{h},p_{h})\nonumber\\
&~~~
+\overline{B_{h}(\mathbf{u},p^{*})}-\overline{B_{h}(\mathbf{u},q_{h})}
-\overline{B_{h}(\mathbf{u}_{h},p^{*})}
+\overline{B_{h}(\mathbf{u}_{h},q_{h})}\nonumber\\
&~~~-\lambda(\mathbf{u},\mathbf{u}^{*})+\lambda(\mathbf{u},\mathbf{v}_{h})+\lambda(\mathbf{u}_{h},\mathbf{u}^{*})
-\lambda(\mathbf{u}_{h},\mathbf{v}_{h})\nonumber\\
&=0+0-\lambda(\mathbf{u}_{h},\mathbf{v}_{h})+\lambda_{h}(\mathbf{u}_{h},\mathbf{v}_{h})
-(1+\theta)\sum\limits_{F\in\mathcal{E}_{h}}\int_{F}\mu\overline{\{{\nabla\mathbf{u}}^{*}\}}:[\![\underline{\mathbf{u}_{h}}]\!]ds\nonumber\\
&=(\lambda_{h}-\lambda)(\mathbf{u}_{h},\mathbf{v}_{h})
-(1+\theta)\sum\limits_{F\in\mathcal{E}_{h}}\int_{F}\mu\overline{\{{\nabla\mathbf{u}}^{*}\}}:[\![\underline{\mathbf{u}_{h}}]\!]ds.
\end{align*}
Applying the same argument as in Lemma 5 in \cite{Wang2024}, we deduce that $|(\mathbf{u}_{h},\mathbf{v}_{h})|$ has a positive lower bound uniformly with respect to $h$. Therefore,
from the above equation we obtain (\ref{s2.65}).~~~$\Box$\\
\indent Denote $\widehat{\lambda}_{h}=\frac{1}{q}\sum_{j=k}^{k+q-1}\lambda_{j,h}$.\\
\indent
For non-self-adjoint eigenvalue problems, when the algebraic multiplicity $q$ of an eigenvalue is greater than its geometric multiplicity $m$,
the eigenvalue is defective. When $q=m$, the eigenvalue is non-defective. When $\lambda$ is a defective eigenvalue, the error of the single
 approximate eigenpair $(\lambda_{h},\mathbf{u}_{h},p_{h})$ fails to achieve the optimal convergence order. Therefore, in Theorems \ref{th2-9} and \ref{th2-10} below, we provide error
 estimates for the two cases $q\geq m$ and $q=m$ separately.\\
\indent Define the spectral projection $\mathbb{E}$
associated with $\mathbf{T}$ and the eigenvalue $\lambda^{-1}$, and $\mathbb{E}_{h}$
associated with $\mathbf{T}_{h}$ and the eigenvalues of $\mathbf{T}_{h}$ which lie in $\Gamma$
 as follows:
\begin{equation*}
\mathbb{E}=\frac{1}{2\pi \mathrm{i}}\int_{\Gamma}R(\mathbf{T},z)\mathrm{d}z,~~~\mathbb{E}_{h}=\frac{1}{2\pi \mathrm{i}}\int_{\Gamma}R(\mathbf{T}_{h},z)\mathrm{d}z,
\end{equation*}
where $R(\mathbf{T},z)=(z-\mathbf{T})^{-1}$, $R(\mathbf{T}_{h},z)=(z-\mathbf{T}_{h})^{-1}$, and $\Gamma$ is a circle in the complex plane centered at $\lambda^{-1}$ which lies in the resolvent set $\rho(\mathbf{T})$ of $\mathbf{T}$ and
encloses no other eigenvalues of $\mathbf{T}$.
\begin{lemma}\label{lem2-8}
Let $\lambda$ and $\lambda_{h}$ be the $j$th eigenvalue of (\ref{s2.2})-(\ref{s2.3}) and (\ref{s2.8})-(\ref{s2.9}), respectively,
and let $[\mathbf{u}_{h}, p_{h}]$ with $\|\mathbf{u}_{h}\|_{0}=1$ be an eigenfunction corresponding to $\lambda_{h}$, $\mathbf{u}=\frac{\mathbb{E}\mathbf{u}_{h}}{\|\mathbb{E}\mathbf{u}_{h}\|_{0}}$.
Then, when $h$ is small enough, there holds
\begin{align}\label{s2.66}
&\|\mathbf{u}_{h}-\mathbf{u}\|_{0}\lesssim \|(\mathbf{T}-\mathbf{T}_{h})\mathbf{u}\|_{0}.
\end{align}
\end{lemma}
\noindent {\bf Proof.} Using the arguments as in
Theorem 7.1 of \cite{Babuska1991book},  we deduce that
\begin{align*}
&\|\mathbf{u}_{h}-\mathbf{u}\|_{0}\leq 2 \|\mathbb{E}_{h}\mathbf{u}_{h}-\mathbb{E}\mathbf{u}_{h}\|_{0}
=2\|\frac{1}{2\pi \mathrm{i}}\int_{\Gamma}(R(\mathbf{T}_{h},z)-R(\mathbf{T},z))\mathbf{u}_{h}\mathrm{d}z\|_{0}\nonumber\\
&~~~=2\|\frac{1}{2\pi \mathrm{i}}\int_{\Gamma}(R(\mathbf{T},z)(\mathbf{T}-\mathbf{T}_{h})R(\mathbf{T}_{h},z))\mathbf{u}_{h}\mathrm{d}z\|_{0}\nonumber\\
&~~~=2\|\frac{1}{2\pi \mathrm{i}}\int_{\Gamma}(R(\mathbf{T},z)(\mathbf{T}-\mathbf{T}_{h})\frac{1}{z-\lambda_{h}^{-1}}\mathbf{u}_{h}\mathrm{d}z\|_{0}\nonumber\\
&~~~\lesssim\|(\mathbf{T}-\mathbf{T}_{h})\mathbf{u}_{h}\|_{0}\nonumber\\
&~~~\lesssim\|(\mathbf{T}-\mathbf{T}_{h})\mathbf{u}\|_{0}+\|(\mathbf{T}-\mathbf{T}_{h})\|_{0}\|\mathbf{u}_{h}-\mathbf{u}\|_{0}.
\end{align*}
Then from (\ref{s2.64}) we obtain (\ref{s2.66}).~~~$\Box$
\begin{theorem}\label{th2-9}
Assume that $\mathbf{R(\Omega)}$ holds, $\mathbb{M}(\lambda),  \mathbb{M}^{*}(\lambda)\subset (H^{1+s}(\Omega))^{d}\times H^{s}(\Omega)$, $\frac{1}{2}< s\leq k$,
$h$ is small enough, then
 for $\theta=-1$, when $q\geq m$ there hold
\begin{align}\label{s2.67}
&\widehat{\delta}(t(\lambda),t_{h}(\lambda)) \lesssim h^{r+s},\\\label{s2.68}
&|\lambda-\widehat{\lambda}_{h}| \lesssim h^{2s};
\end{align}
 when $q=m$, let $[\mathbf{u}_{h}, p_{h}]$ with $\|\mathbf{u}_{h}\|_{0}=1$ be an eigenfunction corresponding to $\lambda_{h}$,
and let $\mathbf{u}=\frac{\mathbb{E}\mathbf{u}_{h}}{\|\mathbb{E}\mathbf{u}_{h}\|_{0}}$ and $p= \lambda\mathbf{S}\mathbf{u}$, then there hold
\begin{align}\label{s2.69}
&\|\mathbf{u}-\mathbf{u}_{h}\|_{0}\lesssim h^{r}(|||\mathbf{u}-\mathbf{u}_{h}|||+\|p-p_{h}\|_{0}),\\\label{s2.70}
&|||\mathbf{u}-\mathbf{u}_{h}|||+\|p-p_{h}\|_{0}\lesssim h^{s},\\\label{s2.71}
&\|\mathbf{u}^{*}-\mathbf{u}_{h}^{*}\|_{0}\lesssim h^{r}(|||\mathbf{u}^{*}-\mathbf{u}_{h}^{*}|||+\|p^{*}-p_{h}^{*}\|_{0}),\\\label{s2.72}
&|||\mathbf{u}^{*}-\mathbf{u}_{h}^{*}|||+\|p^{*}-p_{h}^{*}\|_{0}\lesssim h^{s},\\\label{s2.73}
&|\lambda-\lambda_{h}|\lesssim h^{2s}.
\end{align}
\end{theorem}
\noindent {\bf Proof.}
  Since $\|\mathbf{T}_{h}-\mathbf{T}\|_{0}\rightarrow 0$ as $h\rightarrow 0$,
from Theorem 7.1 and Theorem 7.2 in \cite{Babuska1991book},
  we obtain
\begin{align}\label{s2.74}
&\widehat{\delta}(t(\lambda),t_{h}(\lambda)) \lesssim \|(\mathbf{T}-\mathbf{T}_{h})|_{t(\lambda)}\|_{0},\\\label{s2.75}
&|\widehat{\lambda}-\widehat{\lambda}_{h}| \lesssim \sum_{i,j=k}^{k+q-1}((\mathbf{T}-\mathbf{T}_{h})\mathbf{u}_{j}, \mathbf{u}_{i}^{*})+\|(\mathbf{T}-\mathbf{T}_{h})|_{t(\lambda)}\|_{0}\|(\mathbf{T}^{*}-\mathbf{T}_{h}^{*})|_{t^{*}(\lambda)}\|_{0}.
\end{align}
From (\ref{s2.74}), (\ref{s2.63}) and (\ref{s2.60}) we obtain (\ref{s2.67}).\\
\indent From (\ref{s2.57}), (\ref{s2.29}) and (\ref{s2.60})-(\ref{s2.61}), we deduce
\begin{align}\label{s2.76}
&|((\mathbf{T}-\mathbf{T}_{h})\mathbf{u}_{j}, \mathbf{u}_{i}^{*})|=|A_{h}((\mathbf{T}-\mathbf{T}_{h})\mathbf{u}_{j}, \mathbf{T}^{*}\mathbf{u}_{j}^{*})\nonumber\\
&~~~+
\overline{B_{h}((\mathbf{T}-\mathbf{T}_{h})\mathbf{u}_{j}, \mathbf{S}^{*}\mathbf{u}_{i}^{*})}
+B_{h}(\mathbf{T}^{*}\mathbf{u}_{i}^{*}, \mathbf{S}^{*}\mathbf{u}_{i}^{*}- \mathbf{S}_{h}^{*}\mathbf{u}_{i}^{*})|\nonumber\\
&~~~=|A_{h}((\mathbf{T}-\mathbf{T}_{h})\mathbf{u}_{j}, \mathbf{T}^{*}\mathbf{u}_{i}^{*}-\mathbf{T}_{h}^{*}\mathbf{u}_{i}^{*})\nonumber\\
&~~~+B_{h}(\mathbf{T}^{*}\mathbf{u}_{i}^{*}-\mathbf{T}_{h}^{*}\mathbf{u}_{i}^{*}, \mathbf{S}^{*}\mathbf{u}_{i}^{*}- \mathbf{S}_{h}^{*}\mathbf{u}_{i}^{*})\nonumber\\
&~~~+\overline{B_{h}((\mathbf{T}-\mathbf{T}_{h})\mathbf{u}_{j}, \mathbf{S}^{*}\mathbf{u}_{i}^{*}-\mathbf{S}_{h}^{*}\mathbf{u}_{i}^{*})}|
\lesssim h^{2s}.
\end{align}
Substituting (\ref{s2.76}) into (\ref{s2.75}), and up to higher order term $\|(\mathbf{T}-\mathbf{T}_{h})|_{t(\lambda)}\|_{0}\|(\mathbf{T}^{*}-\mathbf{T}_{h}^{*})|_{t^{*}(\lambda)}\|_{0}$ we obtain (\ref{s2.68}).\\
\indent
When $\lambda$ is a non-defective eigenvalue,
using (\ref{s2.66}) and (\ref{s2.63}),
 we derive
\begin{align}\label{s2.77}
\|\mathbf{u}-\mathbf{u}_{h}\|_{0}&\lesssim \|\mathbf{T}\mathbf{u}-\mathbf{T}_{h}\mathbf{u}\|_{0}\nonumber\\
&\lesssim h^{r}(|||\mathbf{T}\mathbf{u}-\mathbf{T}_{h}\mathbf{u}|||+\|\mathbf{S}\mathbf{u}-\mathbf{S}_{h}\mathbf{u}\|_{0}).
\end{align}
From (\ref{s2.65}) with $[\mathbf{v}_{h}, q_{h}]=[\Pi_{h}\mathbf{u}^{*},\vartheta_{h}p^{*}]$
, we derive
\begin{align}\label{s2.78}
|\lambda-\lambda_{h}|\lesssim h^{r}(|||\mathbf{u}-\mathbf{u}_{h}|||
+\|p-p_{h}\|_{0})
\end{align}
Applying the triangle inequality and (\ref{s2.45}), we obtain
\begin{align}\label{s2.79}
&|~|||\mathbf{u}-\mathbf{u}_{h}|||-|||\lambda \mathbf{T}\mathbf{u}-\lambda \mathbf{T}_{h}\mathbf{u}|||~|
\leq |||\lambda \mathbf{T} \mathbf{u}-\lambda_{h} \mathbf{T}_{h} \mathbf{u}_{h}-\lambda \mathbf{T} \mathbf{u}+\lambda \mathbf{T}_{h} \mathbf{u}||| \nonumber\\
&~~~\leq \|\lambda \mathbf{T}_{h}\mathbf{u}-\lambda_{h}\mathbf{T}_{h}\mathbf{u}_{h}\|_{h}
\lesssim \|\lambda\mathbf{u}-\lambda_{h}\mathbf{u}_{h}\|_{0},\\\label{s2.80}
&|~||p-p_{h}||_{0}-||\lambda \mathbf{S}\mathbf{u}-\lambda \mathbf{S}_{h}\mathbf{u}||_{0}~|
\leq \|\lambda \mathbf{S} \mathbf{u}-\lambda_{h} \mathbf{S}_{h} \mathbf{u}_{h}-\lambda \mathbf{S} \mathbf{u}+\lambda \mathbf{S}_{h} \mathbf{u}\|_{0}\nonumber\\
&~~~\leq \|\lambda \mathbf{S}_{h}\mathbf{u}-\lambda_{h}\mathbf{S}_{h}\mathbf{u}_{h}\|_{0}
\lesssim \|\lambda\mathbf{u}-\lambda_{h}\mathbf{u}_{h}\|_{0}.
\end{align}
From the above two estimates, (\ref{s2.77}) and (\ref{s2.78}), we obtain
\begin{align}\label{s2.81}
|||\mathbf{u}-\mathbf{u}_{h}|||+||p-p_{h}||_{0}\simeq |||\lambda \mathbf{T}\mathbf{u}-\lambda \mathbf{T}_{h}\mathbf{u}|||
+||\lambda \mathbf{S}\mathbf{u}-\lambda \mathbf{S}_{h}\mathbf{u}||_{0}.
\end{align}
From (\ref{s2.81}) and (\ref{s2.77}), we obtain (\ref{s2.69}).
From (\ref{s2.81}) and (\ref{s2.60}),  we obtain (\ref{s2.70}).
By the same argument as for (\ref{s2.69}) and (\ref{s2.70}), we can obtain  (\ref{s2.71}) and (\ref{s2.72}).
Combining (\ref{s2.65}) with (\ref{s2.69}), (\ref{s2.70}), (\ref{s2.71}) and (\ref{s2.72}), we derive (\ref{s2.73}).~~~$\Box$
\begin{theorem}\label{th2-10}
Under the conditions of Theorem \ref{th2-9},
for any $\theta\in \{0,1\}$, when $q\geq m$ there hold
\begin{align}\label{s2.82}
&\widehat{\delta}(t(\lambda),t_{h}(\lambda)) \lesssim h^{s},\\\label{s2.83}
&|\lambda-\widehat{\lambda}_{h}| \lesssim h^{s};
\end{align}
when $q= m$ there hold
\begin{align}\label{s2.84}
&|||\mathbf{u}-\mathbf{u}_{h}|||+\|p-p_{h}\|_{0}\lesssim h^{s},\\\label{s2.85}
&|\lambda-\lambda_{h}| \lesssim h^{s}.
\end{align}
\end{theorem}
\noindent {\bf Proof.} From Theorems 7.1 and 7.2 in \cite{Babuska1991book}, and from
(\ref{s2.60}) we obtain (\ref{s2.82}) and (\ref{s2.83}).
From (\ref{s2.79}), (\ref{s2.80})
and (\ref{s2.60}) we obtain (\ref{s2.84}).
From (\ref{s2.84})
and (\ref{s2.65}) we obtain (\ref{s2.85}).
~~~$\Box$\\

\noindent{\bf Remark 2.1.}  By the proof method of Theorem 3.1 and Theorem 4.1 in \cite{Badia2014}, we can deduce that if $[\mathbf{u}^{f}, p^{f}]$ and $[\mathbf{u}^{f*}, p^{f*}]\in H^{1+s}(\Omega)^{d}\times H^{s}(\Omega)$, $0<s\leq k$, then
\begin{align*}
&\|\mathbf{u}^{f}-\mathbf{u}^{f}_{h}\|_{h}+\|p^{f}-p^{f}_{h}\|_{0}\lesssim h^{s}(|\mathbf{u}^{f}|_{1+s}+|p^{f}|_{s})
+h\|\Pi_{h}\boldsymbol{f}-\boldsymbol{f}\|_{0}~~(\theta\in \{-1,0,1\}),\\
&\|\mathbf{u}^{f*}-\mathbf{u}^{f*}_{h}\|_{h}+\|p^{f*}-p^{f*}_{h}\|_{0}\lesssim h^{s}(|\mathbf{u}^{f*}|_{1+s}+|p^{f*}|_{s})
+h\|\Pi_{h}\boldsymbol{f}-\boldsymbol{f}\|_{0}~~(\theta=-1).
\end{align*}
We replace (\ref{s2.60})-(\ref{s2.61}) in Theorem \ref{th2-5} with the above two estimates. Then, under the condition that the regularity result $\mathbf{R(\Omega)}$ holds for some $r\in (0,1]$,
 we deduce that all results in this paper hold after replacing $h^{r+s}$ with $h^{s}$ on the right-hand side of (\ref{s2.67}), apart from Theorem 2.6 and the estimates (\ref{s2.69}) and (\ref{s2.71}) in Theorem \ref{th2-9}.
\section {A posteriori error estimate }
\indent For the primal problem (\ref{s2.8})-(\ref{s2.9}) and the adjoint problem (\ref{s2.49}),
we introduce the estimators as follows.\\
\indent Let $\lambda$ be the $j$th eigenvalue of (\ref{s2.1}) with multiplicity $q=m\geq 1$, $(\mathbf{u},p)$  be an eigenfunction corresponding to $\lambda$,
$(\lambda_{h}, \mathbf{u}_{h}, p_{h})$ be the $j$th eigenpair of (\ref{s2.8})-(\ref{s2.9}). We introduce the estimator for the problem
\begin{align*}
&\eta^{2}_{R_{\tau}}=h_{\tau}^{2}\|\lambda_{h}\mathbf{u}_{h}+\mu\Delta\mathbf{u}_{h}-(\boldsymbol{\beta}\cdot\nabla)\mathbf{u}_{h}-\nabla p_{h}\|^{2}_{0,\tau}+\|div\mathbf{u}_{h}\|_{0,\tau}^{2},\\
&\eta^{2}_{F_{\tau}}=\frac{1}{2}\sum\limits_{F\subset\partial\tau\setminus\partial\Omega}h_{F}\|[\![(p_{h}\mathbf{I}-\mu \nabla\mathbf{u}_{h}+\mathbf{u_{h}}\otimes\boldsymbol{\beta})]\!]\|^{2}_{0,F},
\end{align*}
where $\mathbf{I}$ denotes the $d\times d$ $(d=2, 3)$ identity matrix. Next, we introduce the following estimator $\eta_{J_{\tau}}$ to measure the jump of the approximate solution $\mathbf{u}_{h}$:
\begin{align*}
\eta^{2}_{J_{\tau}}=\sum\limits_{F\subset\partial\tau, F\in\mathcal{E}^{i}_{h}} \gamma h_{F}^{-1}||[\![\underline{\mathbf{u}_{h}}]\!]||_{0,F}^{2}+ \sum\limits_{F\subset\partial\tau, F\in\mathcal{E}^{b}_{h}}\gamma h_{F}^{-1}||\mathbf{u}_{h}\otimes\mathbf{n}||_{0,F}^{2}.
\end{align*}
The local error estimator is defined as
\begin{align*}
\eta^{2}_{\tau}(\mathbf{u}_{h},p_{h})=\eta^{2}_{R_{\tau}}+\eta^{2}_{F_{\tau}}+\eta^{2}_{J_{\tau}},
\end{align*}
then the global a posteriori error estimator is defined as
\begin{align*}
\eta_{h}(\mathbf{u}_{h},p_{h})=(\sum\limits_{\tau\in\mathcal{T}_{h}}
\eta^{2}_{\tau}(\mathbf{u}_{h},p_{h}))^{\frac{1}{2}}.
\end{align*}
\indent Let $(\lambda^{*}_{h}, \mathbf{u}^{*}_{h}, p^{*}_{h})\in \mathbb{C}\times\mathbb{X}_{h}\times \mathbb{W}_{h}$ be an approximate eigenpair
of the adjoint problem (\ref{s2.49}), we introduce the  estimator  as follows
\begin{align*}
&\eta^{*2}_{R_{\tau}}=h_{\tau}^{2}\|\lambda^{*}_{h}\mathbf{u}^{*}_{h}+\mu\Delta\mathbf{u}^{*}_{h}+(\boldsymbol{\beta}\cdot\nabla)\mathbf{u}^{*}_{h}
-\nabla p^{*}_{h}\|^{2}_{0,\tau}+\|div\mathbf{u}^{*}_{h}\|_{0,\tau}^{2},\\
&{\eta}^{*2}_{F_{\tau}}= \frac{1}{2}\sum\limits_{F\subset\partial\tau\setminus\partial\Omega}h_{F}\|[\![{p}^{*}_{h}\mathbf{I}-\mu\nabla \mathbf{{u}}^{*}_{h}-\mathbf{{u}}^{*}_{h}\otimes \boldsymbol{\beta}]\!]\|^{2}_{0,F},\\
&\eta^{*2}_{J_{\tau}}=\sum\limits_{F\subset\partial\tau, F\in\mathcal{E}^{i}_{h}}\gamma h_{F}^{-1}||[\![\underline{\mathbf{u}^{*}_{h}}]\!]||_{0,F}^{2}+ \sum\limits_{F\subset\partial\tau, F\in\mathcal{E}^{b}_{h}} \gamma h_{F}^{-1}||\mathbf{u}^{*}_{h}\otimes\mathbf{n}||_{0,F}^{2}.
\end{align*}
The local error estimator and global estimator of problem (\ref{s2.49}) are defined as follows, respectively:
\begin{align*}
\eta^{*2}_{\tau}(\mathbf{u}^{*}_{h},p^{*}_{h})=\eta^{*2}_{R_{\tau}}+\eta^{*2}_{F_{\tau}}+\eta^{*2}_{J_{\tau}},~~~
\eta^{*}_{h}(\mathbf{u}^{*}_{h},p^{*}_{h})=(\sum\limits_{\tau\in\mathcal{T}_{h}}\eta^{*2}_{\tau}(\mathbf{u}^{*}_{h},p^{*}_{h}))^{\frac{1}{2}}.
\end{align*}
\indent We denote $\omega(\tau)=$int$\{\bigcup\limits_{\overline{\tau}_{i}\cap \overline{\tau}\not=\emptyset}\bar{\tau}_{i}, \tau_{i}\in\mathcal{T}_{h}\}$ for $\tau\in\mathcal{T}_{h}$, and use $\omega(F)$ to represent the set of all elements which share at least one node with face $F$.
We denote by $\mathbf{z}^{I}$ the Scott-Zhang interpolation function \cite{ScottZhang1990}, then $\mathbf{z}^{I}\in \mathbb{X}\cap\mathbb{X}_{h}$ and
\begin{eqnarray}
&&\|\mathbf{z}-\mathbf{z}^{I}\|_{0,\tau}+h_{\tau}\|\nabla(\mathbf{z}-\mathbf{z}^{I})\|_{0,\tau}\lesssim h_{\tau}\|\nabla\mathbf{z}\|_{0,\omega(\tau)},~~~~~~~\forall \tau\in\mathcal{T}_{h},\label{s3.1}\\
&&\|\mathbf{z}-\mathbf{z}^{I}\|_{0,F} \lesssim h_{F}^{\frac{1}{2}}\|\nabla\mathbf{z}\|_{0,\omega(F)},~~~~~\forall F\subset\partial\tau.\label{s3.2}
\end{eqnarray}
\indent
We define an auxiliary sesquilinear form $\widetilde{A_{h}}(\cdot,\cdot):\mathbb{X}(h)\times \mathbb{X}(h)\longrightarrow \mathbb{C}$ by
\begin{align}\label{s3.3}
\widetilde{\mathcal{A}_{h}}([\mathbf{w},r], [\mathbf{v},q])=\widetilde{A}_{h}(\mathbf{w},\mathbf{v})+B_{h}'(\mathbf{v},r)
+\overline{B_{h}'(\mathbf{w},q)},~~~ \forall [\mathbf{w},r], [\mathbf{v},q]\in \mathbb{X}(h)\times \mathbb{W},
\end{align}
where
\begin{align}\label{s3.4}
&\widetilde{A}_{h}(\mathbf{u}_{h},\mathbf{v}_{h})=\sum\limits_{\tau\in\mathcal{T}_{h}}\int_{\tau}\mu\nabla\mathbf{u}_{h}:\overline{\nabla\mathbf{v}}_{h}dx
+ \sum\limits_{\tau\in\mathcal{T}_{h}}\int_{\tau}(\boldsymbol{\beta}\cdot\nabla)\mathbf{u}_{h}\cdot \overline{\mathbf{v}}_{h}dx\nonumber\\
&~~~
+\sum\limits_{F\in\mathcal{E}_{h}}\int_{F}\gamma h_{F}^{-1}[\![\underline{\mathbf{u}_{h}}]\!]:
\overline{[\![{\underline{\mathbf{v}}_{h}}]\!]}ds
-\frac{1}{2}\sum\limits_{\tau\in\mathcal{T}_{h}}\sum_{F\in\partial\tau}\int_{F} [\![\mathbf{u}_{h}\otimes\boldsymbol{\beta}]\!]\cdot\overline{\{ \mathbf{v}_{h}\}} ds,\\\label{s3.5}
&B_{h}'(\mathbf{v}_{h},q_{h})=
-\sum\limits_{\tau\in\mathcal{T}_{h}}\int_{\tau}q_{h}\overline{\nabla\cdot\mathbf{v}_{h}}dx.
\end{align}
With the definition, one finds that $\widetilde{\mathcal{A}_{h}}(\cdot,\cdot)=\mathcal{A}(\cdot,\cdot)$
on $\mathbb{V}\times \mathbb{V}$,  $\widetilde{A}_{h}(\cdot,\cdot)=A(\cdot,\cdot)$ on $\mathbb{X}\times \mathbb{X}$
and $B_{h}'(\cdot,\cdot)=B(\cdot,\cdot)$ on $\mathbb{X}\times \mathbb{W}$,
and that
\begin{align}\label{s3.6}
&|\widetilde{A}_{h}(\mathbf{u}_{h},\mathbf{v}_{h})|\lesssim \|\mathbf{u}_{h}\|_{h}\|\mathbf{v}_{h}\|_{h},~~~
\forall [\mathbf{u}_{h},\mathbf{v}_{h}]\in \mathbb{X}(h)\times\mathbb{X}(h),\\\label{s3.7}
&|B_{h}'(\mathbf{v}_{h},q_{h})|\lesssim
\|\mathbf{v}_{h}\|_{h}\|q_{h}\|_{0},~~~\forall [\mathbf{v}_{h},q_{h}]\in \mathbb{X}(h)\times \mathbb{W}.
\end{align}
\begin{theorem}\label{th3-1}
Let $(\lambda_{h},\mathbf{u}_{h}, p_{h})$  and $(\lambda^{*}_{h},\mathbf{u}^{*}_{h}, p^{*}_{h})$ be
   the $j$th eigenpairs of (\ref{s2.8})-(\ref{s2.9})
and (\ref{s2.49}), respectively.
Let $(\lambda,\mathbf{u}, p)$ and $(\lambda^{*},\mathbf{u}^{*}, p^{*})$
 be the $j$th eigenpair of (\ref{s2.2})-(\ref{s2.3})
and (\ref{s2.48}), respectively.  Then for any $\theta\in \{-1,0,1\}$, there hold
\begin{align}
&~~~\|\mathbf{u}-\mathbf{u}_{h}\|_{h}+\| p-p_{h}\|_{0}\label{s3.8}\\
&\lesssim\sup\limits_{\mathbf{0}
\neq [\mathbf{v},q]\in\mathbb{V}}\frac{\left|\widetilde{A}_{h}(\mathbf{u}-\mathbf{u}_{h}, \mathbf{v})
+B'_{h}(\mathbf{v},p- p_{h})+\overline{B'_{h}(\mathbf{u}-\mathbf{u}_{h},q)}\right|}{\|[\mathbf{v},q]\|_{\mathbb{V}}}
+\inf\limits_{\mathbf{v}\in \mathbb{X}}\|\mathbf{u}_{h}-\mathbf{v}\|_{h}, \nonumber\\
&~~~~~~~~~~~~\nonumber\\
&~~~\|\mathbf{u}^{*}-\mathbf{u}^{*}_{h}\|_{h}+\| p^{*}-p^{*}_{h}\|_{0}\label{s3.9}\\
&\lesssim\sup\limits_{\mathbf{0}\neq [\mathbf{v},q]\in\mathbb{V}}
\frac{\left|\widetilde{A}_{h}(\mathbf{v},\mathbf{u}^{*}-\mathbf{u}_{h}^{*})
+\overline{B'_{h}(\mathbf{v},p^{*}- p_{h}^{*})}+B'_{h}(\mathbf{u}^{*}-\mathbf{u}_{h}^{*},q)\right|}{\|[\mathbf{v},q]\|_{\mathbb{V}}}
+\inf\limits_{\mathbf{v}\in \mathbb{X}}\|\mathbf{u}^{*}_{h}-\mathbf{v}\|_{h}. \nonumber
\end{align}
\end{theorem}
\noindent{\bf Proof.}
For any $[\mathbf{w},z]\in \mathbb{V}$, due to the fact that $\widetilde{\mathcal{A}_{h}}=\mathcal{A}$ on $\mathbb{V}\times \mathbb{V}$,
From the inf-sup condition (\ref{s2.7}) we obtain
\begin{align}\label{s3.10}
 \alpha\|[\mathbf{u}-\mathbf{w},p-z]\|_{\mathbb{V}}\leq\sup\limits_{[\mathbf{v},q]\in \mathbb{V}}\frac{|\widetilde{\mathcal{A}_{h}}
([\mathbf{u}-\mathbf{w},p-z], [\mathbf{v},q])|}{\|[\mathbf{v},q]\|_{\mathbb{V}}}.
\end{align}
From (\ref{s3.6}) and (\ref{s3.7}), we obtain
\begin{align}\label{s3.11}
&|\widetilde{\mathcal{A}_{h}}([\mathbf{u}-\mathbf{w},p-z], [\mathbf{v},q])|=|\widetilde{A}_{h}(\mathbf{u}-\mathbf{w},\mathbf{v})+B'_{h}(\mathbf{v},p-z)
+B'_{h}(\mathbf{u}-\mathbf{w},q)|\nonumber\\
&~~~=|\widetilde{A}_{h}(\mathbf{u}-\mathbf{u}_{h},\mathbf{v})+B'_{h}(\mathbf{v},p-p_{h})
+B'_{h}(\mathbf{u}-\mathbf{u}_{h},q)\nonumber\\
&~~~+\widetilde{A}_{h}(\mathbf{u}_{h}-\mathbf{w},\mathbf{v})+B'_{h}(\mathbf{v},p_{h}-z)
+B'_{h}(\mathbf{u}_{h}-\mathbf{w},q)|\nonumber\\
&~~~\lesssim|\widetilde{A}_{h}(\mathbf{u}-\mathbf{u}_{h},\mathbf{v})+B'_{h}(\mathbf{v},p-p_{h})
+B'_{h}(\mathbf{u}-\mathbf{u}_{h},q)|\\
&~~~+\|\mathbf{u}_{h}-\mathbf{w}\|_{h}\|\mathbf{v}\|_{1}+\|\mathbf{v}\|_{1}\|p_{h}-z\|_{0}
+\|\mathbf{u}_{h}-\mathbf{w}\|_{h}\|q\|_{0},
~~~ \forall [\mathbf{w},z], [\mathbf{v},q]\in \mathbb{V}.\nonumber
\end{align}
Substituting (\ref{s3.11}) into (\ref{s3.10}),
using the triangle inequality, noting that $\mathbf{w}$ and $ z$ are arbitrary and $\inf\limits_{z\in \mathbb{W}}\|z- p_{h}\|_{0}=0$, then (\ref{s3.8}) is valid.
 By the same argument as for (\ref{s3.8}), then (\ref{s3.9}) holds.~~~~$\Box$
\begin{lemma}\label{lem3-2}
Under the conditions of Theorem \ref{th3-1}, there hold
\begin{align}
&~~~|\widetilde{A}_{h}(\mathbf{u}-\mathbf{u}_{h}, \mathbf{v})+B'_{h}(\mathbf{v}, p-p_{h})|
+|B'_{h}(\mathbf{u_{h}}, q) |\nonumber\\
&\lesssim  \left(\sum\limits_{\tau\in\mathcal{T}_{h}}(\eta_{R_{\tau}}^2+\eta_{F_{\tau}}^2 +|\theta|\eta_{J_{\tau}}^2 )\right)^{\frac{1}{2}}\|\nabla\mathbf{v}\|_{0}+(\sum\limits_{\tau\in\mathcal{T}_{h}}\|div\mathbf{u}_{h}\|_{0,\tau}^2)^{\frac{1}{2}}\|q\|_{0}\nonumber\\
&~~~+\|\lambda \mathbf{u}-\lambda_{h} \mathbf{u}_{h}\|_{0}\|\nabla\mathbf{v}\|_{0},
~~~\forall [\mathbf{v},q]\in \mathbb{V},\label{s3.12}\\
&~~~|\widetilde{A}_{h}(\mathbf{v},\mathbf{u}^{*}-\mathbf{u}^{*}_{h})+B'_{h}(\mathbf{v},p^{*}-p^{*}_{h})|
+|B'_{h}(\mathbf{u}_{h}^{*}, q) |\nonumber\\
&\lesssim \left(\sum\limits_{\tau\in\mathcal{T}_{h}}(\eta^{*2}_{R_{\tau}}+\eta^{*2}_{F_{\tau}} +\eta^{*2}_{J_{\tau}})\right)^{\frac{1}{2}}\|\nabla\mathbf{v}\|_{0}
+(\sum\limits_{\tau\in\mathcal{T}_{h}}\|div\mathbf{u}_{h}^{*}\|_{0,\tau}^2)^{\frac{1}{2}}\|q\|_{0}\nonumber\\
&~~~+\|\lambda^{*} \mathbf{u}^{*}-\lambda^{*}_{h} \mathbf{u}^{*}_{h}\|_{0}\|\nabla\mathbf{v}\|_{0},
~~~\forall [\mathbf{v},q]\in \mathbb{V}.\label{s3.13}
\end{align}
\end{lemma}
\noindent {\bf Proof.}
 Using (\ref{s2.2}), (\ref{s3.4}), (\ref{s3.5}), (\ref{s2.11}) and Green's formula we deduce
\begin{align*}
&~~~~\widetilde{A}_{h}(\mathbf{u}-\mathbf{u}_{h},\mathbf{v})
+B'_{h}(\mathbf{v},p-p_{h})\nonumber\\
&=\widetilde{A}_{h}(\mathbf{u},\mathbf{v})-\widetilde{A}_{h}(\mathbf{u}_{h},\mathbf{v})
+B'_{h}(\mathbf{v},p)-B'_{h}(\mathbf{v},p_{h})\nonumber\\
&=\lambda (\mathbf{u},\mathbf{v})-\widetilde{A}_{h}(\mathbf{u}_{h},\mathbf{v})-B'_{h}(\mathbf{v},p_{h})\nonumber\\
&=(\lambda\mathbf{u}-\lambda_{h}\mathbf{u}_{h},\mathbf{v})
+\sum\limits_{\tau\in\mathcal{T}_{h}}\int_{\tau}(\lambda_{h}\mathbf{u}_{h}+\mu\Delta\mathbf{u}_{h})\cdot\overline{\mathbf{v}}dx
-\sum\limits_{\tau\in\mathcal{T}_{h}}\sum\limits_{F\subset\partial\tau}\int_{F}\mu \nabla\mathbf{u}_{h}\mathbf{n}  \cdot\overline{\mathbf{v}}ds\nonumber\\
&~~~ -\sum\limits_{\tau\in\mathcal{T}_{h}}\int_{\tau}(\boldsymbol{\beta}\cdot\nabla)\mathbf{u}_{h}\cdot\overline{\mathbf{v}}dx
-\sum\limits_{\tau\in\mathcal{T}_{h}}\int_{\tau}
\nabla {p_{h}}\cdot\overline{\mathbf{v}}dx
+\sum\limits_{\tau\in\mathcal{T}_{h}}\sum\limits_{F\subset\partial\tau}\int_{F}p_{_h}\overline{\mathbf{v}}\cdot\mathbf{n}ds\nonumber\\
&~~~+\frac{1}{2}\sum\limits_{\tau\in\mathcal{T}_{h}}\sum\limits_{F\subset\partial\tau}\int_{F}[\![\mathbf{u}_{h}\otimes\boldsymbol{\beta}]\!]\cdot\overline{\{ \mathbf{v}\}} ds
\end{align*}
\begin{align}
&=(\lambda\mathbf{u}-\lambda_{h}\mathbf{u}_{h},\mathbf{v})
+\sum\limits_{\tau\in\mathcal{T}_{h}}\int_{\tau}(\lambda_{h}\mathbf{u}_{h}+\mu\Delta\mathbf{u}_{h}
-(\boldsymbol{\beta}\cdot\nabla)\mathbf{u}_{h}-\nabla p_{h})\cdot\overline{\mathbf{v}}dx\nonumber\\
&~~~+\frac{1}{2}\sum\limits_{\tau\in\mathcal{T}_{h}}\sum\limits_{F\subset\partial\tau}
\int_{F}[\![p_{_h}\mathbf{I}-\mu\nabla\mathbf{u}_{h}+\mathbf{u}_{h}\otimes\boldsymbol{\beta}]\!]
\cdot\overline{\{\mathbf{v}\}}ds,~~~\forall \mathbf{v}\in \mathbb{X}.\label{s3.14}
\end{align}
Note that the Scott-Zhang interpolation $\mathbf{v}-\mathbf{v}^{I}\in \mathbb{X}$.
From (\ref{s3.14}) replacing $\mathbf{v}$ with $\mathbf{v}-\mathbf{v}^{I}$,
the Cauchy-Schwarz inequality and (\ref{s3.1}) and (\ref{s3.2}), we deduce that
\begin{align}\label{s3.15}
&~~~~|\widetilde{A}_{h}(\mathbf{u}-\mathbf{u}_{h},\mathbf{v}-\mathbf{v}^{I})
+B'_{h}((\mathbf{v}-\mathbf{v}^{I}),p-p_{h})|\nonumber\\
&\lesssim (\sum\limits_{\tau\in\mathcal{T}_{h}}(\eta_{R_{\tau}}^{2}+\eta_{F_{\tau}}^{2}))^{\frac{1}{2}}|\mathbf{v}|_{1}
+h\|\lambda\mathbf{u}-\lambda_{h}\mathbf{u}_{h}\|_{0}|\mathbf{v}|_{1}.
\end{align}
By comparing (\ref{s2.10}) with (\ref{s3.4}), using (\ref{s2.2}),(\ref{s2.8}) and $\mathbf{v}^{I}\in \mathbb{X}\cap \mathbb{X}_{h}$, applying the Cauchy-Schwarz inequality, the trace inequality, and (\ref{s3.15}), we obtain
\begin{align*}
&~~~~|\widetilde{A}_{h}(\mathbf{u}-\mathbf{u}_{h}, \mathbf{v}^{I})+B'_{h}(\mathbf{v}^{I}, p-p_{h})|\nonumber\\
&=|A_{h}(\mathbf{u}-\mathbf{u}_{h}, \mathbf{v}^{I})
-\theta\sum_{F\in\mathcal{E}_{h}}\mu\int_{F} [\![ \underline{\mathbf{u}-\mathbf{u}_{h}}]\!]: \overline{\{ \nabla\mathbf{v}^{I} \}} ds
+B'_{h}(\mathbf{v}^{I}, p-p_{h})|\nonumber\\
&=|-\theta\sum_{F\in\mathcal{E}_{h}}\int_{F}\mu [\![ \underline{\mathbf{u}-\mathbf{u}_{h}}]\!]: \overline{\{ \nabla\mathbf{v}^{I} \}} ds
+(\lambda \mathbf{u}-\lambda_{h} \mathbf{u}_{h}, \mathbf{v}^{I})|\nonumber\\
&\lesssim
|\theta|\sum\limits_{\tau\in \mathcal{T}_{h}}\sum_{F\subset\partial\tau\backslash \partial\Omega}\mu h_{F}^{-\frac{1}{2}}\|[\![\underline{\mathbf{u}_{h}}]\!]\|_{0,F}\|\overline{\nabla\mathbf{v}^{I}}\|_{0,\tau}
+\|\lambda \mathbf{u}-\lambda_{h} \mathbf{u}_{h}\|_{0} \|\mathbf{v}^{I}\|_{0}.
\end{align*}
Then,
\begin{align}\label{s3.16}
&~~~~|\widetilde{A}_{h}(\mathbf{u}-\mathbf{u}_{h}, \mathbf{v})+B'_{h}(\mathbf{v}, p-p_{h})|\nonumber\\
&=|\widetilde{A}_{h}(\mathbf{u}-\mathbf{u}_{h}, \mathbf{v}-\mathbf{v}^{I})+B'_{h}(\mathbf{v}-\mathbf{v}^{I}, p-p_{h})
+\widetilde{A}_{h}(\mathbf{u}-\mathbf{u}_{h}, \mathbf{v}^{I})+B'_{h}(\mathbf{v}^{I}, p-p_{h})|\nonumber\\
&\lesssim|\widetilde{A}_{h}(\mathbf{u}-\mathbf{u}_{h}, \mathbf{v}-\mathbf{v}^{I})+B'_{h}(\mathbf{v}-\mathbf{v}^{I}, p-p_{h})|\nonumber\\
&~~~+|\theta|
\sum\limits_{\tau\in \mathcal{T}_{h}}\sum_{F\subset\partial\tau\backslash \partial\Omega}h_{F}^{-\frac{1}{2}}\|[\![\underline{\mathbf{u}_{h}}]\!]\|_{0,F}\|\overline{\nabla\mathbf{v}^{I}}\|_{0,\tau}
+\|\lambda \mathbf{u}-\lambda_{h} \mathbf{u}_{h}\|_{0} \|\mathbf{v}^{I}\|_{0}.
\end{align}
Substituting (\ref{s3.15}) into (\ref{s3.16}), noting that $|B'_{h}(\mathbf{u_{h}}, q)|\lesssim \sum\limits_{\tau\in \mathcal{T}_{h}}\|div\mathbf{u}_{h}\|_{0,\tau}\|q\|_{0} $,
 we obtain (\ref{s3.12}). Using an argument similar to that for (\ref{s3.12}), we obtain (\ref{s3.13}).~~~~$\Box$\\
\indent By virtue of the enriching operator $E_{h}:\mathbb{X}_{h}\to \mathbb{X}_{h}\cap\mathbb{X}$ \cite{Brenner2003,Ohannes2003}, we have the following lemma.
\begin{lemma}\label{lem3-3}
~The following estimates are valid:
\begin{eqnarray}
 &\|\mathbf{u}_{h}-E_{h}\mathbf{u}_{h}\|_{h}^{2}\lesssim \gamma\sum\limits_{F\in\mathcal{E}_{h}^{i}} h_{F}^{-1}||[\![\underline{\mathbf{u}_{h}}]\!]||_{0,F}^{2}+\gamma \sum\limits_{F\in\mathcal{E}_{h}^{b}}
 h_{F}^{-1}||\mathbf{u}_{h}\otimes\mathbf{n}||_{0,F}^{2},\label{s3.17}\\
 &\|\mathbf{u}^{*}_{h}-E_{h}\mathbf{u}^{*}_{h}\|_{h}^{2}\lesssim \gamma\sum\limits_{F\in\mathcal{E}_{h}^{i}} h_{F}^{-1}||[\![\underline{\mathbf{u}^{*}_{h}}]\!]||_{0,F}^{2}+ \gamma\sum\limits_{F\in\mathcal{E}_{h}^{b}} h_{F}^{-1}||\mathbf{u}^{*}_{h}\otimes\mathbf{n}||_{0,F}^{2}.\label{s3.18}
\end{eqnarray}
\end{lemma}
\begin{theorem}\label{th3-4}
 Suppose that the conditions of Theorem \ref{th3-1} hold, then
\begin{align}
&\|\mathbf{u}-\mathbf{u}_{h}\|_{h}+\|p-p_{h}\|_{0}
\lesssim\eta_{h}(\mathbf{u}_{h},p_{h})+\|\lambda_{h}\mathbf{u}_{h}-\lambda \mathbf{u}\|_{0},\label{s3.19}\\
&\|\mathbf{u}^{*}-\mathbf{u}^{*}_{h}\|_{h}+\|p^{*}-p^{*}_{h}\|_{0}
\lesssim\eta^{*}_{h}(\mathbf{u}_{h}^{*},p_{h}^{*})+\|\lambda^{*}_{h}\mathbf{u}^{*}_{h}-\lambda^{*} \mathbf{u}^{*}\|_{0}.\label{s3.20}
\end{align}
\end{theorem}
\noindent{\bf Proof.}
Substituting (\ref{s3.12}) and (\ref{s3.17}) into (\ref{s3.8}), we obtain (\ref{s3.19}). Similarly, (\ref{s3.20}) can be obtained.   ~~~  $\Box$\\
\indent Using the standard bubble argument developed by Verf\"{u}rth \cite{Verfurth2013}, we can deduce the following local bounds.
\begin{lemma}\label{lem3-5}
Under the conditions of Theorem \ref{th3-1}, there hold
\begin{align}\label{s3.21}
&\eta_{R_{\tau}}\lesssim
 \|\nabla (\mathbf{u}-\mathbf{u}_{h})\|_{0,\tau}+\|p-p_{h}\|_{0,\tau}
+h_{\tau}\|\lambda_{h}\mathbf{u}_{h}-\lambda \mathbf{u}\|_{0,\tau},\\\label{s3.22}
&\eta_{F_{\tau}}
\lesssim \|\nabla (\mathbf{u}-\mathbf{u}_{h})\|_{0,\omega(\tau)}+\|p-p_{h}\|_{0,\omega(\tau)}
+( \sum\limits_{\tau\in\omega(\tau)}h_{\tau}^{2} \|\lambda\mathbf{u}-\lambda_{h}\mathbf{u}_{h}\|_{0,\tau}^{2} )^{\frac{1}{2}},\\\label{s3.23}
&\eta^{2}_{J_{\tau}}=\sum\limits_{F\subset\partial\tau, F\in\mathcal{E}^{i}_{h}}\gamma h_{F}^{-1}|[\![{\mathbf{u}_{h}-\mathbf{u}}]\!]|_{0,F}^{2}+ \sum\limits_{F\subset\partial\tau, F\in\mathcal{E}^{b}_{h}}\gamma h_{F}^{-1}||(\mathbf{u}_{h}-\mathbf{u})\otimes\mathbf{n}||_{0,F}^{2},\\\label{s3.24}
&\eta^{*}_{R_{\tau}}\lesssim
 \|\nabla (\mathbf{u}^{*}-\mathbf{u}^{*}_{h})\|_{0,\tau}+\|p^{*}-p^{*}_{h}\|_{0,\tau}
+h_{\tau}\|\lambda^{*}_{h}\mathbf{u}^{*}_{h}-\lambda^{*}\mathbf{u}^{*}\|_{0,\tau},\\\label{s3.25}
&\eta^{*}_{F_{\tau}}
\lesssim \|\nabla (\mathbf{u}^{*}-\mathbf{u}^{*}_{h})\|_{0,\omega(\tau)}+\|p^{*}-p^{*}_{h}\|_{0,\omega(\tau)}
+( \sum\limits_{\tau\in\omega(\tau)}h_{\tau}^{2} \|\lambda^{*}\mathbf{u}^{*}-\lambda^{*}_{h}\mathbf{u}^{*}_{h}\|_{0,\tau}^{2}  )^{\frac{1}{2}},\\\label{s3.26}
&\eta^{*2}_{J_{\tau}}=\sum\limits_{F\subset\partial\tau, F\in\mathcal{E}^{i}_{h}}\gamma h_{F}^{-1}|[\![{\mathbf{u}^{*}_{h}-\mathbf{u}^{*}}]\!]|_{0,F}^{2}+ \sum\limits_{F\subset\partial\tau, F\in\mathcal{E}^{b}_{h}}\gamma h_{F}^{-1}||(\mathbf{u}^{*}_{h}-\mathbf{u}^{*})\otimes\mathbf{n}||_{0,F}^{2}.
\end{align}
\end{lemma}
\noindent {\bf Proof.}
Let $b_{\tau}\in H_{0}^{1}(\tau)$ and $b_{F}\in H_{0}^{1}(\omega(F))$  be the standard bubble function
on element $\tau$ and face $F$ ($d=3$) or edge $F$ ($d=2$) of $\tau$, respectively.
For any $\tau\in \mathcal{T}_{h}$, define the function $R$ and $K$ locally by
\begin{align*}
R|_{\tau}=\lambda_{h}\mathbf{u}_{h}+\mu\triangle\mathbf{u}_{h}-(\boldsymbol{\beta}\cdot\nabla)\mathbf{u}_{h}-\nabla p_{h},~~~K|_{\tau}=h_{\tau}^{2}Rb_{\tau}.
\end{align*}
From Proposition 3.37 in \cite{Verfurth2013}, using $\lambda\mathbf{u}+\mu\Delta \mathbf{u}-(\boldsymbol{\beta}\cdot\nabla)\mathbf{u}-\nabla p=0$ and integrating by parts, together with $\overline{K}|_{\partial\tau}=0$
we deduce that
\begin{align*}
&Ch_{\tau}^{2}\|R\|_{0,\tau}^{2}\leq \int_{\tau}R\cdot \overline{(h_{\tau}^{2}Rb_{\tau})}dx
=\int_{\tau}(\lambda_{h}\mathbf{u}_{h}+\mu\Delta \mathbf{u}_{h}-(\boldsymbol{\beta}\cdot\nabla)\mathbf{u}_{h}-\nabla p_{h})\cdot \overline{K}dx\\
&~~~= -\int_{\tau}\nabla (\mathbf{u}-\mathbf{u}_{h})\cdot \overline{\nabla K}dx+\int_{\tau}(p_{h}-p)
\overline{div K}dx\\
&~~~+\int_{\tau}(\lambda_{h}\mathbf{u}_{h}-\lambda\mathbf{u})\cdot \overline{K}dx
-\int_{\tau}(\boldsymbol{\beta}\cdot\nabla)(\mathbf{u}_{h}-\mathbf{u})\cdot \overline{K}dx.
\end{align*}
Applying the Cauchy-Schwarz inequality, Proposition 3.37 in \cite{Verfurth2013},
and noting $\|\nabla\cdot \mathbf{u}_{h}\|_{0}=\|\nabla\cdot (\mathbf{u}_{h}-\mathbf{u})\|_{0}$, we obtain (\ref{s3.21}).\\
\indent  For any interior edge/face $F\in\mathcal{E}_{h}^{i}$, let the functions $R$ and $\Theta$ be such that
\begin{align*}
R|_{F}=[\![p_{h}\mathbf{I}-\mu\nabla\mathbf{u}_{h} +\mathbf{u}_{h}\otimes\boldsymbol{\beta}]\!]|_{F}~~\rm{and}~~ \Theta=\newblock  {\it h_{F}Rb_{F}}.%
\end{align*}
Using Proposition 3.37 in \cite{Verfurth2013} and $\int_{F}[\![p\mathbf{I}- \mu\nabla\mathbf{u}+\mathbf{u}\otimes\boldsymbol{\beta}]\!]\cdot\Theta ds=0$, we get
\begin{align*}
Ch_{F}\|R\|_{0,F}^{2}\leq\int_{F}R\cdot \overline{(h_{F}Rb_{F})}ds=\int_{F}[\![(p_{h}-p)\mathbf{I}-\mu\nabla(\mathbf{u}_{h}-\mathbf{u})+(\mathbf{u}_{h}-\mathbf{u})\otimes\boldsymbol{\beta}]\!]\cdot\overline{\Theta} ds.
\end{align*}
Applying Green's formula over each element of $\omega(F)$, we derive
\begin{align*}
&Ch_{E}\|R\|_{0,F}^{2}\leq\int_{F}[\![((p_{h}-p)\mathbf{I}-\mu\nabla(\mathbf{u}_{h}-\mathbf{u}))]\!]\cdot\overline{\Theta} ds
+\int_{F}[\![(\mathbf{u}_{h}-\mathbf{u})\otimes\boldsymbol{\beta}]\!]\cdot\overline{\Theta} ds\\
&~~~=\sum\limits_{\tau\in\omega(F)}\int_{\tau}(-\mu\Delta (\mathbf{u}_{h}-\mathbf{u})+\nabla (p_{h}-p))\cdot\overline{\Theta} dx\\
&~~~-\sum\limits_{\tau\in\omega(F)}\int_{\tau}(\mu\nabla (\mathbf{u}-\mathbf{u}_{h})-(p-p_{h})\mathbf{I}):\overline{\nabla\Theta} dx+\int_{F}[\![(\mathbf{u}_{h}-\mathbf{u})\otimes\boldsymbol{\beta}]\!]\cdot\overline{\Theta} ds.
\end{align*}
Using $-\mu\Delta\mathbf{u}+\nabla p=\lambda\mathbf{u}-(\boldsymbol{\beta}\cdot\nabla)\mathbf{u}$, we deduce
\begin{align}\label{s3.27}
&Ch_{F}\|R\|_{0,F}^{2}\leq\sum\limits_{\tau\in\omega(F)}\int_{\tau}(\lambda_{h}\mathbf{u}_{h}+\mu\Delta\mathbf{u}_{h}
-(\boldsymbol{\beta}\cdot\nabla)\mathbf{u}_{h}-\nabla p_{h})\cdot\overline{\Theta} dx\nonumber\\
&~~~+\sum\limits_{\tau\in\omega(F)}\int_{\tau}(\lambda\mathbf{u}-\lambda_{h}\mathbf{u}_{h})\cdot\overline{\Theta} dx
+\sum\limits_{\tau\in\omega(F)}\int_{\tau}(-\mu\nabla (\mathbf{u}-\mathbf{u}_{h})+(p-p_{h})\mathbf{I}):\overline{\nabla\Theta} dx\nonumber\\
&~~~-\sum\limits_{\tau\in\omega (F)}\int_{\tau}((\boldsymbol{\beta}\cdot\nabla)(\mathbf{u}-\mathbf{u}_{h}))\cdot\overline{\Theta }dx
+\int_{F}[\![(\mathbf{u}_{h}-\mathbf{u})\otimes\boldsymbol{\beta}]\!]\cdot\overline{\Theta} ds.
\end{align}
By using the Cauchy-Schwarz inequality for each term on the right-hand side of (\ref{s3.27}),
and by (\ref{s3.21}) and Proposition 3.37 in \cite{Verfurth2013}, we get (\ref{s3.22}).
By using $[\![\mathbf{u}]\!]=0$ for any $F\in \mathcal{E}_{h}^{i}$, and $\mathbf{u}=0 $ for any $F\subset\partial\Omega$, we get (\ref{s3.23}).\\
\indent Similarly, we can prove (\ref{s3.24}), (\ref{s3.25}) and (\ref{s3.26}).
~~~~~~$\Box$\\
\indent From Lemma \ref{lem3-5}, we have the following theorem.
\begin{theorem}\label{th3-6}
Under the conditions of Theorem \ref{th3-1}, there hold
\begin{align}
&\eta_{h}(\mathbf{u}_{h},p_{h})^{2}\lesssim \|\mathbf{u}-\mathbf{u}_{h}\|_{h}^{2}+\|p-p_{h}\|_{0}^{2}+\sum\limits_{\tau\in\mathcal{T}_{h}}h_{\tau}^{2}\|\lambda \mathbf{u}-\lambda_{h}\mathbf{u}_{h}\|_{0,\tau}^{2},\label{s3.28}\\
&\eta^{*}_{h}(\mathbf{u}_{h}^{*},p_{h}^{*})^{2}\lesssim \|\mathbf{u}^{*}-\mathbf{u}^{*}_{h}\|_{h}^{2}+\|p^{*}-p^{*}_{h}\|_{0}^{2}+\sum\limits_{\tau\in\mathcal{T}_{h}}h_{\tau}^{2}
\|\lambda^{*} \mathbf{u}^{*}-\lambda^{*}_{h}\mathbf{u}^{*}_{h}\|_{0,\tau}^{2}.\label{s3.29}
\end{align}
\end{theorem}
\begin{theorem}\label{th3-7}
Assume that $\mathbf{R(\Omega)}$ holds, $\mathbb{M}(\lambda),  \mathbb{M}^{*}(\lambda)\subset (H^{1+s}(\Omega))^{d}\times H^{s}(\Omega)$, $\frac{1}{2}< s\leq k$, then
it is holds that
\begin{align}\label{s3.30}
|\lambda-\lambda_{h}|&\lesssim \eta_{h}^{2}(\mathbf{u}_{h},p_{h})+\sum\limits_{\tau\in\mathcal{T}_{h}}h_{\tau}^{2s}(||\mathbf{u}-\Pi_{h}\mathbf{u}||_{1+s,\tau}^{2}
+\|p-\vartheta_{h}p\|_{s}^{2})\nonumber\\
&+\eta_{h}^{*}(\mathbf{u}_{h}^{*},p_{h}^{*})^{2}+\sum\limits_{\tau\in\mathcal{T}_{h}}h_{\tau}^{2s}(||\mathbf{u}^{*}-\Pi_{h}\mathbf{u}^{*}||_{1+s,\tau}^{2}
+\|p^{*}-\vartheta_{h}p^{*}\|_{s}^{2})\nonumber\\
&+|(1+\theta)\sum\limits_{F\in\mathcal{E}_{h}}\int_{F}\mu\overline{\{{\nabla\mathbf{u}}^{*}\}}:[\![\underline{\mathbf{u}_{h}}]\!]ds|.
\end{align}
\end{theorem}
\noindent{\bf Proof.}
 Choosing  $ [\mathbf{v}_{h},q_{h}]=
[\mathbf{u}_{h}^{*},p_{h}^{*}]$ in (\ref{s2.65}), from (\ref{s2.15}),  (\ref{s2.16}), (\ref{s2.58}),  (\ref{s2.59}), (\ref{s3.19}) and (\ref{s3.20}), we obtain (\ref{s3.30}).~~~$\Box$\\
\indent
When eigenvalues are multiple, the discrete eigenfunctions obtained on different meshes may not approximate the same eigenfunction (see Theorem 7.4 in \cite{Babuska1991book}). Hence,
the a posteriori error estimate and adaptive computation of multiple and cluster eigenvalues have become an attractive issue (see, e.g. \cite{Dai2014,Gallistl2015,Boffi2017,Cances2020,yang2024}).\\
\indent  Denote the index set $\mathbb{J}=\{k,k+1,\cdots,k+m-1\}$. Let $\{\lambda_{j}\}_{j\in \mathbb{J}}$ be an eigenvalue cluster composed of eigenvalues of (\ref{s2.2})-(\ref{s2.3}), $\lambda_{k}\neq \lambda_{k-1}$ and $\lambda_{k+m-1}\neq \lambda_{k+m}$.
\begin{corollary}\label{coro3-8}
When the eigenvalue cluster $\{\lambda_{j}\}_{j\in\mathbb{J}}$ is composed of non-defective eigenvalues, there hold
\begin{align*}
&\sum\limits_{j\in \mathbb{J}}\left(\|{\mathbf{u}_{j}}-{\mathbf{u}}_{j,h}\|^{2}_{h}+\|{p_{j}}-{p_{j,h}}\|^{2}_{0}    \right)
\lesssim \sum\limits_{j\in \mathbb{J}}(\eta_{h}(\mathbf{u}_{j,h},p_{j,h})^{2}+\|\lambda \mathbf{u}_{j}-\lambda_{j,h}\mathbf{u}_{h}\|_{0}^{2}),\\
&\sum\limits_{j\in \mathbb{J}}\left(\|{\mathbf{u}_{j}}-{\mathbf{u}}_{j,h}\|^{2}_{h}+\|{p_{j}}-{p_{j,h}}\|^{2}_{0}    \right)
\gtrsim \sum\limits_{j\in \mathbb{J}}\eta_{h}(\mathbf{u}_{j,h},p_{j,h})^{2},\\
&\sum\limits_{j\in \mathbb{J}}\left(\|{\mathbf{u}_{j}^{*}}-{\mathbf{u}}_{j,h}^{*}\|^{2}_{h}+\|p_{j}^{*}-p_{j,h}^{*}\|^{2}_{0}    \right)
\lesssim\sum\limits_{j\in \mathbb{J}}(\eta_{h}^{*}(\mathbf{u}_{j,h}^{*},p_{j,h}^{*})^{2}+\|\lambda^* \mathbf{u}_{j}^*-\lambda_{j,h}^*\mathbf{u}_{h}^*\|_{0}^{2}),\\
&\sum\limits_{j\in \mathbb{J}}\left(\|{\mathbf{u}_{j}^{*}}-{\mathbf{u}}_{j,h}^{*}\|^{2}_{h}+\|p_{j}^{*}-p_{j,h}^{*}\|^{2}_{0}    \right)
\gtrsim\sum\limits_{j\in \mathbb{J}}\eta_{h}^{*}(\mathbf{u}_{j,h}^{*},p_{j,h}^{*})^{2},\\
&\left|{\lambda_{i}}-\lambda_{i,h}\right|\lesssim \sum\limits_{j\in \mathbb{J}}\left(\eta_{h}(\mathbf{u}_{j,h},p_{j,h})^{2}
+\eta^{*}_{h}(\mathbf{u}^{*}_{j,h},p^{*}_{j,h})^{2}+(1+\theta)\eta_{h}(\mathbf{u}_{j,h},p_{j,h})\right)\nonumber\\
&+\sum\limits_{\tau\in\mathcal{T}_{h}}h_{\tau}^{2r}(|\mathbf{u}-\Pi_{h}\mathbf{u}|_{1+r,\tau}^{2}+\|p-\vartheta_{h}p\|_{r}^{2}
+|\mathbf{u}^{*}-\Pi_{h}\mathbf{u}^{*}|_{1+r,\tau}^{2}
+\|p^{*}-\vartheta_{h}p^{*}\|_{r}^{2})
~~i\in \mathbb{J}.
\end{align*}
\end{corollary}
\noindent{\bf Proof.}
From Theorem \ref{th3-4} and Theorem \ref{th3-6}, we derive the first four estimates in the corollary.
From (\ref{s2.14}) and definition of $\eta_{h}$, we obtain
\begin{align*}
|\sum\limits_{F\in\mathcal{E}_{h}}\int_{F}\mu\overline{\{{\nabla\mathbf{u}_{j}}^{*}\}}:[\![\underline{\mathbf{u}_{j,h}}]\!]ds|
&\lesssim \sum\limits_{F\in\mathcal{E}_{h}}h_{F}^{\frac{1}{2}}\|\nabla\mathbf{u}_{j}^{*}\|_{0,F}h_{F}^{-\frac{1}{2}}\|[\![\underline{\mathbf{u}}_{j,h}]\!]\|_{0,F}\nonumber\\
&\lesssim \|\mathbf{u}_{j}^{*}\|_{1+r}\eta_{h}(\mathbf{u}_{j,h},p_{j,h}).
\end{align*}
Then, from Theorem \ref{th3-7}, we obtain the fifth estimate in the corollary.
~~~$\Box$
\indent For the {\bf SIP}  method ($\theta=-1$), Theorem \ref{th2-9} and Corollary \ref{coro3-8}
show that the reliability and efficiency of the estimator for approximate eigenfunctions, as well as the reliability of the estimator for approximate eigenvalues.
For the {\bf NIP} ($\theta=1$) and {\bf IIP} ($\theta=0$) methods, Theorem \ref{th2-10} and Corollary \ref{coro3-8} imply the a posteriori error estimates are not sharp.
\section{Numerical Experiments}
\indent Based on the standard adaptive algorithm, we obtain the following algorithm for multiple and cluster eigenvalues .\\
\indent{\bf Algorithm 1.}  The adaptive algorithm.\\
\indent Choose the parameter $\vartheta\in(0,1)$.\\
\indent {\bf Step 1.} Set $l=0$ and pick an initial mesh $\mathcal{T}_{h_{l}}$ with the mesh size $h_{l}$.\\
\indent {\bf Step 2.} Solve Eqs.(\ref{s2.8})-(\ref{s2.9}) and (\ref{s2.49}) on $\mathcal{T}_{h_{l}}$ for the discrete solution $\{(\lambda_{j,h_{l}}, \mathbf{u}_{j,h_{l}}, p_{j,h_{l}} )\}_{j\in\mathbb{J}}$ with $\|\mathbf{u}_{j,h_{l}}\|_{0}=1$ and $\{(\lambda^{*}_{j,h_{l}}, \mathbf{u}^{*}_{j,h_{l}}, p^{*}_{j,h_{l}} )\}_{j\in\mathbb{J}}$ with $\|\mathbf{u}^{*}_{j,h_{l}}\|_{0}=1$.\\
\indent {\bf Step 3.} Compute the local indicators $\{\eta_{\tau}(\mathbf{u}_{j,h_{l}}, p_{j,h_{l}})^{2}+ \eta^{*}_{\tau}(\mathbf{u}^{*}_{j,h_{l}}, p^{*}_{j,h_{l}})^{2}\}_{j\in\mathbb{J}}$ for all $\tau\in \mathcal{T}_{h_{l}}$.\\
\indent {\bf Step 4.} Construct $\widehat{\mathcal{T}}_{h_{l}}\subset \mathcal{T}_{h_{l}}$
by D\"{o}rfler \cite{Dorfler1996} {\bf Marking strategy }
\begin{align*}
&~~~\sum\limits_{\tau\in \widehat{\mathcal{T}}_{h_{l}}}(\sum\limits_{j\in\mathbb{J}}\eta_{\tau}(\mathbf{u}_{j,h_{l}}, p_{j,h_{l}})^{2}
+\eta^{*}_{\tau}(\mathbf{u}^{*}_{j,h_{l}}, p^{*}_{j,h_{l}})^{2})\nonumber\\
&\geq \vartheta \sum\limits_{\tau\in {\mathcal{T}}_{h_{l}}}(\sum\limits_{j\in\mathbb{J}}\eta_{\tau}(\mathbf{u}_{j,h_{l}}, p_{j,h_{l}})^{2}+ \eta^{*}_{\tau}(\mathbf{u}^{*}_{j,h_{l}}, p^{*}_{j,h_{l}})^{2} ).
\end{align*}
\indent {\bf Step 5.} Refine $\mathcal{T}_{h_{l}}$ to get a new mesh $\mathcal{T}_{h_{l+1}}$ by procedure {\bf Refine}.\\
\indent {\bf Step 6.} $l\Leftarrow l+1$ and go to Step 2.\\
 \indent  We demonstrate the efficiency of our approach by applying Algorithm 1 with the $P_{k}-P_{k-1}(k=2,3 )$ elements. The program is implemented using the $i$FEM package \cite{Chen2009}, and the matrix eigenvalue problem is solved using
the command $'eigs'$ in MATLAB 2019a on a DESKTOP-CP97IS3
with 3.00GHz CPU and 64GB RAM.
Numerical experiments are conducted in four domains: a 2D square domain, a 2D L-shaped domain, a 3D thick L-shaped domain, and the 3D unit cube domain. \\
\indent The choice of the penalty parameter $\gamma$ needs to satisfy conditions (2.22)-(\ref{s2.23}).
In the computations, the stabilization parameter $\gamma$ in the bilinear form $A_{h}(\cdot , \cdot )$ in problem (\ref{s2.8})-(\ref{s2.9})
 will be chosen proportionally to the square of the polynomial degree $k$ as
$\gamma=\widehat{\gamma}k^{2}$ with $\widehat{\gamma}> 0$.
We choose $\mu=1$ in all numerical experiments. The marking parameter $\vartheta$ is set to 0.5 for all domains.    \\
\indent The initial mesh size $h_{0}$ is set as follows:\\
\indent $h_{0}=\frac{\sqrt{2}}{16}$ for the L-shaped domain;\\
\indent $h_{0}=\frac{\sqrt{3}}{4}$ for the thick L-shaped domain and the unit cube domain.\\
 \indent We adapt the following symbols in our tables and figures:\\
\indent $\lambda_{j}$: the $j$th eigenvalue.\\
\indent $\lambda_{j,h_{l}}$: the $j$th approximate eigenvalue in the $l$th iteration.\\
\indent $dof$: the degrees of freedom in the $l$th iteration.\\
\indent $\eta_{j,h}$=$\eta(\mathbf{u}_{j,h_{l}},p_{j,h_{l}})$, $\eta_{j,h}^{*}$=$\eta^{*}(\mathbf{u}^{*}_{j,h_{l}},p^{*}_{j,h_{l}})$.\\
\indent  Let $\hat{\lambda}=\frac{1}{4}\sum\limits_{j=1}^{4}\lambda_{j}$ and $\hat{\lambda}_{h}=\frac{1}{4}\sum\limits_{j=1}^{4}\lambda_{j,h}$. Let
 $\text{Eff}=\frac{1}{4}(\sum\limits_{j=1}^{4}\eta_{j,h}^{2}+(\eta^{*}_{j,h})^{2})/|\hat{\lambda}-\hat{\lambda}_{h}|$ denote the effectivity index of $\hat{\lambda}_{h}$.\\
\indent Since the exact eigenvalues on the two-dimensional square domain are unknown, we adopt the formula
$order(\lambda_{j,h})\approx \frac{1}{(lg2)}lg\left|\frac{\lambda_{j,h}-\lambda_{j,\frac{h}{2}}}{\lambda_{j,\frac{h}{2}}-\lambda_{j,\frac{h}{4}}}\right|$
 to calculate the convergence order of the eigenvalues.
\subsection{The {\bf SIP} method on 2D domains}
\subsubsection{A 2D square domain}
\indent Let us begin with the two-dimensional square domain $\Omega=(-1,1)^{2}$. We take $\boldsymbol{\beta}(x,y)=(1,0)^{t}$ and $\widehat{\gamma}=10$.\\
\indent Table \ref{tab-1} depict the convergence history on uniform meshes for the first four computed eigenvalues. In all cases presented, the convergence order of the approximate eigenvalues tends to be $h^{2k}$ for $P_{k}-P_{k-1}(k=2,3)$ elements.
This indicates that the eigenfunctions are sufficiently smooth on the square domain.
 \begin{table}[htpb]
\caption{ \ The first four approximate eigenvalues and convergence orders obtained on uniform meshes for the square domain with $\boldsymbol{\beta}=(1,0)^{t}$ and $\widehat{\gamma}=10$.}\label{tab-1}
\scriptsize
\label{tab:1}
\begin{center}
\begin{tabular}{{c|c|cc|cc|cc|ccccccccccc}}
  \hline
       &  $h$ & $\lambda_{1,h}$ & $order$ & $\lambda_{2,h}$ & $order$ & $\lambda_{3,h}$ & $order$ & $\lambda_{4,h}$ & $order$  \\ \hline
       &  $\sqrt{2}/8$   & 13.659632  & ~      & 23.338104  & ~     & 23.652585  & ~     & 33.074986  &   \\
       &  $\sqrt{2}/16$  & 13.613052  & ~      & 23.145206  & ~     & 23.437585  & ~     & 32.354295  &   \\
     $k=2$  &  $\sqrt{2}/32$  & 13.609817  & 3.85   & 23.130753  & 3.74  & 23.423914  & 3.98  & 32.301805  & 3.78   \\
     {\bf SIP}  &  $\sqrt{2}/64$  & 13.609606  & 3.94   & 23.129813  & 3.94  & 23.423034  & 3.96  & 32.298369  & 3.93   \\
       &  $\sqrt{2}/128$ & 13.609593  & 3.98   & 23.129753  & 3.98  & 23.422979  & 3.98  & 32.298151  & 3.98   \\
  \hline
       &  $\sqrt{2}/8$   & 13.610408  & ~      & 23.134251  & ~     & 23.427119  & ~     & 32.320444  &   \\
     $k=3$  &  $\sqrt{2}/16$  & 13.609611  & ~      & 23.129842  & ~     & 23.423061  & ~     & 32.298602  &   \\
     {\bf SIP}  &  $\sqrt{2}/32$  & 13.609593  & 5.44   & 23.129751  & 5.59  & 23.422977  & 5.60  & 32.298145  & 5.58   \\
       &  $\sqrt{2}/64$  & 13.609592  & 5.52   & 23.129749  & 5.67  & 23.422975  & 5.62  & 32.298136  & 5.68   \\
  \hline
\end{tabular}\end{center}
\end{table}
\subsubsection{A 2D L-shaped domain}
\begin{figure}[htpb]
\centering
\begin{minipage}[t]{0.45\textwidth}
  \centering
 \includegraphics[width=6.5cm]{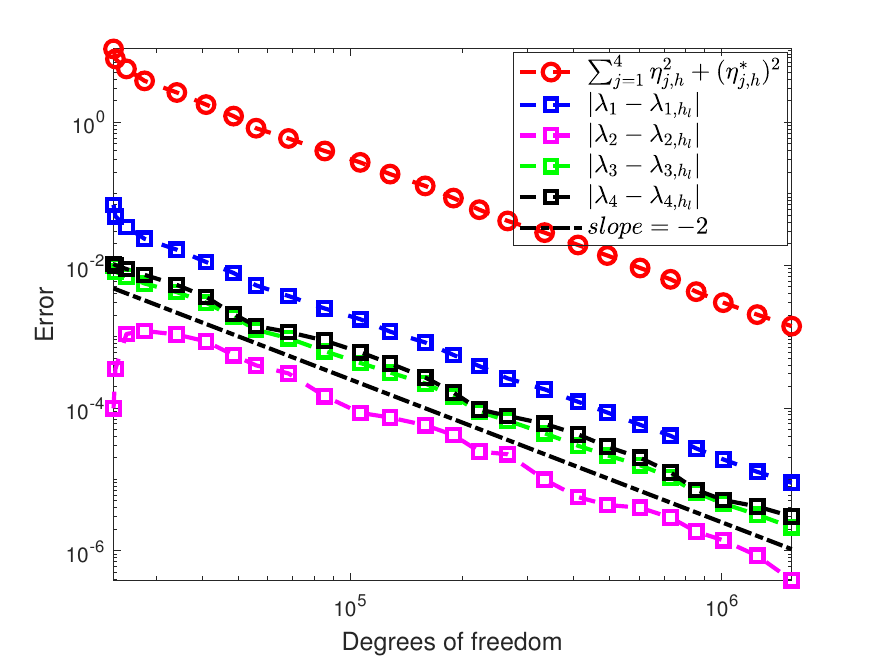}
\end{minipage}
\begin{minipage}[t]{0.45\textwidth}
  \centering
 \includegraphics[width=6.5cm]{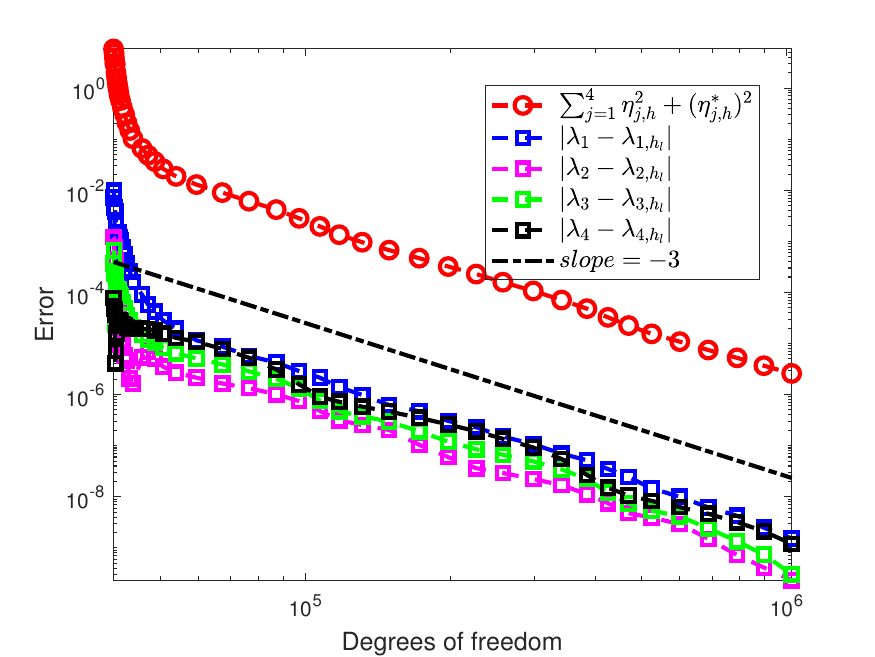}
\end{minipage}
\caption{The error curves on adaptive refinement meshes of the first four eigenvalues on L-shaped domain for $\boldsymbol{\beta}(x,y)=\boldsymbol{\beta}_{1}(x,y)$ (left: DG $P_{2}-P_{1}$; right: DG $P_{3}-P_{2}$).}\label{fig-1}
\end{figure}

\begin{figure}[htpb]
\centering
\begin{minipage}[t]{0.45\textwidth}
  \centering
 \includegraphics[width=6.5cm]{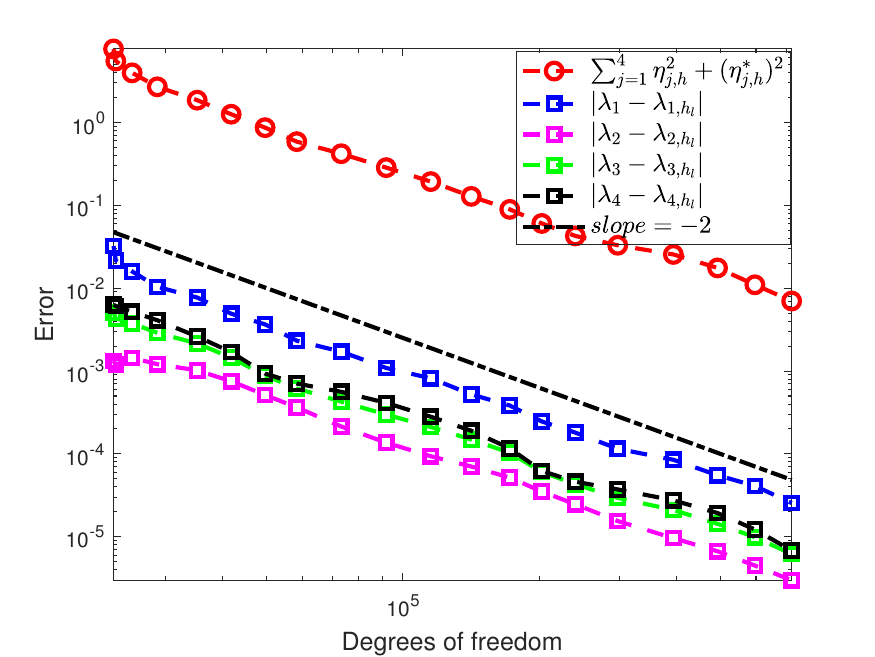}
\end{minipage}
\begin{minipage}[t]{0.45\textwidth}
  \centering
 \includegraphics[width=6.5cm]{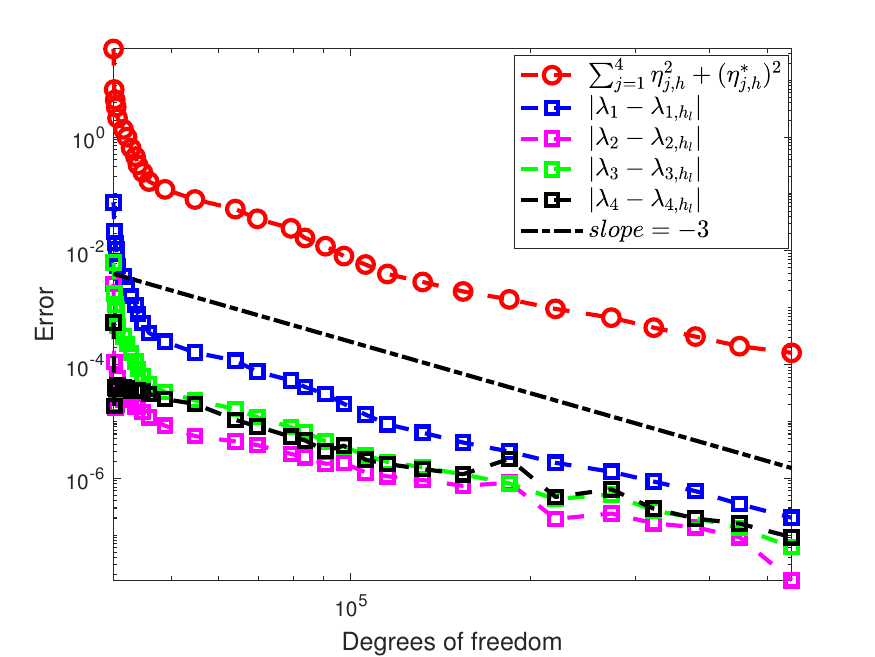}
\end{minipage}
\caption{The error curves on adaptive refinement meshes of the first four eigenvalues on L-shaped domain for $\boldsymbol{\beta}(x,y)=\boldsymbol{\beta}_{2}(x,y)$ (left: DG $P_{2}-P_{1}$; right: DG $P_{3}-P_{2}$).}\label{fig-2}
\end{figure}

\begin{figure}[htpb]
\centering
\begin{minipage}[t]{0.45\textwidth}
  \centering
 \includegraphics[width=6.5cm]{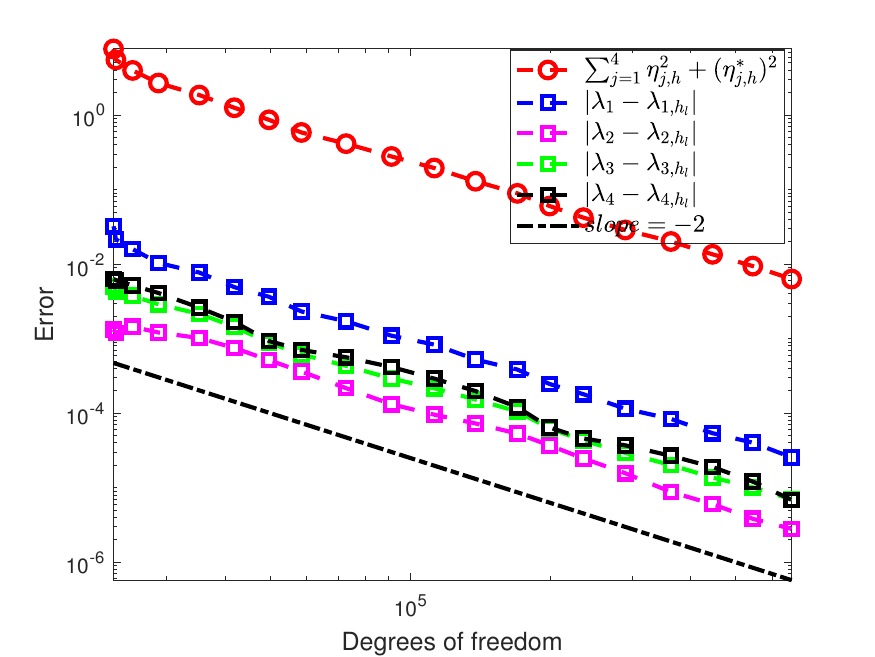}
\end{minipage}
\begin{minipage}[t]{0.45\textwidth}
  \centering
 \includegraphics[width=6.5cm]{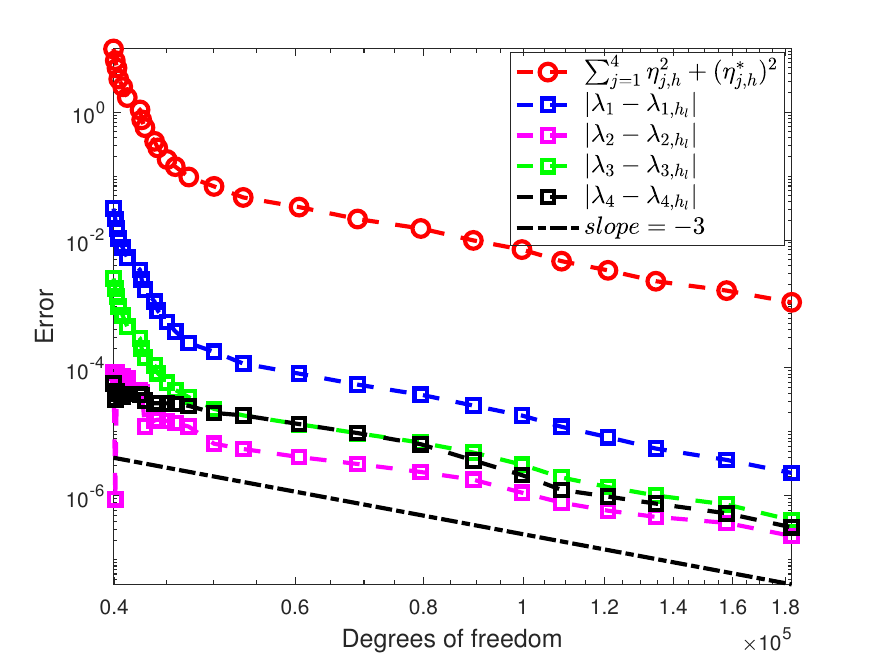}
\end{minipage}
\caption{The error curves on adaptive refinement meshes of the first four eigenvalues on L-shaped domain for $\boldsymbol{\beta}(x,y)=\boldsymbol{\beta}_{3}(x,y)$ (left: DG $P_{2}-P_{1}$; right: DG $P_{3}-P_{2}$).}\label{fig-3}
\end{figure}

\begin{figure}[htpb]
\centering
\begin{minipage}[t]{0.45\textwidth}
  \centering
 \includegraphics[width=6.5cm]{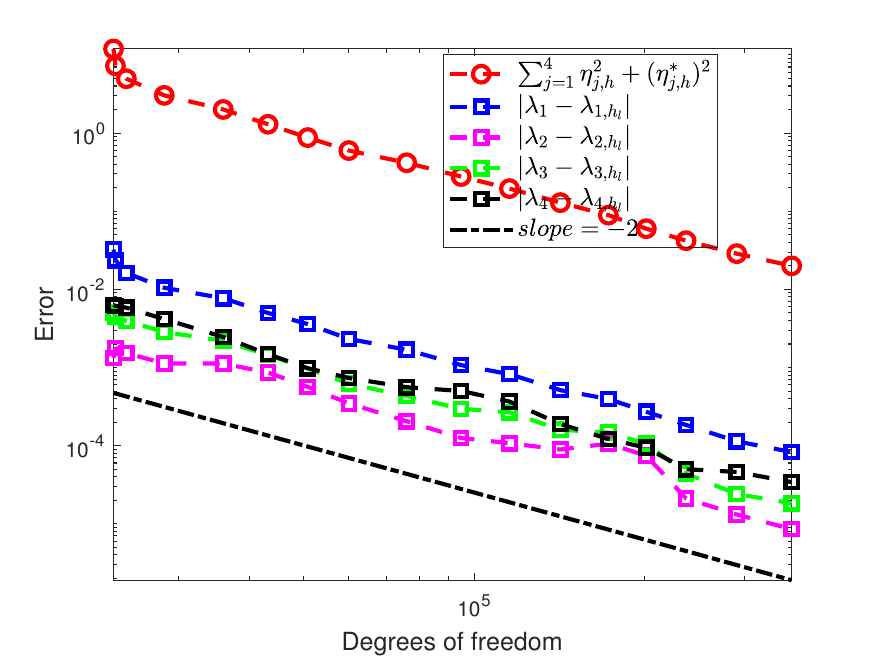}
\end{minipage}
\begin{minipage}[t]{0.45\textwidth}
  \centering
  \includegraphics[width=6.5cm]{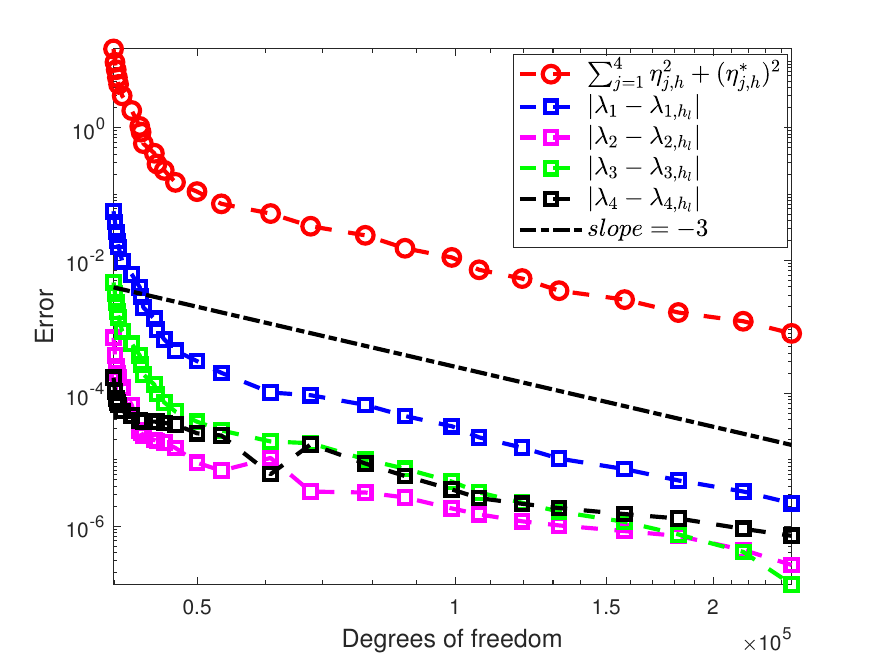}
\end{minipage}
\caption{The error curves on adaptive refinement meshes of the first four eigenvalues on L-shaped domain for $\boldsymbol{\beta}(x,y)=\boldsymbol{\beta}_{4}(x,y)$ (left: DG $P_{2}-P_{1}$; right: DG $P_{3}-P_{2}$).}\label{fig-4}
\end{figure}

\indent  We present the numerical experiments on the L-shaped domain $\Omega=(-1,1)^{2}\setminus (-1,0)^{2}$. We select $\widehat{\gamma}=10$.
Referring to \cite{Lepe2024},
 we consider four cases for the steady flow velocity $\boldsymbol{\beta}$:\\
\indent $\boldsymbol{\beta}_{1}(x,y)=(1,0)^{t}$, \\
\indent $\boldsymbol{\beta}_{2}(x,y)=(\cos(\pi x)\sin(\pi y),-\sin(\pi x)\cos(\pi y))^{t}$,\\
\indent $\boldsymbol{\beta}_{3}(x,y)=(y,-x)^{t}$,\\
\indent $\boldsymbol{\beta}_{4}(x,y)=(\frac{\partial\phi}{\partial y} , -\frac{\partial\phi}{\partial x}  )^{t}$
with $\phi(x,y)=1000(1-x^{2})^{2}(1-y^{2})^{2}$.\\
The functions are normalized in the $L^{\infty}$-norm as $\boldsymbol{\beta}=\boldsymbol{\beta}/\|\boldsymbol{\beta}\|_{\infty,\Omega}$.\\
 \indent The reference values for different $\boldsymbol{\beta}$ values are listed in Table \ref{tab-2}, which are obtained using the adaptive DG $P_{3}-P_{2}$ element with as many degrees of freedom as possible.
The effectivity indices are listed in Table \ref{tab-3}.  The error curves of the first four eigenvalues obtained using adaptive
  DG $P_{k}-P_{k-1}$ ($k=2,3$) elements are shown in Figures \ref{fig-1}-\ref{fig-4}.
From the error curves, both the eigenvalue error curves and the estimator error
 curves are approximately parallel to a straight line with slope $-k~(k=2,3)$, indicating that the error estimator is reliable and efficient,
  and that highly accurate approximate eigenvalues can be obtained.
 \begin{table}[htpb]
\caption{ \ The reference values for the first four eigenvalues with different $\boldsymbol{\beta}$ on the L-shaped domain.}\label{tab-2}
\scriptsize
\label{tab:1}
\begin{center}
\begin{tabular}{{c|ccccccccccccccccc}}
  \hline
              & $\lambda_{1}$  &$\lambda_{2}$  &$\lambda_{3}$  &$\lambda_{4}$ & $dof$\\ \hline
        $\boldsymbol{\beta}_{1}(x,y)$ &   32.9600408 & 37.1171925 & 42.3976456  & 49.2536801 &1348646 \\
        $\boldsymbol{\beta}_{2}(x,y)$ &   32.1326948 & 37.0183349 & 41.9398325  & 48.9835841 &936052 \\
        $\boldsymbol{\beta}_{3}(x,y)$ &   32.1326936 & 37.0183325 & 41.9398199  & 48.9835513 &1037972 \\
        $\boldsymbol{\beta}_{4}(x,y)$ &   32.1326948 & 37.0183349 & 41.9398325  & 48.9835841 &1185808 \\
\hline
\end{tabular}\end{center}
\end{table}
 \begin{table}[htpb]
\caption{ \ The effectivity indices Eff for $\hat{\lambda}_{h}$ obtained by the adaptive DG $P_{k}-P_{k-1}(k=2,3)$  element on L-shaped domain. }\label{tab-3}
\scriptsize
\label{tab:1}
\begin{center}
\begin{tabular}{{c|ccccccccccccccccc}}
  \hline
          $dof$ & 42765 & 59250 & 74700 & 113640 & 201555 & 367680 & 445950 & 540225  \\
        $P_{2}-P_{1}$ & 437  & 454  & 435  & 438  & 456  & 437  & 445  & 423   \\ \hline
        $dof$ & 46904 & 53274 & 69316 & 90194 & 99528 & 110162 & 136682 & 184418  \\
        $P_{3}-P_{2}$ & 1383  & 1362  & 1317  & 1323  & 1317  & 1353  & 1332  & 1315   \\
\hline
\end{tabular}\end{center}
\end{table}
\subsection{The {\bf SIP} method on 3D domains}
\subsubsection{A 3D unit cube domain}
\begin{figure}[htpb]
\centering
\begin{minipage}[t]{0.45\textwidth}
  \centering
 \includegraphics[width=6.5cm]{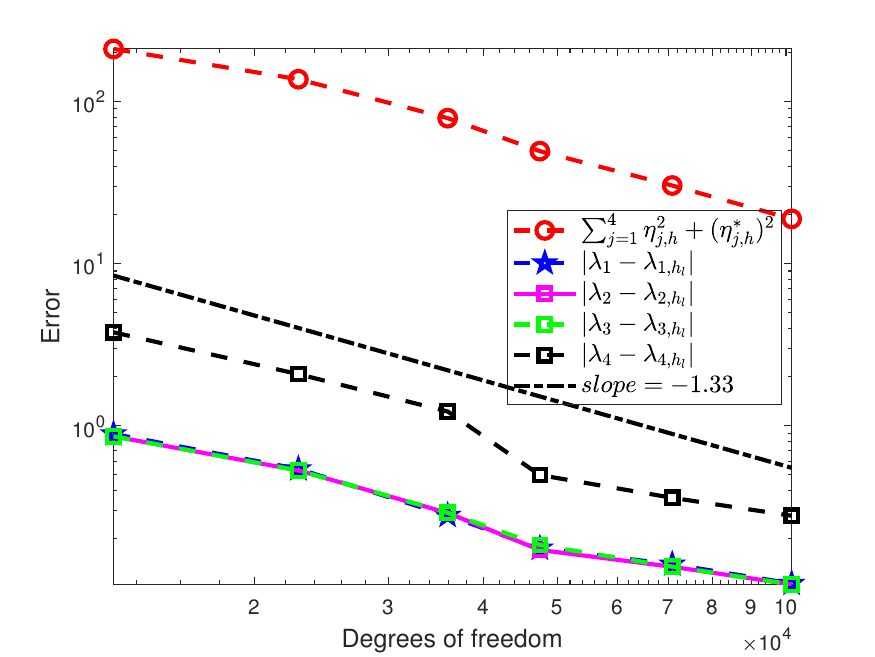}
\end{minipage}
\begin{minipage}[t]{0.45\textwidth}
  \centering
 \includegraphics[width=6.5cm]{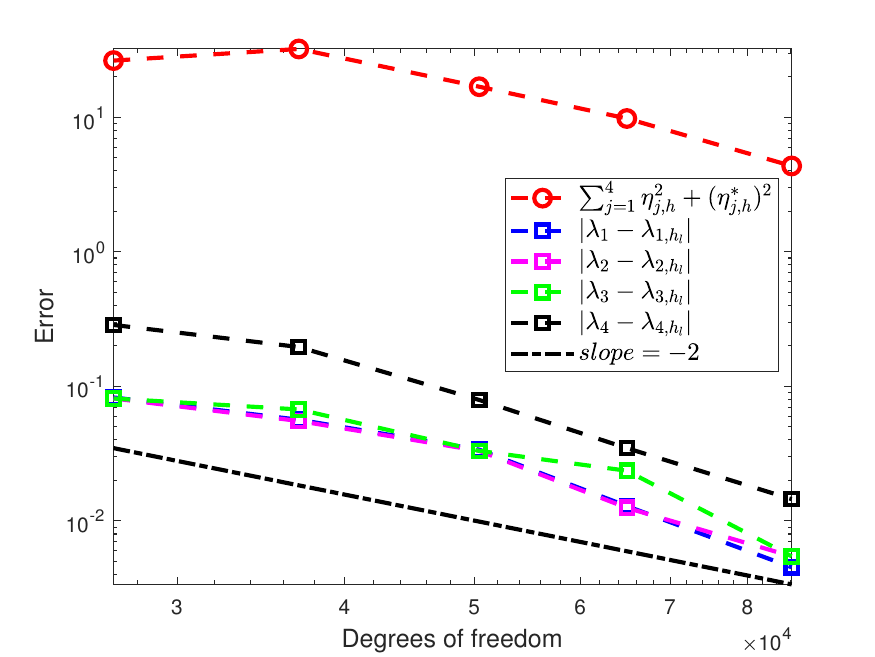}
\end{minipage}
\caption{The error curves on adaptive refinement meshes for the first four eigenvalues on the unit cube domain for $\boldsymbol{\beta}=(0,0,1)^{t}$ (left: DG $P_{2}-P_{1}$; right: DG $P_{3}-P_{2}$).}\label{fig-5}
\end{figure}

\indent In this experiment, we consider the unit cube domain $\Omega=(0,1)^{3}$. We take $\boldsymbol{\beta}=(0,0,1)^{t}$ and $\widehat{\gamma}=10$ for the Oseen eigenvalue problem. The extrapolated values $\lambda_{1}=62.4253$, $\lambda_{2}=62.7107$, $\lambda_{3}=62.7107$ and $\lambda_{4}=91.8801$ (from Table 5 in \cite{Lepe2024new}) are adopted as reference values for the first four eigenvalues.
We use the formula:
$order(\lambda_{j,h_{i}})\approx \frac{1}{(lg\frac{h_{i}}{h_{i+1}})}lg\left|\frac{\lambda_{j,h_{i}}-\lambda_{j}}{\lambda_{j,h_{i+1}}-\lambda_{j}}\right|$
to compute the approximate convergence order and the numerical solutions of the first four eigenvalues on uniform meshes are listed in Table \ref{tab-4}.
From the table, we observe that the numerical solutions achieve the optimal convergence order.\\
\indent The error curves of the first four eigenvalues on adaptive refinement meshes, computed by Algorithm 1, are shown in Figure \ref{fig-5}. From the figure, it can be observed that the error curves are almost parallel to the straight line with a slope of $-\frac{2k}{3}(k=2,3)$. This indicates that the numerical results obtained by the adaptive method basically achieve the optimal convergence order.
   \begin{table}[htpb]
\caption{ \ The first four eigenvalues on uniform meshes for the cube domain with $\boldsymbol{\beta}=(0,0,1)^{t}$.}\label{tab-4}
\scriptsize
\label{tab:1}
\begin{center}
\begin{tabular}{{c|c|cc|cc|cc|ccccccccccc}}
  \hline
     &  $h$ & $\lambda_{1,h}$ & $order$ & $\lambda_{2,h}$ & $order$ & $\lambda_{3,h}$ & $order$ & $\lambda_{4,h}$ & $order$  \\ \hline
     &   $\sqrt{3}/4$ & 64.529422  & ~ & 65.271012  & ~ & 65.445292  & ~ & 97.321863  &   \\
   $k=2$  &   $\sqrt{3}/6$ & 62.998928  & 3.21 & 63.271020  & 3.75 & 63.372723  & 3.50 & 93.402728  & 3.14  \\
   {\bf SIP}  &   $\sqrt{3}/8$ & 62.630616  & 3.57 & 62.896457  & 3.84 & 62.938185  & 3.71 & 92.434393  & 3.51  \\
     &   $\sqrt{3}/10$ & 62.514097  & 3.76 & 62.789536  & 3.84 & 62.807989  & 3.81 & 92.122364  & 3.71  \\
     &   $\sqrt{3}/12$ & 62.469473  & 3.83 & 62.749684  & 3.86 & 62.758883  & 3.85 & 92.001338  & 3.80  \\
  \hline
     &  $\sqrt{3}/4$ & 62.504343  & ~ & 62.781408  & ~ & 62.795545  & ~ & 92.161174  &   \\
   $k=3$  &   $\sqrt{3}/6$ & 62.433502  & 5.59 & 62.717994  & 5.49 & 62.719422  & 5.61 & 91.910104  & 5.52  \\
    {\bf SIP} &   $\sqrt{3}/8$ & 62.426993  & 5.49 & 62.712205  & 5.60 & 62.712479  & 5.53 & 91.886366  & 5.44  \\
  \hline
\end{tabular}\end{center}
\end{table}
\subsubsection{ A 3D thick L-shaped domain}
\begin{figure}[htpb]
\centering
\begin{minipage}[t]{0.45\textwidth}
  \centering
 \includegraphics[width=6.5cm]{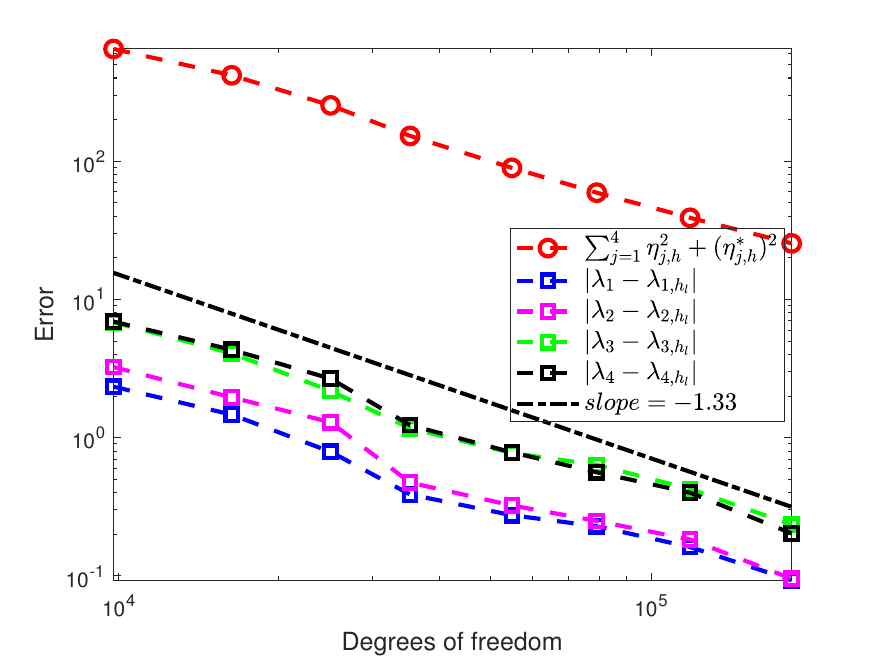}
\end{minipage}
\begin{minipage}[t]{0.45\textwidth}
  \centering
 \includegraphics[width=6.5cm]{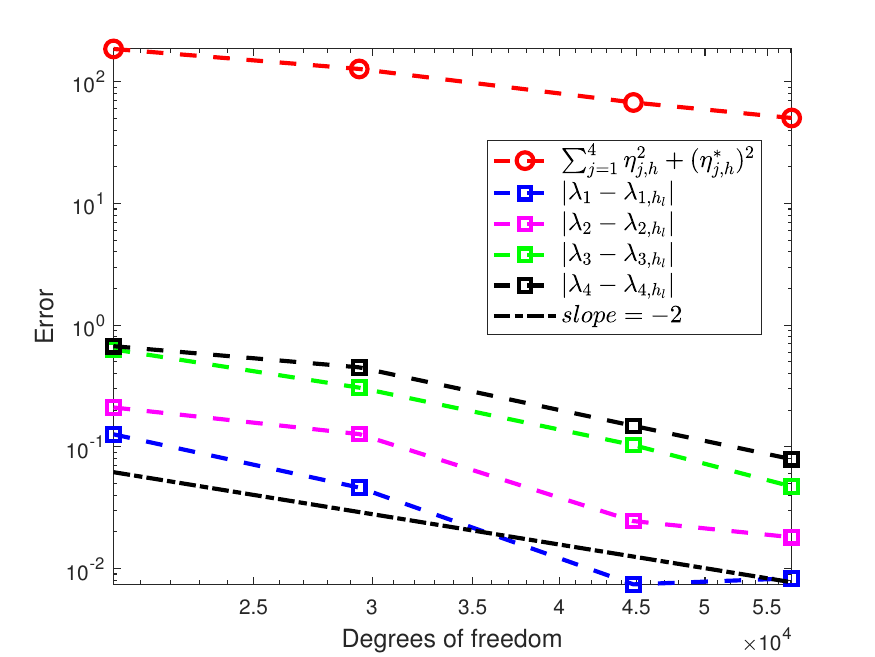}
\end{minipage}
\caption{The error curves on adaptive refinement meshes for the first four eigenvalues on the thick L-shaped domain for $\boldsymbol{\beta}=(0,0,1)^{t}$ (left: DG $P_{2}-P_{1}$; right: DG $P_{3}-P_{2}$).}\label{fig-6}
\end{figure}

\indent In this experiment, we consider the thick L-shaped domain $\Omega=(-0.5,0.5)\times (0,1)\times (-0.5,0.5)\setminus ((0,0.5)\times (0,1)\times (0,0.5))$. We take $\boldsymbol{\beta}=(0,0,1)^{t}$ and $\widehat{\gamma}=10$ for
the Oseen eigenvalue problem.
 The reference values for the first four eigenvalues, obtained by an adaptive DG $P_{3}-P_{2}$ element with 120121 degrees of freedom, are $\lambda_{1}=82.7711$, $\lambda_{2}=88.9986$, $\lambda_{3}=124.8892$ and $\lambda_{4}=125.7274$.\\
\indent The effectivity indices are listed in Table \ref{tab-5}. The error curves on adaptive refinement meshes of the first four eigenvalues are shown in Figure \ref{fig-6}.
From the figure, it can be observed that the eigenvalue error curves and the error estimator curves are approximately parallel to a straight line with a slope of $\frac{2k}{3}(k=2,3)$. This indicates that the error estimator is reliable and effective. \\
 \begin{table}[htpb]
\caption{ \  The effectivity indices $\text{Eff}$ for $\hat{\lambda}_{h}$ obtained by the adaptive DG $P_{2}-P_{1}$  element on the 3D thick L-shaped domain. }\label{tab-5}
\scriptsize
\label{tab:1}
\begin{center}
\begin{tabular}{{c|ccccccccccccccccc}}
  \hline
          $dof$                           & 9792 & 16320    & 25024 & 35292  & 54808 & 79084 & 118388 & 183668  \\
          \hline
        $\mathbb{P}_{2}-\mathbb{P}_{1}$   & 134  & 142      & 146   & 187    & 166   & 141   & 133    & 163   \\
\hline
\end{tabular}\end{center}
\end{table}
\subsection{The comparison among the {\bf SIP}, {\bf NIP} and {\bf IIP} methods}
\begin{figure}[htpb]
\centering
\begin{minipage}[t]{0.45\textwidth}
  \centering
 \includegraphics[width=6.5cm]{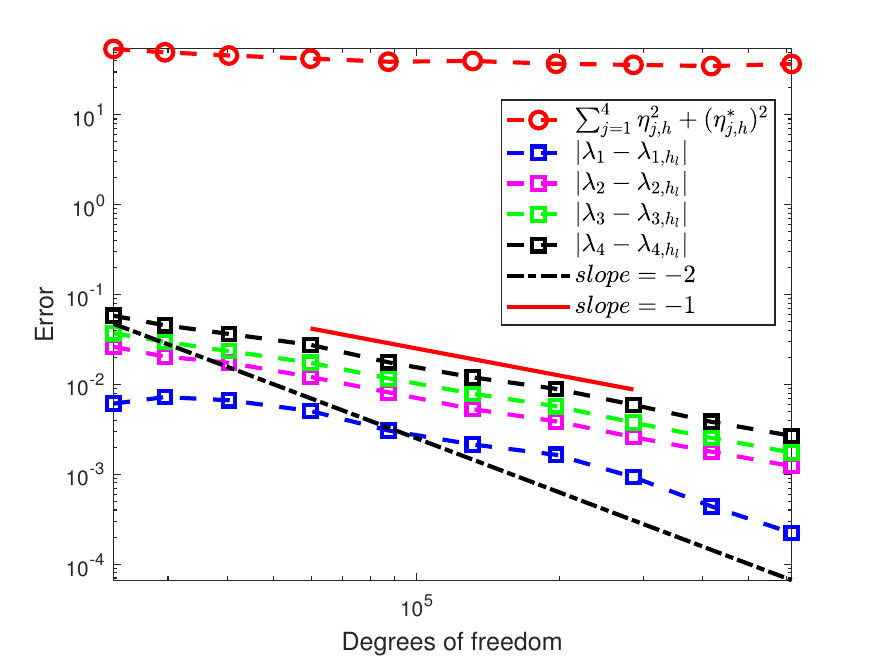}
\end{minipage}
\begin{minipage}[t]{0.45\textwidth}
  \centering
  \includegraphics[width=6.5cm]{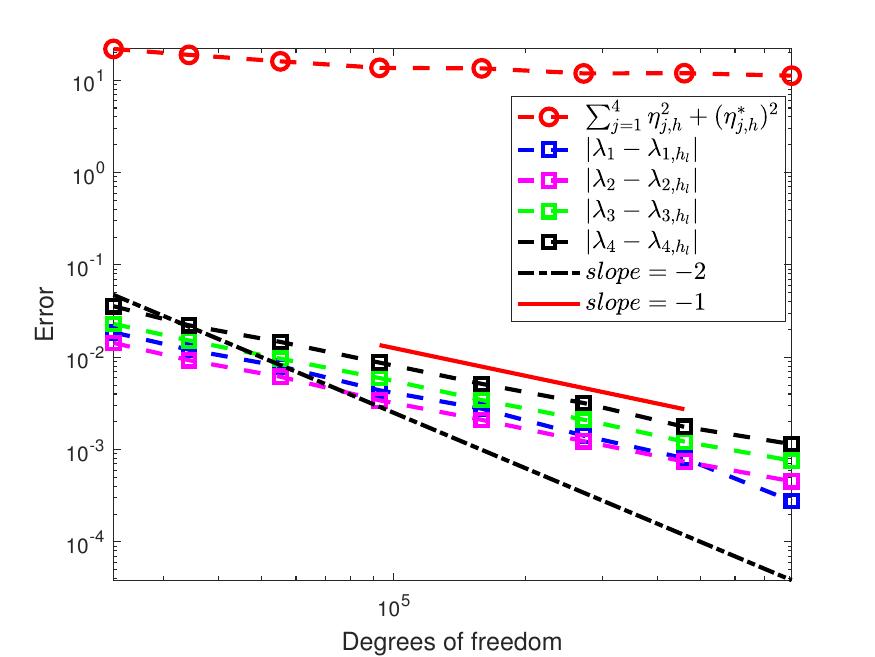}
\end{minipage}
\caption{  The error curves on adaptive refinement meshes of the first four eigenvalues using the $P_{2}-P_{1}$ elements on the 2D L-shaped domain for $\boldsymbol{\beta}(x,y)=\boldsymbol{\beta}_{1}(x,y)$ (left: {\bf NIP}; right: {\bf IIP}).}\label{fig-7}
\end{figure}

\begin{figure}[htpb]
\centering
\begin{minipage}[t]{0.45\textwidth}
  \centering
 \includegraphics[width=6.5cm]{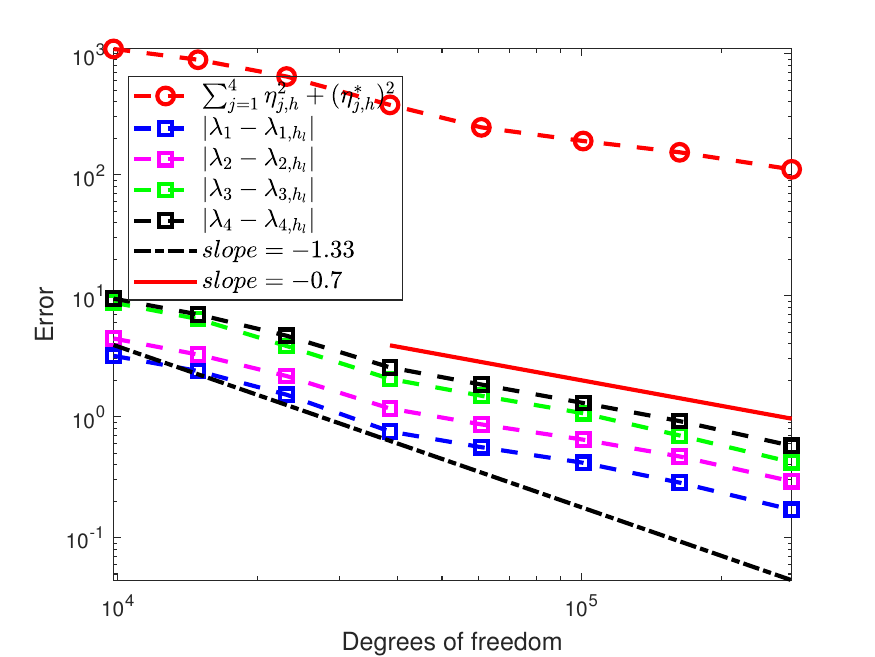}
\end{minipage}
\begin{minipage}[t]{0.45\textwidth}
  \centering
  \includegraphics[width=6.5cm]{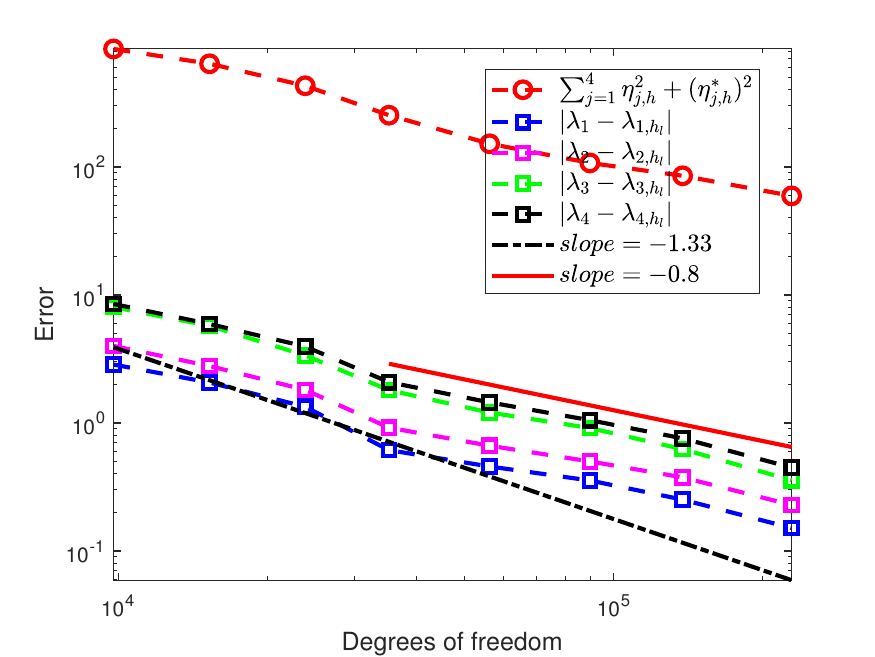}
\end{minipage}
\caption{  The error curves on adaptive refinement meshes of the first four eigenvalues using the $P_{2}-P_{1}$ elements on the 3D thick L-shaped domain for $\boldsymbol{\beta}=(0,0,1)^{t}$ (left: {\bf NIP}; right: {\bf IIP}).}\label{fig-8}
\end{figure}

 \indent  In this section, we compare the {\bf SIP}, {\bf NIP} and {\bf IIP} methods. \\
\indent For 2D square domain, we report the convergence orders in Table \ref{tab-6} for the first four eigenvalues computed by the NIP and IIP methods.
We take $\widehat{\gamma}=2$ and $\boldsymbol{\beta}=(1,0)^{t}$. We consider
uniform meshes with $h=\sqrt{2}/8$, $\sqrt{2}/16$, $\sqrt{2}/32$ and $\sqrt{2}/64$. As shown in Table \ref{tab-6}, both the NIP and IIP methods fail to reach $h^{2k}$, whereas the SIP method (Table \ref{tab-1}) achieves the optimal convergence order.\\
\indent For the 2D L-shaped domain, we list in Table \ref{tab-8} the first four approximate eigenvalues obtained via the SIP discrete method on a fixed uniform mesh with $h=\sqrt{2}/16$ and $\boldsymbol{\beta}=(1,0)^{t}$,
where $\widehat{\gamma}$ takes the values of $1/4, 1/2, 1, 2, 4, 8$. The eigenvalues inside boxes correspond to spurious eigenvalues.
As shown in Table \ref{tab-8}, the SIP method produces spurious eigenvalues when the stabilization parameter $\widehat{\gamma}$ is small enough and
 such spurious eigenvalues gradually disappear from the spectrum with the increase of $\widehat{\gamma}$.\\
\indent In Table \ref{tab-9}, we present the first four computed eigenvalues obtained by the NIP and IIP discrete methods. The computations are carried out on a fixed uniform mesh with
 $h=\sqrt{2}/16$ and $\boldsymbol{\beta}=(1,0)^{t}$ , where $\widehat{\gamma}$ is set to $1/32, 1/16, 1/8, 1/4, 1/2,1$. As shown in Table \ref{tab-9},
 no spurious eigenvalues are generated by the NIP and IIP methods with a small stabilization parameter $\widehat{\gamma}$.\\
\indent For the 3D thick L-shaped domain, we report the first four approximate eigenvalues in Table \ref{tab-10} on uniform mesh with $h=\sqrt{3}/8$, $\boldsymbol{\beta}=(0,0,1)^{t}$
and different values of $\widehat{\gamma}=1/4, 1/2, 1, 2, 4, 8$, which are obtained with the SIP discrete method.
In Table \ref{tab-11}, we present the first four approximate eigenvalues computed by the NIP and IIP discretization methods for different parameter values of $\widehat{\gamma}=1/32, 1/16, 1/8, 1/4, 1/2, 1$.
As shown in the Tables \ref{tab-10} and \ref{tab-11}, the SIP method produces spurious eigenvalues when the stabilization parameter $\widehat{\gamma}$ is small enough and
 such spurious eigenvalues gradually disappear from the spectrum with the increase of $\widehat{\gamma}$. In contract, the NIP and IIP methods need a smaller stabilization parameter $\widehat{\gamma}$ to avoid spurious eigenvalues.\\
\indent The error curves obtained by the nonsymmetric NIP and IIP methods on the two-dimensional L-shaped domain and three-dimensional thick L-shaped domain are presented in Figures \ref{fig-7} and \ref{fig-8}, respectively.
It can be observed that all error curves fail to parallel the straight line with a slope of $\frac{2k}{d}$. In contrast, the symmetric SIP method (see Figures \ref{fig-1} and \ref{fig-6}) achieves the optimal convergence order.\\
\indent Nevertheless, the numerical results demonstrate that although the {\bf NIP} and {\bf IIP} methods cannot achieve the optimal convergence order,
these two methods yield fewer spurious eigenvalues than the {\bf SIP} method for a fixed small penalty parameter $\gamma$.\\
 \begin{table}[htpb]
\caption{ \ The first four approximate eigenvalues and convergence orders obtained on uniform meshes for the {\bf NIP} method on the square domain with $\boldsymbol{\beta}=(1,0)^{t}$ and $\widehat{\gamma}=2$ ( top: $P_{2}-P_{1}$ element; bottom: $P_{3}-P_{2}$ element).}\label{tab-6}
\scriptsize
\label{tab:1}
\begin{center}
\begin{tabular}{{c|c|cc|cc|cc|ccccccccccc}}
  \hline
    &   $h$ & $\lambda_{1,h}$ & $order$ & $\lambda_{2,h}$ & $order$ & $\lambda_{3,h}$ & $order$ & $\lambda_{4,h}$ & $order$  \\ \hline
    &    $\sqrt{2}/8$ & 13.637790  & ~ & 23.229738  & ~ & 23.523373  & ~ & 32.484282  && ~  ~ \\
   $k=2$ &     $\sqrt{2}/16$ & 13.616553  & ~ & 23.154298  & ~ & 23.447612  & ~ & 32.344069  \\
   {\bf NIP} &     $\sqrt{2}/32$ & 13.611325  & 2.02  & 23.135843  & 2.03  & 23.429088  & 2.03  & 32.309605  & 2.02 \\
    &     $\sqrt{2}/64$ & 13.610025  & 2.01  & 23.131268  & 2.01  & 23.424498  & 2.01  & 32.301005  & 2.00 \\
  \hline
    &     $\sqrt{2}/8$ & 13.609748  & ~ & 23.130860  & ~ & 23.424083  & ~ & 32.301498 & ~  \\
   $k=3$ &  $\sqrt{2}/16$  & 13.609602  & ~ & 23.129819  & ~ & 23.423045  & ~ & 32.298350 & ~  \\
   {\bf NIP} &     $\sqrt{2}/32$ & 13.609593  & 3.95  & 23.129753  & 3.99  & 23.422979  & 3.99  & 32.298150 & 3.97  \\
    &     $\sqrt{2}/64$ & 13.609592  & 4.00  & 23.129749  & 4.00  & 23.422975  & 4.00  & 32.298137 & 3.99  \\
  \hline
     &   $\sqrt{2}/8$ & 13.646422  & ~ & 23.256843  & ~ & 23.550372  & ~ & 32.536387  &   \\
   $k=2$  &   $\sqrt{2}/16$& 13.618656  & ~ & 23.161009  & ~ & 23.454250  & ~ & 32.357831  &   \\
    {\bf IIP} &   $\sqrt{2}/32$ & 13.611838  & 2.03  & 23.137483  & 2.03  & 23.430706  & 2.03  & 32.313036  & 1.99   \\
     &   $\sqrt{2}/64$ & 13.610151  & 2.01  & 23.131671  & 2.02  & 23.424896  & 2.02  & 32.301855  & 2.00   \\
  \hline
     &   $\sqrt{2}/8$ & 13.609738  & ~ & 23.130802  & ~ & 23.424021  & ~ & 32.301277  &   \\
  $k=3$   &   $\sqrt{2}/16$& 13.609601  & ~ & 23.129815  & ~ & 23.423040  & ~ & 32.298336  &  \\
  {\bf IIP}   &   $\sqrt{2}/32$& 13.609593  & 3.99  & 23.129753  & 4.00  & 23.422979  & 4.01  & 32.298149  & 3.97   \\
     &   $\sqrt{2}/64$ & 13.609592  & 3.84  & 23.129749  & 4.02  & 23.422975  & 4.01  & 32.298137  & 4.00   \\
    \hline
\end{tabular}\end{center}
\end{table}
 \begin{table}[htpb]
\caption{ \  The first four eigenvalues using $P_{k}-P_{k-1} (k=3)$ element on uniform meshes for the {\bf SIP} method on the 2D L-shaped domain with $\boldsymbol{\beta}=(1,0)^{t}$, $h=\sqrt{2}/16$, and different stabilization parameters $\widehat{\gamma}$.  }\label{tab-8}
\scriptsize
\label{tab:1}
\begin{center}
\begin{tabular}{{c|cccccccccccccccccc}}
  \hline
      &  $\widehat{\gamma}=1/4$ & $\widehat{\gamma}=1/2$ & $\widehat{\gamma}=1$ & $\widehat{\gamma}=2$ & $\widehat{\gamma}=4$ & $\widehat{\gamma}=8$  \\ \hline
      &  \fbox{0.9299}                    & \fbox{$-$12.0892} & \fbox{$-$24.9244}  & 32.9136  & 32.9456  & 32.9746   \\
   $k=3$   &  \fbox{$-$1.8040 }                & \fbox{$-$9.2965}  & \fbox{$-$14.0705}  & 37.1176  & 37.1167  & 37.1157   \\
   {\bf SIP}   &  \fbox{$-$5.9239 $-$ 0.5651i}     & \fbox{$-$8.0163}  & \fbox{$-$2.5208}   & 42.3938  & 42.3963  & 42.3986   \\
      &  \fbox{$-$5.9239 $+$ 0.5651i }    & \fbox{$-$3.5788}  & \fbox{14.6101}     & 49.2533  & 49.2535  & 49.2536   \\
\hline
\end{tabular}\end{center}
\end{table}
 \begin{table}[htpb]
\caption{ \  The first four eigenvalues using $P_{k}-P_{k-1}(k=3)$ element on uniform meshes for the {\bf NIP} and {\bf IIP} methods on the 2D L-shaped domain with $\boldsymbol{\beta}=(1,0)^{t}$, $h=\sqrt{2}/16$, and different stabilization parameters $\widehat{\gamma}$. }\label{tab-9}
\scriptsize
\label{tab:1}
\begin{center}
\begin{tabular}{{c|cccccccccccccccccc}}
  \hline
      &  $\widehat{\gamma}=1/32$ & $\widehat{\gamma}=1/16$ & $\widehat{\gamma}=1/8$ & $\widehat{\gamma}=1/4$ & $\widehat{\gamma}=1/2$ & $\widehat{\gamma}=1$ & $\widehat{\gamma}=2$  \\ \hline
      &  32.9182  & 32.9147  & 32.9088  & 32.9002  & 32.8909  & 32.8867  & 32.8948   \\
    {\bf NIP}  &  37.1142  & 37.1145  & 37.1149  & 37.1157  & 37.1168  & 37.1178  & 37.1183   \\
    $k=3$  &  42.3965  & 42.3961  & 42.3954  & 42.3945  & 42.3934  & 42.3927  & 42.3930   \\
      &  49.2546  & 49.2545  & 49.2543  & 49.2541  & 49.2539  & 49.2538  & 49.2537   \\
    \hline
       & 32.8216  & 32.8198  & 32.8229  & 32.8309  & 32.8446  & 32.8652  & 32.8932   \\
   {\bf IIP}    & 37.1265  & 37.1228  & 37.1209  & 37.1199  & 37.1194  & 37.1191  & 37.1185   \\
     $k=3$  & 42.3953  & 42.3913  & 42.3896  & 42.3892  & 42.3895  & 42.3907  & 42.3925   \\
       & 49.2668  & 49.2598  & 49.2564  & 49.2547  & 49.2539  & 49.2536  & 49.2536   \\
\hline
\end{tabular}\end{center}
\end{table}
 \begin{table}[htpb]
\caption{ \ The first four eigenvalues using $P_{k}-P_{k-1}(k=3)$ element for the {\bf SIP} method on uniform meshes for the 3D thick L-shaped domain with $h=\sqrt{3}/8$, $\boldsymbol{\beta}=(0,0,1)^{t}$, and different stabilization parameters $\widehat{\gamma}$. }\label{tab-10}
\scriptsize
\label{tab:1}
\begin{center}
\begin{tabular}{{c|cccccccccccccccccc}}
  \hline
       & $\widehat{\gamma}=1/4$ & $\widehat{\gamma}=1/2$ & $\widehat{\gamma}=1$ & $\widehat{\gamma}=2$ & $\widehat{\gamma}=4$  & $\widehat{\gamma}=8$ ~ \\ \hline
       & \fbox{$-$0.020}             &\fbox{  $-$0.444}               &\fbox{ 0.066}               &\fbox{ $-$7.035}  &\fbox{ 47.941}  & 82.721   \\
     {\bf SIP}  & \fbox{0.265 }             &\fbox{  $-$1.038}              & \fbox{ $-$0.219}             & \fbox{ $-$5.536}  & 82.657         & 88.988   \\
     $k=3$  & \fbox{$-$0.926}             &\fbox{ 1.626}                &\fbox{  $-$0.468 $-$ 0.684i}  & \fbox{ $-$3.293}  & 88.984         & 124.830   \\
       & \fbox{$-$0.983}             &\fbox{  $-$1.689}               &\fbox{  $-$0.468 + 0.684i}  & \fbox{ $-$2.802}  & 124.742        & 125.719   \\
\hline
\end{tabular}\end{center}
\end{table}
 \begin{table}[htpb]
\caption{ \ The first four eigenvalues using $P_{k}-P_{k-1}(k=3)$ element for the {\bf NIP} and {\bf IIP} methods on uniform meshes for the 3D thick L-shaped domain with $h=\sqrt{3}/8$, $\boldsymbol{\beta}=(0,0,1)^{t}$, and different stabilization parameters $\widehat{\gamma}$. }\label{tab-11}
\scriptsize
\label{tab:1}
\begin{center}
\begin{tabular}{{c|cccccccccccccccccc}}
  \hline
          &  $\widehat{\gamma}=1/32$    &  $\widehat{\gamma}=1/16$ & $\widehat{\gamma}=1/8$ & $\widehat{\gamma}=1/4$ & $\widehat{\gamma}=1/2$ & $\widehat{\gamma}=1$  & ~ \\ \hline
            &\fbox{$-$80.6797}    &  82.6610   & \fbox{$-$31.3178}  & 82.6397   & \fbox{19.7284}   & 82.6060   & ~ \\
   {\bf NIP} &82.6655   & 89.0239   & 82.6529          & 89.0195   & 82.6219          & 89.0095   & ~ \\
   $k=3$     &89.0248   & 124.8164  & 89.0223          & 124.7810  & 89.0151          & 124.7201   & ~ \\
             &124.8238   & 125.8341  & 124.8031         & 125.3276  & 124.7505         & 125.7878   & ~ \\
         \hline
             &83.4753   &82.9131  & 82.6634  & \fbox{41.6644}  & \fbox{34.9280}  & 82.5332   & ~ \\
    {\bf IIP}&90.2407   &89.5762  & 89.2733  & 82.5564  & 82.5232  & 89.0186   & ~ \\
     $k=3$  &128.1954    &126.1580  & 125.2764  & 89.1268  & 89.0543  & 124.6475   & ~ \\
            &129.0871    & 127.2696  & 126.4703  & 124.8741  & 124.7005  & 125.8068   & ~ \\
\hline
\end{tabular}\end{center}
\end{table}

\section*{Data availability}
Data will be made available on request.

\section*{Declaration of competing interest}
The authors declare that they have no known competing financial interests or personal relationships that could have appeared to influence the work reported in this paper.

\section*{Acknowledgments}
The authors cordially thank the editor and the referees for their
valuable comments and suggestions  which lead to the improvement of
this paper.
This work was supported by the High-level Talent Scientific Research Startup Fund(No. 052[2025] Document J of the School-Doctor Cooperation) of Guizhou Medical University and the National Natural Science Foundation of China (Grant Nos. 12261024, 11561014).



\begin{thebibliography}{00}


\bibitem{Benzi2007} M. Benzi, J. Liu, An efficient solver for the incompressible Navier-Stokes equations in rotation form, SIAM J. Sci. Comput. 29 (5) (2007) 1959-1981, https://dx.doi.org/10.1137/060658825.


 \bibitem{Fan2016} H.-T. Fan, X.-Y. Zhu, A modified relaxed splitting preconditioner for generalized saddle point problems from the incompressible Navier-Stokes equations, Appl. Math. Lett. 55 (2016) 18-26, https://dx.doi.org/10.1016/j.aml.2015.11.011.


\bibitem{Lepe2024}F. Lepe, G. Rivera, J. Vellojin, Finite element analysis of the Oseen eigenvalue problem, Comput. Methods Appl. Mech. Engrg. 425 (2024) 116959, https://dx.doi.org/10.1016/j.cma.2024.116959


\bibitem{Adak2025}D. Adak, F. Lepe, G. Rivera, A nonconforming virtual element approximation for the Oseen eigenvalue problem, IMA J. Numer. Anal., 00(2025) 1-37, https://doi.org/10.1093/imanum/drae108.

 \bibitem{Lepe2024new}F. Lepe, G. Rivera, J. Vellojin, A Mixed finite element method for the velocity-pseudostress formulation of the Oseen eigenvalue problem, SIAM. Sci. Comput., 47(6)(2025),A3356-A3382, https://doi.org/10.1137/24M1696950.

 \bibitem{ReedandHill1973}W. Reed and T. Hill, Triangular Mesh Methods for the Neutron Transport Equation, Technical Report LA-UR-73-479, Los Alamos Scientifik Laboratory, 1973.

 \bibitem{Antonio2012} D. A. Di Pietro, A. Ern, Mathematical Aspects of Discntinuous Galerkin Methods, Springer-Verlag, 2012, https://doi.org/10.1007/978-3-642-22980-0.

\bibitem{Riviere2008}B. Rivi\`{e}re,  Discontinuous Galerkin Methods for Solving Elliptic and Parabolic Equations, Theory and Implementation.  SIAM., 2008, https://doi.org/10.1137/1.9780898717440.

\bibitem{Hesthaven2008}J. S. Hesthaven, T. Warburton, Nodal Discontinuous Galerkin Methods: Algorithms, Analysis, and Applications. Springer, 2008, https://doi.org/10.1007/978-0-387-72067-8.

 \bibitem{Cockburn1999} B. Cockburn, G. E. Karniadakis, C. W. Shu, Discontinuous Galerkin methods: theory, computation
and applications, Berlin, Heidelberg: Springer, 2000, https://doi.org/10.1007/978-3-642-59721-3.

\bibitem{Hansbo2002}P. Hansbo, M. G. Larson, Discontinuous Galerkin methods for incompressible and nearly incompressible elasticity by Nitsche's method,
 Comput. Methods Appl. Mech. Engrg. 191 (2002) 1895-1908.

\bibitem{Stokes-flow2010} N. Nguyen, J. Peraire, B. Cockburn, A hybridizable discontinuous Galerkin method for Stokes flow, Comput. Methods Appl. Mech. Engrg., 199(15)(2010) 582-597, https://doi.org/10.1016/j.cma.2009.10.007.

\bibitem{Navier-Stokes2010}K. Chrysafinos, N. Walkington, Discontinuous Galerkin approximations of the Stokes and Navier-Stokes equations, Math.Comp., 79(272) (2010) 2135-2167, https://doi.org/10.1090/S0025-5718-10-02348-3.

\bibitem{Badia2014}S. Badia, R. Codina, T. Gudi, J. Guzm\'{a}n, Error analysis of discontinuous Galerkin methods
for the Stokes problem under minimal regularity, IMA. J. Numer. Anal., 34 (2014) 800-819, https://doi.org/10.1093/imanum/drt022.

\bibitem{Pietro2010}D. A. Di Pietro, A. Ern, Discrete functional analysis tools for discontinuous Galerkin methods with
application to the incompressible Navier-Stokes equations, Math. Comp., 79 (2010) 1303-1330, https://doi.org/10.1090/s0025-5718-10-02333-1.

\bibitem{Houston2005}P. Houston, D. Sch\"{o}tzau, T. P. Wihler, Energy norm a posteriori error estimation for mixed
discontinuous Galerkin approximations of the Stokes problem, J. Sci. Comput., 22 (2005) 347-370, https://doi.org/10.1007/s10915-004-4143-7.


\bibitem{Schotzau}D. Sch\"{o}tzau, C. Schwab, A. Toselli, Mixed hp-DGFEM for incompressible flows, SIAM J. Numer.
Anal., 40 (2002), 2171-2194, https://doi.org/10.1137/s0036142901399124.


 \bibitem{Kanschat2008}G. Kanschat, D. Sch\"{o}tzau, Energy norm a posteriori error estimation for divergence-free discontinuous Galerkin approximations of the Navier-Stokes equations.  Int. J. Numer. Meth. in Fluids, 57 (2010) 1093-1113, https://doi.org/10.1002/fld.1795.

\bibitem{Cockburn2002}B. Cockburn, G. Kanschat, D. Sch\"{o}tzau, C. Schwab, Local discontinuous Galerkin
methods for the Stokes system, SIAM J. Numer. Anal., 40 (2002) 319-343, https://doi.org/10.1137/S0036142900380121.

\bibitem{Antonietti2006}P. F. Antonietti, A. Buffa, I. Perugia, Discontinuous Galerkin approximation of the Laplace eigenproblem. Comput. Methods Appl. Mech. Eng. 195, 3483-3503 (2006), https://doi.org/10.1016/j.cma.2005.06.02330.

\bibitem{Descloux1978a} J. Descloux, N. Nassif, and J. Rappaz, On spectral approximation. I. The problem of convergence,
RAIRO Anal. Numer., 12 (1978), 97-112, iii, https://doi.org/10.1051/m2an/1978120200971.

\bibitem{Descloux1978b} J. Descloux, N. Nassif, and J. Rappaz, On spectral approximation. II. Error estimates for the
Galerkin method, RAIRO Anal. Numer., 12 (1978), 113-119, iii, https://doi.org/10.1051/m2an/
1978120201131.


\bibitem{Cliff2010}K. Cliffe, E. Hall, P. Houston, Adaptive discontinuous galerkin methods for eigenvalue problems arising in incompressible fluid flows, SIAM J. Sci. Comput. 31(6) (2010) 4607-4632, https://doi.org/10.1137/080731918.


 \bibitem{Sun2023Cicp}L. Sun, H. Bi, Y. Yang, A Multigrid Discretization of Discontinuous Galerkin Method for the Stokes Eigenvalue Problem. Commun. Comput. Phys., 34(5) (2023) 1391-1419, https://doi.org/10.4208/cicp.OA-2023-0027.



\bibitem{DG+Gedicke2020}J. Gedicke, A. Khan,  Divergence-conforming discontinuous Galerkin finite elements for Stokes eigenvalue problems. Numer. Math. 144 (2020)585-614, https://doi.org/10.1007/s00211-019-01095-x32.

\bibitem{DG+Lepe2020} F. Lepe, D. Mora, Symmetric and nonsymmetric discontinuous Galerkin methods for a pseudostress
formulation of the Stokes spectral problem. SIAM J. Sci. Comput. 42 (2020) 698-722, https://doi.org/10.1137/19M1259535.


\bibitem{Lepe2023+DG+Stokes} F. Lepe, Interior penalty discontinuous Galerkin methods for the velocity-pressure formulation of the Stokes spectral problem, Adv. Comput. Math. 49 (2023), no. 4, Paper No. 61, 31 pp. https://link. springer.com/article/10.1007/s10444-
023-10062-y.


 \bibitem{Arnold1982+SIP}D. N. Arnold, An interior penalty finite element method with discontinuous elements, SIAM J. Numer. Anal. 19 (1982) 742-760.

\bibitem{Baumann+NIP}C.E. Baumann, J.T. Oden, A discontinuous $hp$ finite element method for convection-diffusion problems, Comput. Methods Appl.
Mech. Engrg. 175 (3-4) (1999) 311-341.

\bibitem{Riviere+NIP}B. Rivi\`{e}re, M. F. Wheeler, V. Girault, Improved energy estimates for interior penalty, constrained and discontinuous Galerkin
methods for elliptic problems, Part I, Comput. Geosci. 3 (1999) 337-360.

\bibitem{Dawson+IIP}C. Dawson, S. Sun, M. F. Wheeler, Compatible algorithms for coupled flow and transport, Comput. Methods Appl. Mech. Engrg.
193 (2004) 2565-2580.

\bibitem{Boffi2010} D. Boffi, Finite element approximation of eigenvalue problems, Acta Numer, (2010) 1-120.

\bibitem{Arnold2002} D. N. Arnold, F. Brezzi, B. Cockburn, L. D. Marini, Unified analysis of
discontinuous Galerkin methods for elliptic problems,  SIAM J. Numer. Anal., 39(2002) 1749-1779.


\bibitem{Nedelec1986} J. C. N\'{e}d\'{e}lec, A new family of mixed finite elements in $\mathbb{R}^{3}$, Numer. Math., 50(1986) 57-81, https://doi.org/10.1007/BF01389668.


\bibitem{Acosta2006}G. Acosta, R. Dur\'{a}n,  M. A. Muschietti, Solutions of the divergence operator on John domains. Adv.
Math., 206(2006) 373-401, https://doi.org/10.1016/j.aim.2005.09.004.


\bibitem{Babuska1991book}I. Babu\v{s}ka, J. Osborn, Eigenvalue problems, in: P.G. Ciarlet, J.L. Lions (Eds.), Finite Element Methods (Part I), in: Handbook of Numerical Analysis, North-Holland: Elsevier Science Publishers, (1991) 641-787.

\bibitem{Babuska1978}I. Babu\v{s}ka, W. C. Rheinboldt,  Error estimates for adaptive finite element computations. SIAM J. Numer. Anal. 15 (1978) 736-754, https://doi.org/10.1137/0715049.

\bibitem{Verfurth2013}R. Verf\"{u}rth,  A Posteriori Error Estimation Techniques for Finite Eleement Methods, Numer. Math. Sci. Comput., Oxford University Press,  Oxford, 2013, https://doi.org/10.1093/acprof:oso/9780199679423.001.0001.

 \bibitem{oden2011}M. Ainsworth, J. Oden, A Posteriori Error Estimation in Finite Element Analysis. Wiley-Interscience, New York, 2000,   https://doi.org/10.1002/9781118032824.

\bibitem{Adaptive2023}L. Chamoin, F. Legoll, An Introductory Review on A Posteriori Error Estimation in Finite Element Computations, SIAM Review. 65(4) 2023, https://doi.org/10.1137/21M1464841.

\bibitem{Adaptive2024}A. Bonito, C. Canuto, R. H. Nochetto, A. Veeser, Adaptive finite element methods. Acta Numerica. 33 (2024)163-485, https://doi.org/10.1017/S0962492924000011.

\bibitem{Boffi2013}D. Boffi, F. Brezzi, M. Fortin, Mixed finite element methods and applications, Springer-Verlag Berlin Heidelberg, 2013, https://doi.org/10.1007/978-3-642-36519-5.


\bibitem{Boffi2025} D. Boffi, A. Khan,  Adaptive mixed FEM for the Stokes eigenvalue problem. Mathematics of Computation. (2025)[Online first]. https://doi.org/10.1090/mcom/4151.


 \bibitem{Brenner2003}S. Brenner, Poincar\'{e}-Friedrichs inequalities for piecewise $H^{1}$ functions,  SIAM J. Numer. Anal., 41 (2003) 306-324, https://doi.org/10.1081/NFA-200042165.

\bibitem{Ohannes2003}O. Karakashian, F. Pascal, A posteriori error estimates for a discontinuous Galerkin approximation of second-order elliptic problems, SIAM J. Numer. Anal.,41 (2004) 2374-2399, https://doi.org/10.1137/s0036142902405217.

\bibitem{John2016}V. John, Finite element methods for incompressible flow problems, in: Springer Series in Computational Mathematics, vol. 51, Springer, Cham, 2016, p.xiii+812, https://dx.doi.org/10.1007/978-3-319-45750-5.

\bibitem{Ern2021}A. Ern, J. Guermond, Finite elements II: Galerkin approximation, elliptic and mixed PDEs, Cham:
Springer, 2021, https://doi.org/10.1007/978-3-030-56923-5.

\bibitem{Brezzi1991}F. Brezzi, M. Fortin, Mixed and Hybrid Finite Element Methods. Springer-Verlag, New York, 1991, https://doi.org/10.1007/978-1-4612-3172-1.

\bibitem{Fabes1988} E.B. Fabes, C.E. Kenig, G.C. Verchota, The Dirichlet problem for the Stokes system on Lipschitz domains, Duke Math. J. 57 (3) (1988) 769-793,http://dx.doi.org/10.1215/S0012-7094-88-05734-1.

 \bibitem{Savare1998} G. Savar\'{e}, Regularity results for elliptic equations in Lipschitz domains, J. Funct. Anal. 152 (1) (1998) 176-201, http://dx.doi.org/10.1006/jfan.1997.3158.

 \bibitem{Kellogg1976}R. B. Kellogg, J. E. Osborn, A regularity result for the Stokes problem in a convex polygon, Journal of Functional Analysis, 21 (4)
(1976), 397-431, http://dx.doi.org/10.1016/0022-1236(76)90035-5.


\bibitem{Wang2024}S. Wang, H. Bi, Y. Yang, The a posteriori error estimates and an adaptive algorithm of the discontinuous Galerkin method for the modified transmission eigenvalue problem with absorbing media, Numer Algor, (2024), https://doi.org/10.1007/s11075-024-01981-y.

 \bibitem{ScottZhang1990}L. Scott, S. Zhang, Finite element interpolation of non-smooth functions satisfying boundary conditions,  Math. Comp., 54 (1990),  483-493, https://doi.org/10.1090/S0025-5718-1990-1011446-7.

\bibitem{Dai2014} X. Dai, L. He, A. Zhou, Convergence and quasi-optimal complexity of adaptive finite
element computations for multiple eigenvalues, IMA J. Numer. Anal., 35 (2014) 1934-1977, https://doi.org/10.1093/imanum/dru059.

\bibitem{Gallistl2015} D. Gallistl, An optimal adaptive FEM for eigenvalue clusters, Numer. Math., 130 (2015) 467-496, https://doi.org/10.1007/s00211-014-0671-8.

\bibitem{Boffi2017} D. Boffi, D. Gallistl, F. Gardini, L. Gastaldi, Optimal convergence of adaptive FEM for eigenvalue clusters in mixed form, Math. Comp., 86 (2017) 2213-2237.

\bibitem{Cances2020} E. Canc\`{e}s, G. Dusson, Y. Maday, B. Stamm, M. Vohral\'{\i}k, Guaranteed a posteriori bounds for eigenvalues and eigenvectors: Multiplicities and clusters, Math. Comp., 89 (2020) 2563-2611, https://doi.org/10.1090/MCOM/3549.

 \bibitem{yang2024} Y. Yang, S. Wang, H. Bi, The a posteriori error estimates of the FE approximation of defecitive eigenvalues for non-self-adjoint eigenvalue problems, SIAM J Numer Anal., 62(6) (2024) 2419-2438, https://doi.org/10.1137/23M162065X.

 \bibitem{Dorfler1996} W. D\"{o}rfler, A convergent adaptive algorithm for Poisson's equation,  SIAM J. Numer. Anal., 33(1996) 1106-1124, https://doi.org/10.1137/0733054.

 \bibitem{Chen2009}L. Chen, $i$FEM: An integrated finite element method package in Matlab, Technical Report, University of California at Irvine, (2009).





\end{thebibliography}
\end{document}